\documentclass[a4paper,oneside,fleqn,11pt]{article}
\usepackage{amsmath,amssymb,amsfonts,amsthm,type1cm,bm,xcolor,mathrsfs}
\usepackage{mathtools} 
\usepackage{graphicx,psfrag,epsf}
\usepackage{rotating}
\usepackage{authblk}
\usepackage{yfonts}
\usepackage{pdflscape}
\usepackage{setspace}
\usepackage{natbib}
\usepackage{float} 
\usepackage{caption}
\usepackage{subcaption}
\RequirePackage[colorlinks,citecolor=blue,urlcolor=pink]{hyperref}
\usepackage{bigstrut} 
\usepackage{comment}
\mathtoolsset{showonlyrefs=true}

\DeclareTextFontCommand{\texttt}{\ttfamily\upshape}
\DeclareMathOperator*\argmin{\arg\min}
\DeclareMathOperator\diag{diag}

\DeclareMathOperator\E{\mathbb{E}}

\DeclareMathOperator\F{F}

\DeclareMathOperator\FDR{FDR}

\DeclareMathAlphabet\mathbfcal{OMS}{cmsy}{b}{n}
\def\b0{\mathbf{0}}
\def\b1{\mathbf{1}}

\def\bA{\mathbf{A}}
\def\bB{\mathbf{B}}
\def\bC{\mathbf{C}}
\def\bD{\mathbf{D}}
\def\bE{\mathbf{E}}
\def\bF{\mathbf{F}}
\def\bG{\mathbf{G}}
\def\bH{\mathbf{H}}
\def\bI{\mathbf{I}}

\def\bK{\mathbf{K}}

\def\bM{\mathbf{M}}
\def\bN{\mathbf{N}}
\def\bP{\mathbf{P}}
\def\bQ{\mathbf{Q}}
\def\bR{\mathbf{R}}
\def\bS{\mathbf{S}}

\def\bU{\mathbf{U}}
\def\bV{\mathbf{V}}
\def\bW{\mathbf{W}}
\def\bX{\mathbf{X}}
\def\bY{\mathbf{Y}}
\def\bZ{\mathbf{Z}}

\def\bb{\mathbf{b}}
\def\bc{\mathbf{c}}

\def\be{\mathbf{e}}
\def\bff{\mathbf{f}}
\def\bg{\mathbf{g}}

\def\bu{\mathbf{u}}
\def\bv{\mathbf{v}}
\def\bw{\mathbf{w}}

\def\bz{\mathbf{z}}

\def\cD{\mathcal{D}}

\def\bbR{\mathbb{R}}

\def\bGamma{\boldsymbol{\Gamma}}
\def\bgamma{\boldsymbol{\gamma}}

\def\bbeta{\boldsymbol{\beta}}

\def\bDelta{\boldsymbol{\Delta}}
\def\bdelta{\boldsymbol{\delta}}
\def\bLambda{\boldsymbol{\Lambda}}

\def\bOmega{\boldsymbol{\Omega}}
\def\bPhi{\boldsymbol{\Phi}}

\def\bSigma{\boldsymbol{\Sigma}}
\def\bGamma{\boldsymbol{\Gamma}}

\def\bepsilon{\boldsymbol{\epsilon}}

\def\bzero{\mathbf{0}}

\theoremstyle{plain}
\newtheorem{thm}{Theorem}
\newtheorem{lem}{Lemma}
\newtheorem{prop}{Proposition}

\theoremstyle{definition}
\newtheorem{ass}{Assumption}
\newtheorem{rem}{Remark}

\newcommand{\CD}{\stackrel{d}{\longrightarrow}}
\newcommand{\CP}{\stackrel{p}{\longrightarrow}}

\makeatletter
\def\section{\@startsection {section}{1}{\z@}{-3.5ex plus -1ex minus-.2ex}{2.3ex plus .2ex}{\large\bf}}
\makeatother

\makeatletter
\def\subsection{\@startsection {subsection}{1}{\z@}{-3.5ex plus -1ex minus-.2ex}{2.3ex plus .2ex}{\normalsize\bf}}
\makeatother

\graphicspath{ {./images/} }

\title{Revisiting Asymptotic Theory for Principal Component Estimators of Approximate Factor Models}

\author[$\,\!$]{\textsc{Peiyun Jiang}$^\dagger$}
\author[$\,\!$]{\textsc{Yoshimasa Uematsu}\thanks{Correspondence: Yoshimasa Uematsu, Department of Social Data Science, Hitotsubashi University, 2-1 Naka, Kunitachi, Tokyo 186-8601, Japan (E-mail: yoshimasa.uematsu@r.hit-u.ac.jp).}}
\author[$\,\!$]{\textsc{Takashi Yamagata}$^\ddagger$}
\affil[$\dagger$]{\textit{Graduate School of Management, Tokyo Metropolitan University}}
\affil[*]{\textit{Department of Social Data Science, Hitotsubashi University}}
\affil[$\ddagger$]{\textit{Department of Economics and Related Studies, University of York}}
\affil[$\ddagger$]{\textit{Institute of Social and Economic Research, Osaka University}}

\begin{document}

\maketitle

\begin{abstract}
It is well known that approximate factor models exhibit rotation indeterminacy. Principal component (PC) estimators are typically analyzed relative to a rotated factor-loading representation, but the commonly used rotation depends on the estimator itself, leaving unclear which fixed population parameters are estimated. We show that any starting representation of the common component can be mapped by a rotation matrix $\bH$, constructed without using the PC estimator or the idiosyncratic errors and unique up to column signs, to a population representation satisfying the same PC normalization as the estimator. Although $\bH$ depends on the starting representation, the resulting representation is invariant to that choice up to simultaneous column sign changes. We call the resulting factors and loadings the pseudo-true (PC-normalized) parameters. Under a general weak factor model allowing signal eigenvalues to diverge at possibly different rates, we establish consistency and asymptotic normality of the PC estimators for these fixed, estimator-independent targets, together with fixed-target asymptotic theory for factor-augmented regressions. The theory thereby justifies confidence intervals for PC-normalized factors, gives loading profiles a fixed-target interpretation, and identifies the population coefficients associated with estimated PC factors in factor-augmented regressions.
\end{abstract}

\noindent%

\section{Introduction}

High-dimensional factor models have become an increasingly important analytical tool for psychology, finance, economics, biology, and so on. They are quite useful in reducing the high dimensionality of data by the low-rank approximation. For example, \cite{CR1983} first introduced the approximate factor model to finance, and estimation and inferential methods are developed by \cite{ConnorKorajczyk1986,ConnorKorajczyk1993}, \cite{StockWatson2002JASA,StockWatson2002JBES}, \cite{BaiNg2002}, \cite{Bai2003}, \cite{FanEtAl2013}, among many others.

\subsection{Factor model and estimator}

Consider a high-dimensional data matrix $\bX\in \bbR^{T\times N}$ generated by the latent factor model:
\begin{align}\label{model:AFM*}
\bX=\bF^{\ast}\bB^{\ast\prime}+\bE,
\end{align}
where $\bF^*=(\bff^*_1, \dots, \bff^*_T)'=(\bF^*_1,\dots,\bF^*_r)$ with $\bff_t^* \in \bbR^r$ and $\bF_k^* \in \bbR^T$ is a matrix of latent factors, $\bB^*=(\bb^*_1,\dots,\bb^*_N)'=(\bB_1^*,\dots,\bB_r^*)$ with $\bb_i^*\in \bbR^r$ and $\bB_k^*\in \bbR^N$ is a matrix of factor loadings, $\bE\in\bbR^{T\times N}$ is an idiosyncratic error matrix with $\E[\bE]=\bzero$, and $N,T\to\infty$ while $r$ is fixed. 
Throughout the paper, we suppose that $\bB^*$ and $\bF^*$ are of full column rank and the largest eigenvalue of $\E[T^{-1}\bE'\bE]$ is uniformly bounded in $N$ and $T$, but we do not require any specific normalization or structure in $(\bF^*,\bB^*)$, such as diagonality of $\bB^{*\prime}\bB^*$ and/or $\bF^{*\prime}\bF^*$, or sparseness of $\bB^*$. 
We also assume that the $r$ largest eigenvalues of the signal part, $T^{-1}\bF^*\bB^{*\prime}\bB^*\bF^{*\prime}$, denoted by $\lambda_1\geq \dots \geq \lambda_r$, diverge as $N\to\infty$. For later use, set $\bLambda=\diag(\lambda_1,\dots,\lambda_r)$.

The factor model in \eqref{model:AFM*} has a rotational indeterminacy. Although the \textit{starting parameter} $(\bF^{*},\bB^{*})$ is arbitrary from the viewpoint of the common component, a particular parameterization may be endowed with a structural interpretation when additional economic or subject-matter information is available. Throughout the paper, we therefore treat $(\bF^{*},\bB^{*})$ as a generic starting parameter; when additional identifying restrictions are imposed, it may also be interpreted as the \textit{structural parameter}.

The \textit{principal component} (PC) estimator, $(\hat{\bF}, \hat{\bB})$, is defined as a minimizer of $\|\bX-\bF\bB'\|_{\F}^2$ subject to the $r^2$ restrictions, $T^{-1}\bF'\bF=\bI_r$ and $\bB'\bB\in\cD(r)$, where $\cD(r)$ is a set of all the $r\times r$ diagonal matrices with $r$ positive diagonal entries in non-increasing order. 
The constraint minimization problem reduces to the eigenvalue problem of $T^{-1}\bX\bX'$; the factor estimator $\hat\bF\in\bbR^{T\times r}$ is obtained as $\sqrt{T}$ times the $r$ eigenvectors associated with the $r$ largest eigenvalues of $T^{-1}\bX\bX'$, and the loading estimator $\hat\bB\in\bbR^{N\times r}$ is $\hat{\bB}=T^{-1}\bX'\hat{\bF}$. 
By construction, $T^{-1}\hat{\bF}'\hat{\bF}=\bI_r$ and $\hat{\bB}'\hat{\bB}\in\cD(r)$ with its $k$th diagonal element equal to $\hat{\lambda}_k$, the $k$th largest eigenvalue of $T^{-1}\bX\bX'$. Setting $\hat{\bLambda}=\diag(\hat{\lambda}_1,\dots,\hat{\lambda}_r)$, we can write $\hat{\bB}'\hat{\bB}=\hat{\bLambda}$.

\subsection{Weak factor models with heterogeneous signal eigenvalues}

The majority of the literature on latent factor models, including \cite{StockWatson2002JASA,StockWatson2002JBES}, \cite{Bai2003}, \cite{BaiNg2002,BaiNg2006,BaiNg2013}, has employed the \textit{strong factor} (SF) models, where the signal eigenvalues satisfy $\lambda_k\asymp N$ for all $k=1,\dots,r$. This requirement is somewhat strong in view of real data; see evidence in \cite{DeMol2008}, \cite{BaileyEtAl2016}, \cite{WangFan2017}, \cite{Freyaldenhoven21JoE}, \cite{UY2019,UY2019inference}, among others. Therefore, in the present study, the new asymptotic properties of the PC estimator are derived for the \textit{weak factor} (WF) models that have possibly different divergence rates for the signal eigenvalues: $\lambda_k\asymp N^{\alpha_k}$ with $0<\alpha_k\leq 1$ for $k=1,\dots,r$. The SF model is a special case of the WF model with $\alpha_r=1$.

The WF models have been gaining increasing attention recently. 
\cite{UY2019,UY2019inference} and \cite{WeiZhang2026} consider the WF models induced by sparse factor loadings. \cite{Onatski2010} and \cite{Freyaldenhoven21JoE} propose methods to determine the number of factors for the WF models. 
\citet[Section 5]{BaiNg2023} consider a model with ``weaker loadings,'' essentially assuming that $\bB^{*\prime}\bB^*$ is asymptotically diagonal with the elements possibly diverging, arriving at a WF model characterized by the magnitude of signal eigenvalues. \cite{ChoiYuan2024} and \cite{FanEtAl2024} develop sharper asymptotic theory by applying leave-one-out and related perturbation techniques.

\subsection{Problem statement: What does the PC estimator estimate?}\label{sec:1.2}

The large body of influential literature, including \cite{Bai2003,BaiNg2002,BaiNg2006,BaiNg2023}, analyses the asymptotic properties of the PC estimator $(\hat{\bF}, \hat{\bB})$ relative to the ``rotated parameter'' $(\bF^*\hat{\bH},\bB^*\hat{\bH}^{\prime -1})$, where\footnote{The data-dependent rotation matrix $\hat\bH$ defined by \eqref{Hhat} is usually denoted without the `hat' in the literature; see \cite{Bai2003,BaiNg2023}, for example.}
\begin{align}\label{Hhat}
\hat{\bH} = {\bB^*}'\bB^*(T^{-1}{\bF^*}'\hat{\bF}) \hat{\bLambda}^{-1}. 
\end{align}
Typically, they discuss the asymptotic approximations,
\begin{align}\label{"consis"}
\hat{\bff}_t=\hat{\bH}'\bff_t^* +o_p(1),~~~~
\hat{\bb}_i=\hat{\bH}^{-1}\bb_i^* +o_p(1),
\end{align}
and their associated asymptotic normality. 
The issue is that the rotated parameter $(\hat{\bH}'\bff_t^*,\hat{\bH}^{-1}\bb_i^*)$ depends on the PC estimator itself through $\hat\bH$. Approximation \eqref{"consis"} literally says that each of $\hat{\bF}$ and $\hat{\bB}$, which is a function of $\bX$, approximates another function of $\bX$. This is not a conventional consistency statement for a fixed parameter and therefore does not reveal which fixed population parameters the PC estimators are estimating. 

Just to obtain a data-independent rotation matrix, it may seem sufficient to consider a probability limit $\bar{\bH}$ of $\hat{\bH}$, when such a limit is available. 
However, this route is not sufficient for our purpose in WF models, because the limiting behavior of $\hat{\bH}$ depends on the normalization and rate structure, especially when the signal eigenvalues diverge at different rates; see \citet[Section 5]{BaiNg2023}. 
Even if $\bar{\bH}$ is well-defined, an approach based on it is indirect; to approximate $(\hat{\bff}_t,\hat{\bb}_i)$ by $(\bar{\bH}'{\bff}_t^{*},\bar{\bH}^{-1}{\bb}_i^*)$, we first consider \eqref{"consis"} and then approximate $(\hat{\bH}'{\bff}_t^{*},\hat{\bH}^{-1}{\bb}_i^*)$ by $(\bar{\bH}'{\bff}_t^{*},\bar{\bH}^{-1}{\bb}_i^*)$. 
Another issue with using $\hat{\bH}$ is that it seems challenging to develop a comprehensive asymptotic theory for approximation \eqref{"consis"} in general WF models. 
In fact, Proposition 7-ii and the subsequent discussion in \cite{BaiNg2023} do not provide the asymptotic normality of $\sqrt{T}(\hat{\bb}_{i}-\hat{\bH}^{-1}\bb_{i}^{*})$ and $\hat{\bff}_t'\hat{\bb}_{i}-{\bff_t^*}'\bb_{i}^{*}$.

To address these issues, it is of great importance to find a rotation matrix, denoted by $\bH$, that depends solely on the chosen starting parameter $(\bF^{*},\bB^{*})$ and can be used to develop a natural and comprehensive asymptotic theory for the PC estimator with respect to $(\bH'\bff_t^{*},\bH^{-1}\bb_i^{*})$.

\subsection{Contributions}\label{sec:1 cont}

The first contribution of this article is to establish that, under very mild conditions on model \eqref{model:AFM*}, for any starting parameter $(\bF^{*},\bB^{*})$ of the common component, there exists a rotation matrix $\bH$, unique up to column signs, such that $(\bF^0,\bB^0):=(\bF^{*}\bH,\bB^{*}\bH^{\prime -1})$ satisfies \begin{align}\label{pc1}
\frac{1}{T}\bF^{0\prime}\bF^0 = \bI_r 
~~\text{ and }~~
\bB^{0\prime}\bB^0 \in \cD(r). 
\end{align}
We call the rotated parameter $(\bF^0,\bB^0)$ the \textit{pseudo-true} parameter. Although $\bH$ depends on the chosen starting parameter, the resulting pseudo-true parameter is invariant to that choice, up to simultaneous column sign changes. This rotation is useful because it automatically yields $\bF^0$ and $\bB^0$ satisfying the $r^2$ restrictions in \eqref{pc1}, which make them separately identifiable. This construction is not intended to recover a particular structural parameter; rather, the pseudo-true parameter is a fixed population parameter satisfying the same normalization as the PC estimator. Thus, it is natural to develop the asymptotic theory of the PC estimator relative to the pseudo-true parameter.

As the second contribution, we establish the consistency and asymptotic normality of the PC estimator relative to the pseudo-true parameter in WF models, thereby providing a fixed-target asymptotic theory in contrast to the conventional theory formulated relative to rotated quantities that themselves depend on the PC estimator. Within this fixed-target framework, we also sharpen the rate and factor-strength conditions substantially. These results suggest that considering $\bH$ rather than $\hat{\bH}$ is indispensable for deriving a sharp fixed-target asymptotic theory.

\subsection{Relation to recent weak-factor PCA theory}

Our contribution is complementary to recent weak-factor PCA theory. \citet{BaiNg2023}, \citet{FanEtAl2024}, and \citet{ChoiYuan2024} develop asymptotic theory for PC estimators relative to data-dependent rotations or alignments. Such devices yield sharp error bounds and distributional approximations, but the resulting rotated quantities are not fixed population objects. By contrast, our rotation $\bH$ depends only on the starting parameter $(\bF^*,\bB^*)$ and yields the fixed PC-normalized parameter $(\bF^0,\bB^0)$.

\citet{FanEtAl2024} and \citet{ChoiYuan2024} obtain sharper rates and allow weaker factors under additional distributional, dependence, or row-wise loading-strength conditions. These results concern approximation after data-dependent alignment, rather than convergence to $(\bF^0,\bB^0)$. Since the rotations, normalizations, assumptions, and error metrics differ, the rates cannot be compared directly. Our objective is instead to identify and analyze the fixed population target of the PC estimator.

This difference also prevents a direct application of the proof techniques in \citet{FanEtAl2024}. Their row-wise perturbation arguments are developed around a sample-dependent Procrustes alignment, whereas fixed-target inference additionally requires controlling its discrepancy from $\bH$ under the componentwise CLT normalizations. Under our Bai--Ng-type expansion, some factor-strength restrictions arise already under data-dependent centering. Our refined rotation bounds ensure that replacing the data-dependent targets by $(\bF^0,\bB^0)$ causes no additional rate loss, but do not remove all restrictions of the underlying expansion. Thus, when $N\asymp T$, our current asymptotic-normality results require $\alpha_r>1/2$, whereas \citet{FanEtAl2024} reach weaker factors under their different assumptions. The two approaches should therefore be viewed as complementary rather than directly competing.

\subsection{Practical relevance}\label{sec:1 prac}

The practical relevance of our approach is not that the pseudo-true parameters recover a particular structural representation. Rather, they provide fixed population counterparts of the quantities routinely used in empirical PC analyses. For example, even when $(\bF^{*},\bB^{*})$ is endowed with a structural interpretation, a test of $H_0^{i,k}: b^0_{i,k}=0$ is not, in general, a test of the structural restriction $b^{*}_{i,k}=0$. It is a test of whether the $i$th variable loads on the $k$th PC-normalized population factor. Hence such tests can be used to characterize the $k$th PC-normalized factor through the variables that significantly load on it; see \cite{UY2019inference} and Section \ref{subsec:addapplication} of Supplementary Material 
for an empirical illustration.

This distinction is important for inference. Confidence intervals based on conventional PC theory are often written for rotated quantities such as the $k$th element of $\hat{\bH}'\bff_t^*$. Since this target depends on the estimator, its interpretation as a fixed population parameter is unclear. In contrast, our theory justifies confidence intervals for fixed pseudo-true quantities such as $f^0_{t,k}$ and $b^0_{i,k}$; see Section \ref{subsec:application} for an empirical illustration.

The pseudo-true is also useful for studying factor strength; our rotation gives the identity, $\lambda_k=\|\bB_k^0\|_2^2$. Thus the $k$th pseudo-true loading vector $\bB_k^0$ carries information about the strength of the $k$th factor. If sparsity is imposed, as in \citet{UY2019}, the number of statistically significant elements in $\bB_k^0$ can be used to estimate the strength exponent $\alpha_k$.

\subsection{Factor-augmented regression}

Our approach also applies to factor augmented regressions, considered in \cite{BaiNg2006}, \cite{StockWatson2002JASA,StockWatson2002JBES} and \cite{LudvigsonNg2009}, among others. They are widely used in situations with a large number of predictors. Assuming that certain unobserved common factors influence the movement of the predictors, the factors extracted from the predictors can be used to forecast a particular series. To illustrate, consider a simple factor augmented predictive regression model, $y_{t+1}=\bgamma^{*\prime}\bff_t^*+\epsilon_{t+1}$, where $\bgamma^{*}$ is the coefficient vector associated with the chosen starting representation, $y_t$ is a variable of interest and $\epsilon_t$ is an error term. As the predictive factor $\bff_t^*$ is not observable, it is usually replaced by the PC estimator $\hat{\bff}_t$ extracted from a larger set of predictors, $\bX$. \cite{BaiNg2006} employ the approximation \eqref{"consis"} so that $\bgamma^{*\prime}\bff_t^*=\bgamma_{\hat{\bH}}'\hat{\bff}_t+o_p(1)$, where $\bgamma_{\hat{\bH}}=\hat{\bH}^{-1}\bgamma^*$. This approximation is useful, but the coefficient $\bgamma_{\hat{\bH}}$ on the PC estimator $\hat{\bff}_t$ is itself estimator-dependent.

Our pseudo-true model gives a fixed population counterpart: $y_{t+1}=\bgamma^{0\prime}\bff_t^0+\epsilon_{t+1}$ with $\bgamma^0=\bH^{-1}\bgamma^*$. Testing $\gamma_k^0=0$ is not equivalent to testing whether the $k$th factor in a particular structural representation has zero predictive content. Instead, it is the population counterpart of the conventional $t$-test on the $k$th estimated PC factor. Testing which PC factors are significant in factor augmented regressions is important because the extracted PC factors are ordered by their importance in the covariability of predictors, which may not correspond to their forecast power for the particular series of interest; see further discussions in \cite{BaiNg2008,BaiNg2009} and \cite{ChengHansen2015}. We provide formal asymptotic theory for this interpretation, which is new to the literature; an empirical illustration is given in the companion paper \citet{jiang2024Mw}.

\subsection{Organization}
The rest of the paper is structured as follows. 
In Section \ref{sec:2} the rotation matrix $\bH$, which is purely a function of signals $(\bF^*, \bB^{*})$, is derived and its uniqueness is proved.  
In Section \ref{sec:3} the weak factor model is formally introduced and assumptions are made, then consistency of the PC factor is shown in Section \ref{sec:4} and its asymptotic normality is proved in Section \ref{sec:5}. Section \ref{sec:6} discusses the asymptotic properties of the estimators of factor augmented regressions. Section \ref{sec: MC} discusses the experimental design and summarizes the results. 
Section \ref{subsec:application} provides an empirical illustration, and
Section \ref{sec:con} contains some concluding remarks.   


\section{Derivation of the Rotation Matrix}\label{sec:2}

We derive the explicit form of rotation matrix $\bH$ that rotates \eqref{model:AFM*} to the identifiable pseudo-true model,
\begin{align}\label{model:AFM0}
\bX=\bF^{0}\bB^{0\prime}+\mathbf{E},
\end{align}
where $\bF^{0}$ and $\bB^{0}$ satisfy \eqref{pc1} under a quite general assumption. For this purpose, we consider the eigenvalue problem of the $r \times r$ matrix ${\bB^*}'\bB^*\left(T^{-1}{\bF^*}'\bF^*\right)$. Let $\bLambda$ and $\bP$ denote the $r\times r$ diagonal matrix containing the eigenvalues of ${\bB^*}'\bB^*\left(T^{-1}{\bF^*}'\bF^*\right)$ in descending order and the $r\times r$ matrix whose columns are composed of the corresponding normalized eigenvectors, respectively. Then, we may write
\begin{align}\label{eigen-decomp}
{\bB^*}'\bB^*\left(T^{-1}{\bF^*}'\bF^*\right) \bP = \bP\bLambda. 
\end{align}
To progress further, we impose minimal conditions on \eqref{eigen-decomp}.

\begin{ass}\normalfont\label{cond:eigen}
(i) The smallest eigenvalues of ${\bB^*}'\bB^*$ and $T^{-1}{\bF^*}'\bF^*$ are bounded away from zero; 
(ii) All the $r$ diagonal elements of $\bLambda$ are distinct. 
\end{ass}

Assumption \ref{cond:eigen} is very mild. 
Under Assumption \ref{cond:eigen}, all the diagonal elements in $\bLambda$ are positive,  bounded away from zero, and distinct. Moreover, Assumption \ref{cond:eigen}(ii) implies the linear independence of the column vectors in $\bP$, so that $\bP^{-1}$ is well-defined. 
Let $\bU=\bP^{-1}{\bB^*}'\bB^*{\bP^{-1}}'$ and 
$\bV=\bP'\left(T^{-1}{\bF^*}'\bF^*\right)\bP$. 
Then \eqref{eigen-decomp} immediately implies $\bU\bV=\bLambda$. 

\begin{lem}\label{lem:UV}
If Assumption \ref{cond:eigen} is true, then $\bU$ and $\bV$ are invertible diagonal matrices.
\end{lem}

\begin{thm}\label{thm:rotation}
If Assumption~\ref{cond:eigen} holds, then there exists a rotation matrix $\bH$ that maps \eqref{model:AFM*} to \eqref{model:AFM0} satisfying \eqref{pc1}. It is unique up to column sign changes and is explicitly given by
\begin{align}\label{H}
\bH := \bP\bV^{-1/2}.
\end{align}
In particular, we have 
$\bB^{0\prime}\bB^0=\bLambda$ in \eqref{pc1}. 
\end{thm}

We emphasize that $\bH$ depends only on the starting parameter $(\bF^*,\bB^*)$, whereas $\hat{\bH}$ defined in \eqref{Hhat} also depends on the PC estimators $(\hat{\bF},\hat{\bB})$. By \eqref{H} and Lemma~\ref{lem:UV}, $\bH$ can be written as
\begin{align}
\bH
=
\bB^{*\prime}\bB^*
(T^{-1}\bF^{*\prime}\bF^0)
\bLambda^{-1},
\end{align}
and may therefore be regarded as a population counterpart of $\hat{\bH}$.

Theorem~\ref{thm:rotation} shows that, under a mild assumption, the starting model \eqref{model:AFM*} can always be rotated to the pseudo-true model \eqref{model:AFM0} satisfying \eqref{pc1}. The two models are observationally equivalent, while the latter provides a fixed population representation satisfying the same normalization as the PC estimator.

\begin{thm}\label{thm:invariance}
Let $(\bF^{\ast},\bB^{\ast})$, $\bH$, and $(\bF^{0},\bB^{0})$ be as in Theorem~\ref{thm:rotation}. Consider any other parameter $(\bF^{\dagger},\bB^{\dagger})$ for which Assumption~\ref{cond:eigen} holds in place of $(\bF^{\ast},\bB^{\ast})$ and $\bF^{\dagger}\bB^{\dagger\prime}=\bF^{\ast}\bB^{\ast\prime}$. Denote by $\bH^{\dagger}$ a rotation matrix constructed from $(\bF^{\dagger},\bB^{\dagger})$ in the same manner as $\bH$ in Theorem~\ref{thm:rotation}. Then $\bF^{\dagger}\bH^{\dagger}=\bF^{0}$ and $\bB^{\dagger}\bH^{\dagger\prime-1}=\bB^{0}$ hold up to column sign changes.
\end{thm}

Theorem~\ref{thm:invariance} clarifies that the dependence of the rotation matrix on the starting parameter does not carry over to the pseudo-true parameter. Any two full-rank parameterizations of the same common component are related by a nonsingular change of coordinates, while the corresponding rotation matrices change inversely so as to offset this change exactly. Hence, up to column sign changes, $(\bF^{0},\bB^{0})$ is the unique representation of the common component satisfying \eqref{pc1}, regardless of the starting parameter.

\section{Weak Factor Models}\label{sec:3}

We formulate the weak factor (WF) models in preparation for deriving the asymptotic theory. We start with introducing the assumptions on the idiosyncratic error term.  

\begin{ass}[Idiosyncratic errors]\normalfont\label{assumption errors}
For some constant $M<\infty$ that does not depend on $N$ and $T$, we have:
(i) $\E[e_{t,i}|\bF^*,\bB^*]=0$ and
$\E[e_{t,i}^{4}| \bF^*,\bB^*]\le M$ for all $i$ and $t$; 
(ii) $\E[e_{s,i}e_{t,j}| \bF^*,\bB^*]=\E[e_{s,i}e_{t,j}]$ for all $s,t,i,j$;
(iii) for all $i$,
$\bigl|\E[e_{s,i}e_{t,i}]\bigr|\le|\gamma_{s,t}|$ for some $\gamma_{s,t}$ such that
$\sum_{t=1}^{T}|\gamma_{s,t}|\le M$;
(iv) for all $t$,
$\bigl|\E[e_{t,i}e_{t,j}]\bigr|\le|\tau_{i,j}|$ for some $\tau_{i,j}$ such that
$\sum_{j=1}^{N}|\tau_{i,j}|\le M$;
(v) $\dfrac{1}{NT}\displaystyle\sum_{t,s=1}^{T}\sum_{i,j=1}^{N}
\bigl|Cov\!\left(e_{t,i}e_{t,j},\,e_{s,i}e_{s,j}\right)\bigr|\le M$; 
(vi) $\|\bE\|_{2}^{2}=O_{p}\!\left(\max\{N,T\}\right)$;
(vii) the minimum eigenvalue of $\bSigma_{e}=\E[\be_{t}\be_{t}']$ is bounded away from zero.
\end{ass}

Assumption \ref{assumption errors}(i) collects the exogeneity and moment conditions. 
The conditioning on the entire $(\bF^*,\bB^*)$, rather than on the single pair 
$(\bff_t^*,\bb_i^*)$, is needed because the pseudo-true signals 
$(\bff_t^0,\bb_i^0)=(\bH'\bff_t^*,\bH^{-1}\bb_i^*)$ depend on the whole panel through 
$\bH$, so that conditioning on $(\bF^*,\bB^*)$ is what delivers $\E[\bff_t^0 e_{t,i}]=\mathbf{0}$ 
and $\E[\bb_i^0 e_{t,i}]=\mathbf{0}$.
Assumption \ref{assumption errors}(ii) adds that conditioning on the panel does not change 
the second moments of the errors, so the conditional variances of the rotated summands 
$\bff_t^0 e_{t,i}$ and $\bb_i^0 e_{t,i}$ reduce to unconditional ones.
The weak cross-sectional and serial correlations in Assumption 
\ref{assumption errors}(iii) and (iv) are similar to \citet[Assumption C]{Bai2003}; they 
imply $N^{-1}\sum_{i=1}^N\sum_{j=1}^N|\E(e_{t,i}e_{t,j})|\le M$ and 
$T^{-1}\sum_{s=1}^T\sum_{t=1}^T|\E(e_{s,i}e_{t,i})|\le M$, respectively. 
Assumption \ref{assumption errors}(v) is a mild fourth-order regularity condition limiting 
the dependence among the squared errors.
Assumption \ref{assumption errors}(i), (vi), 
and (vii) are also frequently imposed.

\begin{ass}[Signal strength]\normalfont\label{assumption factors}
There exist random or non-random variables $d_1,\dots,d_r>0$ and constants $0<\alpha_r\leq \dots \leq \alpha_1 \leq 1$ such that $\lambda_k=d_k^2N^{\alpha_k}$ for $k=1,\dots, r$ with ordered $0<\lambda_r< \dots < \lambda_1$ for large $N$. If $d_k$'s are random, we have $\E[d_k^4] \le M$ for all $k$ and $d_k\ge 1/M $ almost surely hold for some constant $M<\infty$. In addition, $d_1,\dots,d_r$ are well separated, in the sense that 
$\max_{k\ne\ell:\,\alpha_k=\alpha_\ell}\left|d_k^2/(d_k^2-d_\ell^2)\right|=O_p(1)$, 
so that the signal eigenvalues with the same $\alpha_k$ remain asymptotically distinct.
\end{ass}

Hereafter, denote $\bN=\diag(N^{\alpha_1}, \dots,N^{\alpha_r})$ and $\bD=\diag(d_1,\dots,d_r)$, so that we can write $\bLambda=\bD^2 \bN$. 
Under this condition, all the signal eigenvalues, $\lambda_1,\dots,\lambda_r$, are distinct and $\lambda_k\asymp N^{\alpha_k}$. If $\alpha_r=1$, this reduces to the SF model. Compared with Assumption \ref{assumption factors}, \citet[Section 5]{BaiNg2023} and \cite{Freyaldenhoven21JoE} impose the conditions that $\bN^{-1/2}\bB^{*\prime}\bB^*\bN^{-1/2}$ tends to a diagonal matrix. They additionally suppose that $T^{-1}\bF^{*\prime}\bF^*$ converges (in probability) to a positive definite matrix and $\bI_r$, respectively. These conditions lead to $\lambda_k \asymp N^{\alpha_k}$, and thus their models are examples of our WF models. However, their assumptions are much stronger than ours because they directly restrict the structure of $(\bF^*,\bB^*)$ in data generating process \eqref{model:AFM*}. They are unobserved and it does not seem possible to identify them without extra exogenous information. On the contrary, our strategy is to interpret the PC estimators $(\hat\bF,\hat\bB)$ as the consistent estimators of $\bF^0=\bF^* \bH$ and $\bB^0=\bB^* \bH^{\prime -1}$ in the pseudo-true model \eqref{model:AFM0} with \eqref{pc1}. This interpretation is always possible under Assumption \ref{cond:eigen} as considered in Section \ref{sec:2}. If necessary, we rotate back $\bF^0$ to identify $\bF^*$ when sufficient exogenous information exists for the purpose.

\section{Consistency}\label{sec:4}

We reconsider statistical consistency of the PC estimator. As noted in the introduction, the majority of the existing results have considered ``consistency'' for the rotated parameter by $\hat{\bH}$, which is a function of the PC estimator. In this section, we derive consistency for the pseudo-true parameters obtained by $\bH$ in \eqref{model:AFM0}. Accordingly, we will make further assumptions on this pseudo-true model for investigating the asymptotics. 

Write $\bB^*=(\bb^*_1,\dots,\bb^*_N)'=(\bB_1^*,\dots,\bB_r^*)$ and $\bF^*=(\bff^*_1, \dots, \bff^*_T)'=(\bF^*_1,\dots,\bF^*_r)$; the same notational rule applies to the other matrices. Since the starting parameter is arbitrary up to nonsingular rotations, imposing these moment conditions directly on $(\bF^*,\bB^*)$ would make them depend on the chosen coordinates. We therefore formulate them for the uniquely identified pseudo-true parameter.

\begin{ass}[Factors and Loadings]\normalfont\label{assumption BaiNgA3}
For some constant $M < \infty$ that does not depend on $N$ and $T$, we have: 
(i) $\E\|\bff_t^0\|_2^4 \le M$ and $\E\|\bb_i^0\|_2^4 \le M$; 
(ii) $\E\| \bN^{-\frac{1}{2}}\sum_{i=1}^{N} \bb^0_ie_{t,i} \|_2^2 \le M$ for each $t$; 
(iii) $\E\| T^{-\frac{1}{2}}\sum_{t=1}^{T}\bff_t^0e_{t,i}\|_2^2 \le M$ for each $i$; 
(iv)  $\E\| T^{-\frac{1}{2}}\bN^{-\frac{1}{2}}\sum_{t=1}^{T}\sum_{j=1}^{N}\bb_j^0[e_{t,i}e_{t,j}-\E(e_{t,i}e_{t,j})]\|_2^2 \le M$ for each $i$; 
(v) $\E\| (NT)^{-\frac{1}{2}}\sum_{s=1}^{T}\sum_{i=1}^{N}\bff_s^0[e_{s,i}e_{t,i}-\E(e_{s,i}e_{t,i})]\|_2^2 \le M$ for each $t$; 
(vi) the $r \times r$ matrix satisfies $\E\| T^{-\frac{1}{2}}\bN^{-\frac{1}{2}} \sum_{t=1}^{T}\sum_{i=1}^{N}\bb_i^{0}e_{t,i}\bff_t^{0'} \|_{2}^2 \le M$.
\end{ass}
We impose Assumption \ref{assumption BaiNgA3} at a high level, in the style of 
\citet{Bai2003} and \citet{BaiNg2023}. The difference is that we place these conditions on 
the pseudo-true parameter $(\bF^0,\bB^0)$ rather than on a starting representation, so that 
they do not depend on the chosen coordinates. Assumption \ref{assumption errors} ensures 
that the conditional moments of the rotated summands $\bb_i^0 e_{t,i}$ and $\bff_t^0 
e_{t,i}$ are well defined. Conditions (ii), (iii), and (vi) involve the rotated loadings 
and factors jointly with the errors, rather than $\bb_i^0$ and $e_{t,i}$ (or $\bff_t^0$ and 
$e_{t,i}$) separately, because under the weak factor scaling $\bN^{-1/2}$ the separate 
treatment is not sufficient for boundedness.

Following the analysis conducted by \cite{BaiNg2023}, we use other data-dependent rotations of $\bF^*$ and $\bB^*$ than $\hat\bH$, such as 
\begin{align*}
\hat{\bH}_4={\bB^{*}}'\hat{\bB}\hat{\bLambda}^{-1}, 
~
\tilde{\bH}_4={\bB^{0}}'\hat{\bB}\hat{\bLambda}^{-1},
~
\tilde{\bH}={\bB^0}'\bB^0(T^{-1}{\bF^0}'\hat{\bF}) \hat{\bLambda}^{-1}, 
~
\hat{\bQ}= T^{-1}\hat{\bF}'\bF^*, 
~
\tilde{\bQ}= T^{-1}\hat{\bF}'\bF^0, 
\end{align*}
and others appear in Section \ref{subsec:related lemmas} of Supplementary Material. These rotation matrices are theoretically important since the discrepancy between the PC estimator and the parameter rotated by such matrices can converge faster. By the definition, we immediately obtain $\bF^*\hat{\bH}_4=\bF^0\tilde{\bH}_4$, $\bB^*\hat{\bQ}'=\bB^0\tilde{\bQ}'$, and  $\bF^*\hat{\bH}=\bF^0\tilde{\bH}$. Because the PC eigenvectors are unique only up to column signs, we henceforth choose the signs of \((\hat \bF,\hat \bB)\) so that
$
T^{-1}\hat \bF_k^\prime \bF_k^0\ge 0,\; k=1,\ldots,r.
$ Without it, all convergence statements involving \(\tilde \bQ\) and \(\tilde \bH_4\) hold only up to a diagonal sign matrix. Employing these rotation matrices, we achieve the ``consistency'' results for the data-dependent counterparts:

\begin{lem}\label{lem fh4bq}
Suppose that Assumptions \ref{cond:eigen}--\ref{assumption BaiNgA3} hold. If $\frac{N^{1-\alpha_r}}{T}\rightarrow 0$, then we have
\begin{align*}
&(i)~\frac{1}{\sqrt{T}} \|  \hat{\bF} - \bF^* {\hat{\bH}_4}  \|_{\F} 
= \frac{1}{\sqrt{T}} \|  \hat{\bF} - \bF^0 {\tilde{\bH}_4}  \|_{\F}
=O_p\left(\frac{N^{1-\alpha_r}}{T}\right) 
+O_p\left(N^{-\frac{1}{2}\alpha_r}\right), \\
&(ii)~\frac{1}{\sqrt{N}}	 \| \hat{\bB} - \bB^* \hat{\bQ}' \|_{\F}
=\frac{1}{\sqrt{N}}	 \| \hat{\bB} - \bB^0 \tilde{\bQ}' \|_{\F}
=O_p\left(\frac{1}{\sqrt{T}}\right) +
O_p\left(N^{-\frac{1}{2}-\frac{1}{2}\alpha_{r}}\right),\\
 & (iii)~ \tilde{\bH}_4,\tilde{\bQ}\; \text{are} \;O_p(1). \;\text{Moreover,}\;
 \tilde{\bH}_4,\tilde{\bQ} \CP \bI_r.
 \end{align*}
\end{lem}

Lemma \ref{lem fh4bq}(i) is in line with the result in \citet[Proposition 6(i)]{BaiNg2023}. Lemma \ref{lem fh4bq}(ii) sharpens the corresponding results in \citet[Proposition 6(ii)]{BaiNg2023}, but related to \citet[Proposition 3.3]{WeiZhang2026}. Interestingly, Lemma \ref{lem fh4bq}(ii) does not depend on the value of $\alpha_1$, which may lead to the fast convergence rates. The boundedness of the rotation matrices with ``tilde'' implied by Lemma \ref{lem fh4bq}(iii) is key to deriving the asymptotic normality for the pseudo-true parameters in the next section. In contrast, the rotation matrices with ``hat'' are not necessarily bounded unless additional conditions are imposed on $\bF^*$. 

\begin{thm}\label{0} Suppose that Assumptions \ref{cond:eigen}--\ref{assumption BaiNgA3}
hold. If $\frac{N^{1-\alpha_r}}{T}\rightarrow 0$, then we have
\begin{align*}
&(i)~	\frac{1}{\sqrt{T}} \|  \hat{\bF} - \bF^0\|_{\F}  = 
 O_p\left(\frac{N^{1-\alpha_r}}{T}\right)+
  O_p\left(N^{-\frac{1}{2}\alpha_r}\right), \\
&(ii)~ \frac{1}{\sqrt{N}}	 \| \hat{\bB} - \bB^0 \|_{\F}
		= O_p\left(\frac{1}{\sqrt T}\right) + O_p\left(N^{-\frac{1}{2}-\frac{1}{2}\alpha_r}\right), \\
&(iii)~	\frac{1}{\sqrt{NT}} \|  \hat{\bC} - \bC^*\|_{\F}  =
	O_p\left(\frac{1}{\sqrt{T}}\right)+ O_p\left(\frac{1}{\sqrt{N}}\right),
\end{align*}
where $\hat{\bC} =\hat{\bF}\hat{\bB}^{\prime}$ and  $\bC^* = \bF^*\bB^{* \prime}= \bF^0\bB^{0 \prime}$. 
\end{thm} 

This consistency theorem is new as it involves the pseudo-true parameter $(\bF^0,\bB^0)=(\bF^*\bH, \bB^*{\bH}^{\prime-1})$ rather than the data-dependent parameter as in Lemma \ref{lem fh4bq} and \citet{BaiNg2023}. This means that Theorem \ref{0} achieves the \textit{consistency} of the PC estimator for the pseudo-true parameter. 
A notable implication of the sharp bounds for \( (\tilde{\bH}_4-\bI_r,\tilde{\bQ}-\bI_r)\) is that the pseudo-true parameter interpretations do not entail an additional rate loss for the factor and loading estimators. In particular, the convergence rates in Theorem \ref{0}(i)(ii) are the same as those in Lemma \ref{lem fh4bq}(i)(ii). This follows because the additional terms $(\bF^0(\tilde{\bH}_4-\bI_r), \bB^0(\tilde{\bQ}'-\bI_r))$ caused by replacing the data-dependent parameter $(\bF^0\tilde{\bH}_4, \bB^0\tilde{\bQ}')$ with the pseudo-true parameter \((\bF^0, \bB^0)\) are of the same order as the data-dependent approximation error \( ( \hat{\bF} - \bF^0 {\tilde{\bH}_4}, \hat{\bB}-\bB^0\tilde{\bQ}')\). 

In the case of SF models (i.e., $\alpha_r=1$) in Theorem \ref{0}, the rates in (i) and (ii) become $\frac{1}{T}+\frac{1}{\sqrt{N}}$ and $\frac{1}{\sqrt{T}}+\frac{1}{N}$, respectively. They are the same as those in Lemma \ref{lem fh4bq}. This is remarkable because Theorem \ref{0} reveals the consistency of the PC estimators in SF models, whereas Lemma \ref{lem fh4bq} or \cite{BaiNg2002} does not. Moreover, these rates are smaller than those in \cite{BaiNg2002} for the SF models, which are $\frac{1}{\sqrt{T}}+\frac{1}{\sqrt{N}}$.  

Regarding the common component estimator, $\hat{\bC}$, the rate given by Theorem \ref{0}(iii) is faster than that shown by \citet[Proposition 3.3]{WeiZhang2026}.

\section{Asymptotic Normality}\label{sec:5}

We derive the asymptotic normality of the PC estimator centered by the pseudo-true parameter. The significance of this normality lies in justification of testing for restrictions on the pseudo-true parameters that are not data-dependent. We emphasize that data-dependent centering does not provide conventional inference for componentwise fixed population restrictions, though rotation-invariant restrictions, such as joint zero restrictions, are an exception.
To facilitate the discussion, we suppose a central limit theorem.

\begin{ass}[CLT]\normalfont
	\label{assumption dist}
For each $i$ and $t$, we have $\bD^{-1}\bN^{-\frac{1}{2}} \sum_{i=1}^{N} \bb_i^0e_{t,i}\CD N(\mathbf{0}, \bGamma_t)$\, and $ \frac{1 }{\sqrt{T}}\sum_{t=1}^{T}\bff_t^0e_{t,i}\CD N(\mathbf{0}, \bPhi_i)$, where $\bGamma_t=\lim_{N\rightarrow\infty} \sum_{i=1}^{N}\sum_{j=1}^{N}\E[\bD^{-1}\bN^{-1/2}\bb_i^0 e_{t,i}e_{t,j}\bb_j^{0 \prime}\bN^{-1/2}\bD^{-1}]$ and $
\bPhi_i=\lim_{T\rightarrow\infty} \frac{1}{T}\sum_{t=1}^{T}\sum_{s=1}^{T}\E[\bff_t^0 e_{t,i}e_{s,i}\bff_s^{0 \prime}].
$
They are bounded and positive definite.
\end{ass}

We first derive the asymptotic normality of the PC estimator with respect to the rotated parameters by data-dependent rotation matrices. 

\begin{lem}\label{lem dist} Suppose that Assumptions \ref{cond:eigen}--\ref{assumption dist} hold. \\
(i) If $\alpha_r>1/2$ and $\frac{N^{\frac{3}{2}-\alpha_{r}}}{T}  \to 0$, we have
$\bD\bN^{\frac{1}{2}}(\hat{\bff}_t-\hat{\bH}_4'\bff_t^*) \CD N(\mathbf{0}, \bGamma_t)$. \\
(ii) If $\frac{N^{1-\alpha_{r}}}{\sqrt{T}} \to 0$ and $\frac{\sqrt{T}}{N^{\alpha_r}} \to 0$, we have
$\sqrt{T}(\hat{\bb}_i-\hat{\bQ}\bb_i^*) \CD N(\mathbf{0}, \bPhi_i)$. 
\end{lem}

Lemma \ref{lem dist}(i) is similar to \citet[Proposition 7(i)]{BaiNg2023}. Lemma \ref{lem dist}(ii) is new and theoretically useful in comparison to the existing results, such as \citet[Proposition 7(ii)]{BaiNg2023} and \citet[Theorem 3.4(ii)]{WeiZhang2026}.
In the case where $N\asymp T$, the conditions for Lemma \ref{lem dist}(i) and (ii) are simplified to $\alpha_r>1/2$. The lemma leads to the next theorem:

\begin{thm}\label{thm dist} Suppose that Assumptions \ref{cond:eigen}--\ref{assumption dist} hold. \\
(i) If $\alpha_r>1/2$, and $\frac{N^{\frac{3}{2}-\alpha_{r}}}{T}  \to 0$, we have 
$\bD\bN^{\frac{1}{2}}(\hat{\bff}_t-\bff_t^0)\CD N(\mathbf{0}, \bGamma_t)$. \\
(ii) If $\frac{N^{1-\alpha_{r}}}{\sqrt{T}} \to 0$, and $\frac{\sqrt{T}}{N^{\alpha_r}}\to 0$, we have
$\sqrt{T}(\hat{\bb}_i-\bb_i^0) \CD N(\mathbf{0}, \bPhi_i)$. \\
(iii) If $\alpha_r>1/2$, $\frac{N^{\frac{3}{2}-\alpha_{r}}}{T}  \to 0$, and $\frac{\sqrt{T}}{N^{\alpha_r}} \to 0$, we have $\frac{\hat{c}_{t,i}-c_{t,i}^*}{\sqrt{V_{t,i}+U_{t,i}}}
\CD N(0, 1)$, where $V_{t,i}={\bb_i^0}'\bD^{-1}\bN^{-\frac{1}{2}}\bGamma_t\bD^{-1}\bN^{-\frac{1}{2}}\bb_i^0$ and $U_{t,i}=T^{-1}{\bff_t^0}'\bPhi_i\bff_t^0$.
\end{thm}

Theorem \ref{thm dist} allows us to construct confidence intervals for the pseudo-true factors and loadings. On the other hand, Lemma \ref{lem dist}, \citet{Bai2003}, and  \cite{WeiZhang2026}, among others are not generally applicable as the rotated parameters depend on the estimator itself; the only exception is the test for rotation-invariant restrictions, such as joint zero restrictions. 

In comparison to Lemma \ref{lem dist}(i)(ii), Theorem \ref{thm dist}(i)(ii) requires no additional conditions. The conditions for Theorem \ref{thm dist}(iii) are weaker than imposing the conditions for (i) and (ii) jointly because the estimation of $\bC^*$ does not require separate identification of the factors and factor loadings, and Lemma \ref{lem dist} applies directly. When $N\asymp T$, the conditions for Theorem \ref{thm dist}(i)(ii)(iii) are simplified to $\alpha_r>1/2$. When the SF model (i.e., $\alpha_r=1$) is considered, the conditions for Theorem \ref{thm dist}(i) and (ii) reduce to ${\sqrt{N}}/{T}\to 0$ and ${\sqrt{T}}/{N}\to 0$, respectively, which are identical to those in \citet[Theorems 1\&2]{Bai2003} for the case with the data-dependent rotation. 

To conduct statistical inference  based on Theorem \ref{thm dist}, the unknowns should be replaced with estimators. In Theorem \ref{thm dist}(i), the scaling matrix $\bD\bN^{1/2}$  can be consistently estimated by $(\hat{\bB}^{\prime}\hat{\bB})^{1/2}$ because Assumption \ref{assumption factors} gives $(\bB^{0\prime}\bB^0)^{1/2}=\bD \bN^{1/2}$. 
As for $\bGamma_t$ and $\bPhi_i$, the choice of the estimators depends on the dependence structure of $e_{t,i}$. 
In general, we may use $\hat{\bGamma}_t=(\hat{\bB}^{\prime}\hat{\bB})^{-1/2}
\hat{\bOmega}_{t} 
(\hat{\bB}^{\prime}\hat{\bB})^{-1/2}$ with $\hat{\bOmega}_{t}$ appropriately chosen. For example, when cross-sectional independence and time-series heteroskedasticity are allowed, choose $\hat{\bOmega}_{t}=
\sum_{i=1}^{N}\hat{\bb}_i \hat{e}_{t,i}^2\hat{\bb}_i^{\prime}$, where $\hat\bE=\bX-\hat\bF \hat\bB'=(\hat{e}_{t,i})$. For $\bPhi_i$, when $e_{t,i}$ is serially correlated and heteroskedastic over $i$, the heteroskedasticity and autocorrelation consistent (HAC) estimator of \cite{NeweyWest1987} can be employed: 
$\hat{\bPhi}_i=\hat{\bV}_{0,i} + \sum_{\ell=1}^L(1-\ell/(L+1))(\hat{\bV}_{\ell,i}+\hat{\bV}_{\ell,i}^{\prime})$, where $L>0$ is a slowly diverging sequence and 
$\hat{\bV}_{\ell,i}=T^{-1}\sum_{t=\ell+1}^T \hat{\bff}_t \hat{e}_{t,i}\hat{e}_{t-\ell,i} 
\hat{\bff}_{t-\ell}^{\prime}$. 
Finally, for Theorem \ref{thm dist}(iii), we may construct $\hat{V}_{t,i}={\hat\bb_i}'(\hat{\bB}^{\prime}\hat{\bB})^{-1/2}\hat{\bGamma}_t(\hat{\bB}^{\prime}\hat{\bB})^{-1/2}\hat{\bb}_i$ and $\hat{U}_{t,i}=T^{-1}{\hat{\bff}_t}'\hat{\bPhi}_i\hat{\bff}_t$.

\section{Factor-Augmented Regression}\label{sec:6}

We consider factor augmented regression models. Assuming that certain unobserved common factors influence the predictors, the factors extracted from the predictors can be used to forecast a particular series; see \cite{StockWatson2002JASA}. In such an analysis, significance tests for individual factors are important because the estimated factors ordered by the magnitude of the associated eigenvalues may not correspond to their forecast power for a particular series of interest; see further discussions in \cite{BaiNg2008,BaiNg2009} and \cite{ChengHansen2015}.

Consider the factor augmented regression model
\begin{align}\label{augmodel*}
y_{t+h}={\bgamma^*}'\bff_t^*+\bbeta' \bw_t+\epsilon_{t+h},
\end{align}
where $\bff_t^*$ is a vector of latent factors, $\bw_t$ is a vector of exogenous variables, $\epsilon_t$ is an error term, and $({\bgamma^*}',\bbeta')'$ is a coefficient vector. Since $\bff_t^*$ is unobservable, it is replaced by the PC estimator $\hat{\bff}_t$ estimated from a set of many predictors, $\bX$. 
Following \cite{BaiNg2006,BaiNg2023}, we may choose the data-dependent rotation matrix $\hat\bH$. Then, under approximation \eqref{"consis"} and some boundedness condition on $\bgamma_{\hat\bH} := \hat\bH^{-1}\bgamma^*$, we have
$
{\bgamma^*}'\bff_t^*
= \bgamma_{\hat\bH}' \hat{\bff}_t - \bgamma_{\hat\bH}' (\hat{\bff}_t-\hat\bH'\bff_t^*)
= \bgamma_{\hat\bH}' \hat{\bff}_t + o_p(1).
$
Since $\hat{\bff}_t$ is observed, we can estimate its coefficient vector, and denote it as $\hat\bgamma$. We may derive the asymptotic normality for the discrepancy, $\hat\bgamma-\bgamma_{\hat\bH}$, and construct the associated $t$-ratio. However, it is questionable whether using such statistics is justified as the ``parameter'' $\bgamma_{\hat\bH}$ depends on the PC estimator through $\hat\bH$. Despite the problem, such tests are routinely reported in empirical studies; see \cite{LudvigsonNg2009} among many others.

As in the previous section, we consider $\hat{\bff}_t$ as the estimator of ${\bff_t^0}$. To facilitate this perspective, we introduce $\bH$ and rewrite model \eqref{augmodel*} to the pseudo-true model:
\begin{align}\label{augmodel0}
y_{t+h}=\bgamma^{0 \prime}\bff_t^0+\bbeta' \bw_t+\epsilon_{t+h},
\end{align}
where $\bgamma^0:=\bH^{-1}\bgamma^*$ and $\bff_t^0=\bH^{\prime}\bff_t^*$. Observe that $\bgamma^0$ is a function of the parameters, $\bF^*$, $\bB^*$, and $\bgamma^*$. Therefore, tests for general parameter restrictions on $\bgamma^0$ seem justified. Rewrite \eqref{augmodel0} to the matrix form: 
\begin{align}\label{augmodel0_mat}
    \bY={\bF}^0 \bgamma^0 + \bW \bbeta + \bepsilon= \bZ^0 \bdelta^0 + \bepsilon,
\end{align}
where $\bY=(y_{1+h},\dots,y_{T+h})'$, $\bepsilon=(\epsilon_{1+h},\dots,\epsilon_{T+h})'$, $\bW=(\bw_1,\dots, \bw_T)'$, $\bdelta^0=({\bgamma^{0\prime}},\bbeta')'$, and ${\bZ}^0=({\bF}^0,\bW)=({\bz}_1^0,\dots,{\bz}_T^0)'$. We make the assumptions to derive the asymptotic theory.

\begin{ass}\normalfont\label{ass: augreg}
We have: 
(i) $\|\bgamma^{0}\|_{2} \le M$ and $\|\bbeta\|_{2} \le M$; 
(ii) $\E(\epsilon_{t+h}|y_t,\bz^0_t,y_{t-1},\bz^0_{t-1},\cdots)=0$ for any $h>0$; 
(iii) $T^{-1/2} {\bZ^0}'\bepsilon \CD N(\mathbf{0},\bSigma_{z^0\epsilon})$, $T^{-1} {\bZ^0}'\bZ^0 \CP \bSigma_{z^0}$, where $\bSigma_{z^0\epsilon}$ and $\bSigma_{z^0}$ are positive definite and bounded; 
(iv) 
$\E\| T^{-1/2}\bN^{-\frac{1}{2}} \sum_{t=1}^{T}\sum_{i=1}^{N}\bb_i^{0}e_{t,i}\bw_t' \|_2^2 \le M$ and 
$\E\| T^{-1/2}\bN^{-\frac{1}{2}}\sum_{t=1}^{T}\sum_{j=1}^{N}\bb_j^0[\epsilon_{t+h}e_{t,j}-\E(\epsilon_{t+h}e_{t,j})]\|_2^2 \le M$; 
(v) For all $t$, $\left|\E\left(\epsilon_{t+h}e_{t,i} \right)\right| \leq|{\tau}_{i}|$ for some $\tau_{i}$ such that $\sum_{i=1}^N|{\tau}_{ i}| \leq M$. 
\end{ass}
 
Assumption \ref{ass: augreg}(i) and (ii) are standard. 
Assumption \ref{ass: augreg}(iii) is for identification of $\bdelta^0$ and asymptotic normality of its estimator. Assumption \ref{ass: augreg}(iv) and (v) permit weak correlations between the errors in the factor model and in the augmented regression model. They imply 
$\bW'\bE\bB^0=O_p( \sqrt{TN^{\alpha_1}})$ and $\frac{1}{T}\bN^{-\frac{1}{2}}{\bB^0}'{\bE}'\bepsilon =O_p\left(\frac{1}{\sqrt{T}}+\frac{1}{\sqrt{N^{\alpha_r}}}\right)$.

\subsection{Asymptotic normality}

In model \eqref{augmodel0_mat}, unknown $\bff_t^0$ is replaced with the estimated factor $\hat{\bff}_t$. Then \eqref{augmodel0_mat} is written as
\begin{align}\label{model FFR}
    \bY=\hat{\bF} \bgamma^0 + \bW \bbeta + \bu= \hat\bZ \bdelta^0 + \bu,
\end{align}
where $\bu=\bepsilon-(\hat\bF-\bF^0)\bgamma^0$ and $\hat{\bZ}=(\hat{\bF},\bW)=(\hat{\bz}_1,\dots,\hat{\bz}_T)'$. 
The estimator of $\bdelta^0$ is obtained as $\hat\bdelta = (\hat\bZ'\hat\bZ)^{-1}\hat\bZ'\bY$. Using \eqref{model FFR}, we  have
\begin{align}\label{dd0}
\sqrt{T}(\hat\bdelta - \bdelta^0)  = (T^{-1}\hat\bZ'\hat\bZ)^{-1} T^{-\frac{1}{2}}\hat\bZ'\bepsilon
- (T^{-1}\hat\bZ'\hat\bZ)^{-1} T^{-\frac{1}{2}}\hat\bZ'(\hat\bF-\bF^0)\bgamma^0.    
\end{align}
In the right-hand side of \eqref{dd0}, the first term tends to a normal distribution with mean zero, and the second term causes potential bias in $\hat\bdelta$.

The following theorem ensures that the effect of the replacement by the PC estimator is asymptotically negligible and a test for a general restriction on $\bgamma^0$ 
is asymptotically valid. 

\begin{thm}\label{thm:forecast1}
Suppose that Assumptions \ref{cond:eigen}--\ref{ass: augreg} hold.
If $\frac{N^{1-\alpha_r}}{\sqrt{T}} \to  0$ and $\sqrt{T}N^{-\alpha_r} \to 0$, we have 
$\sqrt{T} ( \hat{\bdelta}-\bdelta^0) \CD N(\mathbf{0}, \bSigma_{\delta^0})$, where $\bSigma_{\delta^0}=\bSigma_{z^0}^{-1}\bSigma_{z^0\epsilon}\bSigma_{z^0}^{-1}$.
\end{thm}

In general, $\bSigma_{\delta^0}$ can be estimated by $\hat{\bSigma}_{\delta^0}=(T^{-1}\hat{\bZ}^{\prime}\hat{\bZ})^{-1}
\hat{\bSigma}_{z^0 \epsilon}
(T^{-1}\hat{\bZ}^{\prime}\hat{\bZ})^{-1}
$, where $\hat{\bSigma}_{z^0 \epsilon}$ is a consistent estimator. 
For example, when $\epsilon_t$ is heteroskedastic, we may use $\hat{\bSigma}_{z^0 \epsilon}=
T^{-1}\sum_{t=1+h}^{T+h}\hat\bz_t \hat{\epsilon}_t^2 \hat\bz_t'$, where $\hat{\epsilon}_{t+h}=y_{t+h}-\hat{\bdelta}'\hat{\bz}_t$.

The rate conditions for Theorem \ref{thm:forecast1} (and Theorem \ref{thm:forecast2} below) are identical to the corresponding results for approximations with data-dependent rotations, such as $(\hat{\bF},\hat{\bB}) \approx (\bF^* \hat{\bH}_4, \bB^* \hat{\bQ}')$, where $\hat{\bH}_4\hat{\bQ}\CP \bI_r$ by \citet[Lemma 3]{BaiNg2023}. 
They imply the associated approximations $\hat{\bgamma}\approx \hat{\bH}_4^{-1} \bgamma^{*}$ and $\hat{\bgamma}\approx \hat{\bQ} \bgamma^{*}$, respectively; see Lemma \ref{lem:delta34} of Supplementary Material. 
However, for such approximations the ``parameter'' estimated by $\hat{\bgamma}$ depends on $\bX$, and thus the only justifiable inference is about the joint significance, $H_0: \bgamma^{*}=\mathbf{0}$.

\subsection{Forecasting}

We consider the $h$-step ahead forecast. 
From \eqref{augmodel0}, the expectation of $y_{T+h}$ conditional on $\{\bz_T^0,\dots,\bz_{1}^0\}$ is computed as
$y_{T+h \mid T}:=\mathbb{E}(y_{T+h} \mid \bz_T^0,\dots,\bz_{1}^0)
=\bdelta^{0\prime} \bz_T^0$, 
which is the infeasible $h$-step ahead forecast of $y$ at time $T$. After estimating the forecast regression \eqref{model FFR}, $y_{T+h \mid T}$ is estimated by $\hat{y}_{T+h \mid T}=\hat{\bdelta}^{\prime} \hat{\bz}_T$. 
Thus the estimation error, $\hat{y}_{T+h \mid T}-y_{T+h \mid T}$, is derived as
$
\hat{y}_{T+h \mid T}-y_{T+h \mid T}=(\hat{\bdelta}-\bdelta^0)^{\prime}\hat{\bz}_T
+\bgamma^{\ast \prime} \bH^{-1 \prime}(\hat{\bff}_T -\bH' {\bff}_T^*).
$

\begin{thm}\label{thm:forecast2}
Suppose that Assumptions \ref{cond:eigen}--\ref{ass: augreg} hold. If $\sqrt{T}N^{-\alpha_r} \to 0$, $\frac{1}{2}<\alpha_r$, and $\frac{N^{\frac{3}{2}-\alpha_{r}}}{T}  \to 0$, we have 
$\frac{\hat{y}_{T+h \mid T}-y_{T+h \mid T}}{{\sigma}_{T+h \mid T}} \stackrel{d}{\longrightarrow} N(0,1)$, 
where ${\sigma}_{T+h \mid T}^2=
T^{-1} {\bz}_T^{0\prime}\bSigma_{\delta^0}{\bz}_T^0
+{\bgamma}^{0 \prime} 
\bD^{-1}\bN^{-1/2}\bGamma_T\bN^{-1/2}\bD^{-1}
{\bgamma}^0$.
\end{thm}

We may estimate ${\sigma}_{T+h \mid T}^2$ by
$\hat{\sigma}_{T+h \mid T}^2=
T^{-1} \hat{\bz}_T^{\prime}\hat{\bSigma}_{\delta^0}\hat{\bz}_T
+\hat{\bgamma}^{ \prime} 
(\hat{\bB}'\hat{\bB})^{-1/2}\hat{\bGamma}_T(\hat{\bB}'\hat{\bB})^{-1/2}
\hat{\bgamma}$ with some consistent estimators 
$\hat{\bSigma}_{\delta^0}$ and $\hat{\bGamma}_T$.
The rate conditions for Theorem \ref{thm:forecast2} are identical to those with data-dependent rotations; see Lemma \ref{lem:yh-y34} of Supplementary Material.

Finally, we consider the out-of-sample forecast. By $y_{T+h}=y_{T+h|T}+\epsilon _{T+h}$, we define the out-of-sample forecast error as 
$
v_{T+h|T}=\hat{y}_{T+h|T}-y_{T+h}=\hat{y}_{T+h|T}-y_{T+h|T}-\epsilon _{T+h}.
$
Assuming that $\epsilon _{t}\sim \text{i.i.d.}N\left( 0,\sigma _{\epsilon
}^{2}\right) $, we have $v_{T+h|T}\sim N( 0,\sigma _{\epsilon
}^{2}+\sigma _{T+h|T}^{2}) $ approximately. The variance of $v_{T+h|T}$ can be
estimated by $\hat{\sigma}_{\epsilon }^{2}+\hat{\sigma}_{T+h|T}^{2}$, where $\hat{\sigma}_{\epsilon }^{2}=T^{-1}\sum_{t=1+h}^{T+h}\hat{\epsilon}_{t}^{2}$
and $\hat{\sigma}_{T+h|T}^{2}$ is given above. This result yields the
prediction interval for $y_{T+h|T}$ centered at $\hat{y}_{T+h|T}$. If normality
assumption is not plausible, a bootstrap confidence interval should be
computed; see \cite{GoncalvesPerron2014}.

\section{Monte Carlo Experiments}\label{sec: MC}

\subsection{Approximate factor models}\label{subsec:fac}

We discuss estimation accuracy of the PC estimator $(\hat\bF,\hat\bB)$ against the pseudo-true parameter $(\bF^0,\bB^0)$ and other existing data-dependent rotations, $(\bF^*\hat\bH,\bB^{*}\hat\bH^{-1\prime})$ and $(\bF^*\hat\bH_4,\bB^{*}\bQ')$. We also examine the normal approximation for $\hat{\bff}_t$ and $\hat{\bb}_i$. Furthermore, we investigate the size of tests, such as $H_{0}:f_{t,k}=f_{t,k}^{0}$ and $H_{0}:b_{i,k}=b_{i,k}^{0}$, which are only justified in our proposed approach. 

\subsubsection{Experimental design\label{mc: fac des}}

We choose $r=2$ throughout, and generate $f_{t,k}\sim \text{i.i.d.}U(\mu_{k}-\sqrt{3},\mu_{k}+\sqrt{3})$ with $\mu_{k}=1$ for $t=1,\dots ,T$ and $k=1,2$, to form $\mathbf{F}=\left(f_{t,k}\right)$. 
Applying the Gram-Schmidt procedure to $\mathbf{F}$, we obtain $\mathbf{F}^{0}$ such that $T^{-1}\mathbf{F}^{0\prime}\mathbf{F}^{0}=\mathbf{I}_{2}$. 
We next generate two types of loading matrices: 
(i) Non-sparse factor loadings: $\mathbf{B}^{0}=\left(  \mathbf{B}%
_{1}^{0},\mathbf{B}_{2}^{0}\right)  =\left(  N^{(\alpha_{1}-1)/2}%
{g}_{i,1},N^{(\alpha_{2}-1)/2}{g}_{i,2}\right)  $, where $g_{i,1}$
is $2$ for $i=1,\dots,N/2$ and $1$ for the rest, and $g_{i,2}$ is $1/2$ for
$i=1,\dots,N/2$ and  $-1$ for the rest; 
(ii) Sparse factor loadings: $\mathbf{B}^{0}=\left(  \mathbf{B}_{1}^{0}%
,\mathbf{B}_{2}^{0}\right)  $, where $b_{i,1}^{0}$ is 2 for $i=1,\dots,N_{1}$ and
$0$ for the rest, and $b_{i,2}^{0}$ is $1$ for $i=1,\dots,N_{2}/2$, $-1$ for
$i=N_{2}/2+1,\dots ,N_{2}$, and zero for the rest, where $N_{k}$ are the smallest even
number not less than $N^{\alpha_{k}}$ for $k=1,2$. These designs ensure that $\mathbf{B}_{k}^{0\prime}\mathbf{B}_{k}^{0}\asymp
N^{\alpha_{k}}$. We construct $\mathbf{F}^{\ast}=\mathbf{F}^{0}\mathbf{H}^{-1}$ and $\mathbf{B}^{\ast}=\mathbf{B}^{0}\mathbf{H}^{\prime}$, where $\mathbf{H=}\bigl(
\begin{smallmatrix}
   1 & 0.5 \\
   0.5 & 2
\end{smallmatrix}
\bigl)$. Then the data is generated by $\mathbf{X}=\mathbf{F}^{\ast}\mathbf{B}^{\ast\prime}+\mathbf{E}$,
where $\mathbf{E}=\left(  e_{t,i}\right)  $ with $e_{t,i}\sim \text{i.i.d.}N(0,\sigma
_{e}^{2})$ and $\sigma_{e}^{2}=0.5$. 
We consider six models: $\left(  \alpha_{1},\alpha_{2}\right)  =\left(  1.0,1.0\right)  $,
$\left(  1.0,0.9\right)  $, $\left(  1.0,0.8\right)  $, $\left(0.9,0.7\right)  $, $\left(  0.8,0.6\right)  $, and $\left(  0.7,0.5\right)$, where $\left(  1.0,1.0\right)  $ is the SF model. We consider sample sizes, $T=N=50,100,200,500$. All the results are based on 50,000 replications. Using the correlation between
$\mathbf{F}^{0}$ and $\mathbf{\hat{F}}$, the signs and the orders of the $r$
columns of $\mathbf{\hat{F}}$\ and $\mathbf{\hat{B}}$ are adjusted to fit to
those of $\mathbf{F}^{0}$, when necessary, so that $\mathbf{\hat{F}}$\ and
$\mathbf{\hat{B}}$ can be regarded as the estimates of $\mathbf{F}^{0}$ and
$\mathbf{B}^{0}$ (i.e. no sign indeterminacy). This adjustment is for experimental purposes only; this may be irrelevant in practice. To save space, we only report the results for the models with
non-sparse loadings. The results with sparse loadings are very similar and are available in Section \ref{subsec:addMC} of the Supplementary Materials. 

We also include, as an additional comparison target, a Procrustes-type rotation in the spirit of \citet{FanEtAl2024}, expressed here in our notation. 
Recall $\bC^{*}=\bF^{*}\bB^{*\prime}=\bF^{0}\bB^{0\prime}$ together with the PC normalization $T^{-1}\bF^{0\prime}\bF^{0}=\bI_r$ and $\bB^{0\prime}\bB^{0}=\bLambda$. 
Hence, with $\bU^{0}=\bF^{0}/\sqrt{T}$ and $\bV^{0}=\bB^{0}\bLambda^{-1/2}$, the scaled common component has the singular value decomposition $\bC^{*}/\sqrt{T}=\bU^{0}\bLambda^{1/2}\bV^{0\prime}$, where $\bU^{0\prime}\bU^{0}=\bV^{0\prime}\bV^{0}=\bI_r$ and the singular values are $\diag(\bLambda^{1/2})$. 
This decomposition is unique up to column signs because the
$\lambda_k$ are distinct. Let $\hat{\bU}=\hat{\bF}/\sqrt{T}$ collect the sample left
singular vectors implied by the PC estimator, so that $\hat{\bU}'\hat{\bU}=\bI_r$. The
orthogonal Procrustes alignment of $\hat{\bU}$ to $\bU^{0}$ is 
$\hat{\bR}_U
=\argmin_{\bR:\,\bR'\bR=\bI_r}\big\|\hat{\bU}\bR-\bU^{0}\big\|_{\F}
=\argmin_{\bR:\,\bR'\bR=\bI_r}\tfrac{1}{\sqrt{T}}\big\|\hat{\bF}\bR-\bF^{0}\big\|_{\F}
=\bW_1\bW_2'$, 
where $\hat{\bU}'\bU^{0}=T^{-1}\hat{\bF}'\bF^{0}=\tilde{\bQ}=\bW_1\bSigma_{\mathrm{sv}}\bW_2'$
is a singular value decomposition; equivalently, $\hat{\bR}_U$ is the orthogonal polar
factor of $\tilde{\bQ}$. Because $\sqrt{T}\,\bU^{0}=\bF^{0}=\bF^{*}\bH$, the factor- and
loading-side rotations of the starting parameter implied by this alignment are
$\hat{\bR}_{\mathrm{Fan},F}=\bH\hat{\bR}_U'$ and $\hat{\bR}_{\mathrm{Fan},B}=\bH^{-1\prime}\hat{\bR}_U'$, 
so that $\bF^{*}\hat{\bR}_{\mathrm{Fan},F}=\bF^{0}\hat{\bR}_U'$ and
$\bB^{*}\hat{\bR}_{\mathrm{Fan},B}=\bB^{0}\hat{\bR}_U'$ are the pseudo-true parameters
realigned by the common orthogonal matrix $\hat{\bR}_U'$, and the common component is
preserved: $(\bF^{*}\hat{\bR}_{\mathrm{Fan},F})(\bB^{*}\hat{\bR}_{\mathrm{Fan},B})'
=\bF^{0}\bB^{0\prime}=\bC^{*}$. 
Like $\hat{\bH}$, $\hat{\bH}_4$, and $\hat{\bQ}$, the pair $(\hat{\bR}_{\mathrm{Fan},F},\hat{\bR}_{\mathrm{Fan},B})$ is data-dependent and is recomputed in every replication.

\subsubsection{Estimation accuracy and normal approximation}
\begin{figure}[h!]
	\centering
	\includegraphics[clip,width=14.3cm]{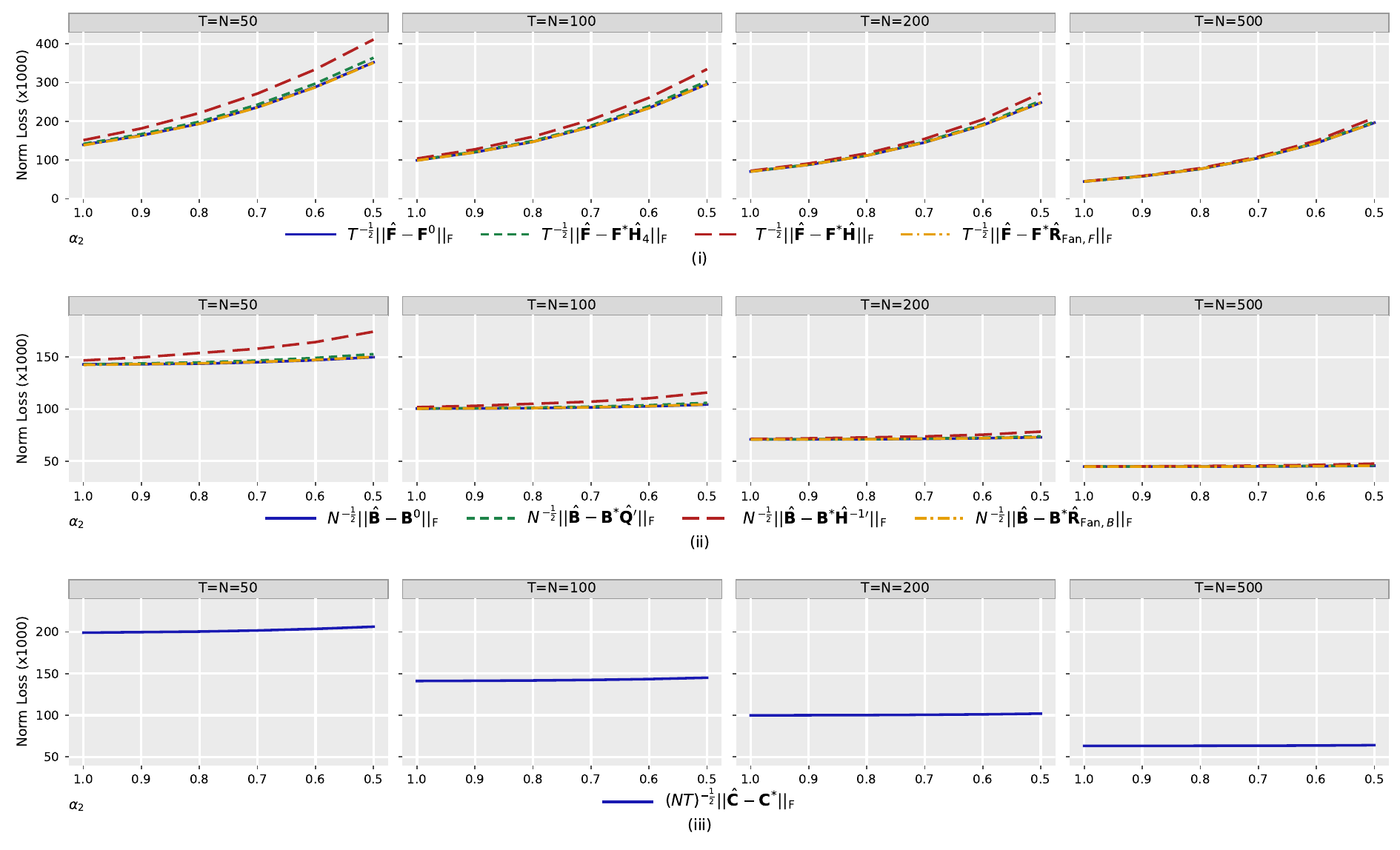}
	\vspace{-5mm}
\caption{Average norm losses of (i)$\hat{\bF}$, (ii)$\hat{\bB}$, and (iii)$\hat{\bC}$
relative to the pseudo-true and data-dependent targets based on
$\hat{\bH}_4$, $\hat{\bQ}$, $\hat{\bH}$, and the Procrustes rotation of
\citet{FanEtAl2024}.}
 \label{fig:1}
\end{figure}

To assess the convergence results for the PC estimator
$(\hat{\bF},\hat{\bB})$ and their product, Figure \ref{fig:1} summarizes
the averages over the replications of the relevant norm losses.
Figure \ref{fig:1}(i) reports
$T^{-\frac{1}{2}}\|\hat{\bF}-\bF^{0}\|_{\F}$,
$T^{-\frac{1}{2}}\|\hat{\bF}-\bF^{*}\hat{\bH}_{4}\|_{\F}$,
$T^{-\frac{1}{2}}\|\hat{\bF}-\bF^{*}\hat{\bH}\|_{\F}$, and
$T^{-\frac{1}{2}}\|\hat{\bF}-\bF^{*}\hat{\bR}_{\mathrm{Fan},F}\|_{\F}$,
while Figure \ref{fig:1}(ii) compares
$N^{-\frac{1}{2}}\|\hat{\bB}-\bB^{0}\|_{\F}$,
$N^{-\frac{1}{2}}\|\hat{\bB}-\bB^{*}\hat{\bQ}'\|_{\F}$,
$N^{-\frac{1}{2}}\|\hat{\bB}-\bB^{*}\hat{\bH}^{-1\prime}\|_{\F}$, and
$N^{-\frac{1}{2}}\|\hat{\bB}-\bB^{*}\hat{\bR}_{\mathrm{Fan},B}\|_{\F}$.
Figure \ref{fig:1}(iii) shows
$(NT)^{-\frac{1}{2}}\|\hat{\bC}-\bC^{*}\|_{\F}$.

Each block in the figure is the result for a particular sample size.
Within the block, the vertical axis shows the magnitude of the average
norm loss, and the horizontal axis indicates the value of $\alpha_2$,
which identifies one of the six models considered, hence, the weakness of
the model -- the smaller the value, the weaker the factor model.

It is clear from Figure \ref{fig:1}(i) and (ii) that the norm loss of the
PC estimator $(\hat{\bF},\hat{\bB})$ against $(\bF^0,\bB^0)$ is the
smallest compared with those against
$(\bF^{*}\hat{\bH}_4,\bB^{*}\hat{\bQ}')$ and
$(\bF^{*}\hat{\bH},\bB^{*}\hat{\bH}^{-1\prime})$, closely followed by the
former approximation. The approximation by $(\bF^0,\bB^0)$ improves as the
model weakens especially when the sample size is small. This superiority
eventually disappears as the sample size increases. These patterns are
consistent with Theorem \ref{0}: under the simulation design with $N=T$,
the rates for the pseudo-true targets coincide with those for the
data-dependent targets in Lemma \ref{lem fh4bq}. A similar comment applies to
(iii). 

The losses relative to the targets considered by \citet{FanEtAl2024},
$\bF^{*}\hat{\bR}_{\mathrm{Fan},F}$ and
$\bB^{*}\hat{\bR}_{\mathrm{Fan},B}$, are very close to
those relative to the pseudo-true parameters in every model and at every
sample size; the two losses differ by less than one percent throughout.
Tables~\ref{tab:D_rotations} and \ref{tab:D_rotations_sd} in
Section~\ref{sec:D_fan} of the Supplementary Materials show that the Monte
Carlo means of $\hat{\bR}_{\mathrm{Fan},F}$ differ from the fixed rotation
$\bH$ by less than $5\times10^{-4}$ entrywise in every model and at every
sample size, while the across-replication standard deviations of its entries
are at most about $4\times10^{-2}$ at $T=N=50$ and fall to about
$5\times10^{-3}$ at $T=N=500$ 
(the paired loading-side rotation $\hat{\bR}_{\mathrm{Fan},B}=\hat{\bR}_{\mathrm{Fan},F}^{-1\prime}$
is correspondingly close to $\bH^{-1\prime}$). These findings indicate that, in
our design, the population counterpart of the Procrustes alignment is
effectively the fixed rotation that maps $(\bF^{*},\bB^{*})$ to the pseudo-true
representation $(\bF^{0},\bB^{0})$, and the sample Procrustes rotations are
tightly concentrated around this population target. This is because the scaled
common component factors as
$T^{-1/2}\bC^*=(T^{-1/2}\bF^0)\bLambda^{1/2}(\bB^0\bLambda^{-1/2})'$, which is
its singular value decomposition; since the $\lambda_k$, $k\in[r]$, are distinct,
it recovers the PC-normalized population representation $(\bF^0,\bB^0)$ up to
column signs; see Section~\ref{sec:D_fan}. The full set of results under the
\citet{FanEtAl2024}-type rotation, including the size of the $t$-tests under all
centerings and the average and standard deviation of the rotation matrices, is
reported there.

\subsubsection{Tests of linear restrictions}

We report the performance of testing for $H_{0}:f_{1,k}=f_{1,k}^{0}$,  $H_{0}:b_{1,k}=b_{1,k}^{0}$, and $H_0: c_{1,1}=c_{1,1}^*$ using the statistics, 
$
z_{f,k}  = \sqrt{\lambda_{k}(\bGamma_{1}^{-1})_{k,k}}(\hat{f}_{1,k}-f_{1,k}^{0}),~~~
z_{b,k}  = \sqrt{T(\bPhi_{1}^{-1})_{k,k}}(  \hat{b}_{1,k}-b_{1,k}^{0}), ~~~
z_{c}=\frac{ \hat{c}_{1,1}-c_{1,1}^{\ast}}{\sigma_{c}},
$
respectively, where $\lambda_{k}=(\bB^{0\prime}\bB^{0})_{kk}$ and $\sigma_{c}^{2}=\sigma_{e}^{2}[\mathbf{b}_{1}^{0\prime}(\mathbf{B}^{0\prime}\mathbf{B}^{0})^{-1}\mathbf{b}_{1}^{0}+T^{-1}\mathbf{f}_{1}^{0\prime}\mathbf{f}_{1}^{0}]$. 
We evaluate the empirical size of the nominal 5\% tests by computing the
frequencies with which the absolute values of the statistics exceed 1.96
over the Monte Carlo replications. These tests are fixed-target tests for
the PC-normalized population factors, loadings, and common component.

\begin{figure}[h!]	
	\centering
	\includegraphics[clip,width=14.3cm]{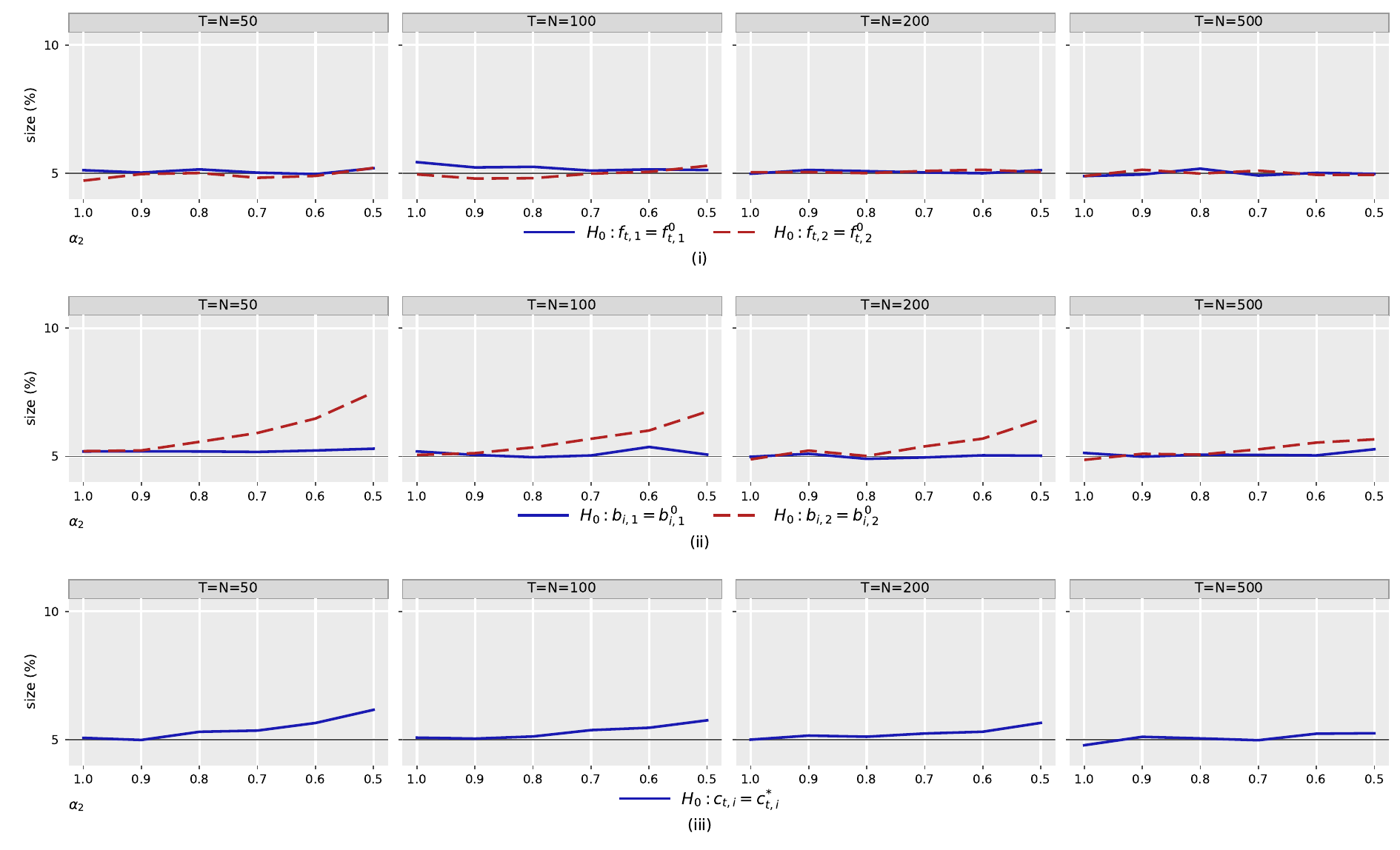}
	\vspace{-5mm}
	\caption{Empirical size of the 5\% level tests for (i) $H_0: f_{1,k}=f_{1,k}^0$;
 (ii)  $H_0: b_{1,k}=b_{1,k}^0$;
(iii) $H_0: c_{1,1}=c_{1,1}^*$.}
\label{fig:3}
\end{figure}

Figure \ref{fig:3}(i) shows that the empirical size of the proposed factor
test is close to the nominal level across the designs. Figure \ref{fig:3}(ii)
shows moderate over-rejection for the loading test in the weakest models, but
the distortion decreases as the sample size increases. 
A similar pattern is observed for the common-component test in Figure \ref{fig:3}(iii).

The joint chi-square diagnostics, which underlie Wald-type tests for joint
restrictions on \(\bff_t^0\) or \(\bb_i^0\), are reported in
Appendix~\ref{app:joint-normal-centerings}. They show that centering at the
pseudo-true parameters gives the most accurate joint normal approximation,
with the Procrustes centering of \citet{FanEtAl2024} behaving very similarly.
The conventional data-dependent centerings, especially the one based on
\(\hat{\bH}\), deteriorate as the factor model weakens.

\subsection{Factor augmented regressions}\label{subsec:facaug}

We examine the normal approximation of the coefficient estimates and the size of
componentwise tests, such as $H_{0}:\gamma_{k}=\gamma_{k}^{0}$. These tests are
justified only under our proposed fixed-target interpretation, because $\gamma_k^0$
is the coefficient on the $k$th PC-normalized population factor, rather than on a
data-dependently rotated factor. The forecast intervals derived in
Theorem~\ref{thm:forecast2} are also important practical implications of the theory,
but their Monte Carlo evaluation concerns coverage of the fitted conditional mean
$y_{T+h|T}$ and the actual future value $y_{T+h}$, which aggregate factor and
coefficient estimation error through the forecast equation. We therefore report the
corresponding forecast-coverage results in Appendix~\ref{app:forecast-interval-mc}.

\subsubsection{Experimental design}

Using $\mathbf{f}_{t}^{0}$ and $\mathbf{H}$ in Section \ref{mc: fac des}, we generate
\begin{align}
y_{t+h}=\mathbf{f}_{t}^{0\prime}\bgamma^{0}+\mathbf{w}_{t}
^{\prime}\boldsymbol{\beta}+\epsilon_{t+h},~~~t=1,\dots ,T,
\label{mc: faug0}
\end{align}
where $\bgamma^{0}=\left(  1,\dots ,1\right)'\in\bbR^r$, $\boldsymbol{\beta}=\left(
1,\dots ,1\right)'\in\bbR^L$, $\epsilon_{t}\sim \text{i.i.d.}N\left(
0,\sigma_{\epsilon}^{2}\right)  $, and $\mathbf{w}_{t}=\left(  w_{t,1},\dots ,w_{t,L}\right)  ^{\prime}$ with
$w_{t,\ell}=\mu_{w}+\sum_{k=1}^{r}\rho_{\ell}[  f_{t,k}^{0}-\E\left(  f_{k}^{0}\right)  ]  +\varepsilon_{w,t,\ell}$, where $\varepsilon_{w,t,\ell
}\sim \text{i.i.d.}N(  0,1-r\rho_{\ell}^2)  $, for $\ell=1,\dots ,L$, so that 
$\rho_\ell$ controls the association between $f_{t,k}^0$ and $w_{t,\ell}$. 
Note that $\mathbf{f}_{t}^{0\prime}\bgamma%
^{0}=\mathbf{f}_{t}^{\ast\prime}\bgamma^{\ast}$\ with
$\bgamma^{\ast}=\mathbf{H}\bgamma^{0}$ and $\mathbf{f}_{t}^{\ast
\prime}=\mathbf{f}_{t}^{0\prime}\mathbf{H}^{-1}$ in \eqref{mc: faug0}.
We choose $L=1$, $\mu_{w}=1$, $\rho_{\ell}=0.5$ for all $\ell$, and
$\sigma_{\epsilon}=1$, for which the population $R^{2}$ of the augmented
regression (\ref{mc: faug0}) is approximately $3/4$. We have chosen $\alpha=0.05$ and $h=1$.
All the results are obtained by 50,000 replications.

The infeasible estimator of $\bgamma^{0}$ is defined as 
$\hat{\bgamma}^{0}=(  \sum_{t=1}^{T-h}\mathbf{\tilde{f}}_{t}^{0}\mathbf{\tilde{f}}_{t}^{0\prime})  ^{-1}\sum_{t=1}%
^{T-h}\mathbf{\tilde{f}}_{t}^{0}y_{t+h}$, where $\mathbf{\tilde{f}}%
_{t}^{0\prime}$ is $t$th rows
of the $T\times r$ matrix $\mathbf{M}_{w}\mathbf{F}^{0}$ with $\mathbf{M}_{w}=\mathbf{I}_{T}-\mathbf{W}%
(\mathbf{W}^{\prime}\mathbf{W})^{-1}\mathbf{W}^{\prime}$ and $\bW=(\bw_{1},\dots ,\bw_{T})'$. The feasible estimator of $\bgamma^{0}$ is defined by replacing $\bF^0$ with $\hat\bF$ in $\hat{\bgamma}^{0}$, and denoted by $\hat{\bgamma}$.

\subsubsection{Estimation accuracy and normal approximation}

\begin{figure}[h!]
  \centering
  \includegraphics[clip,width=14.3cm]{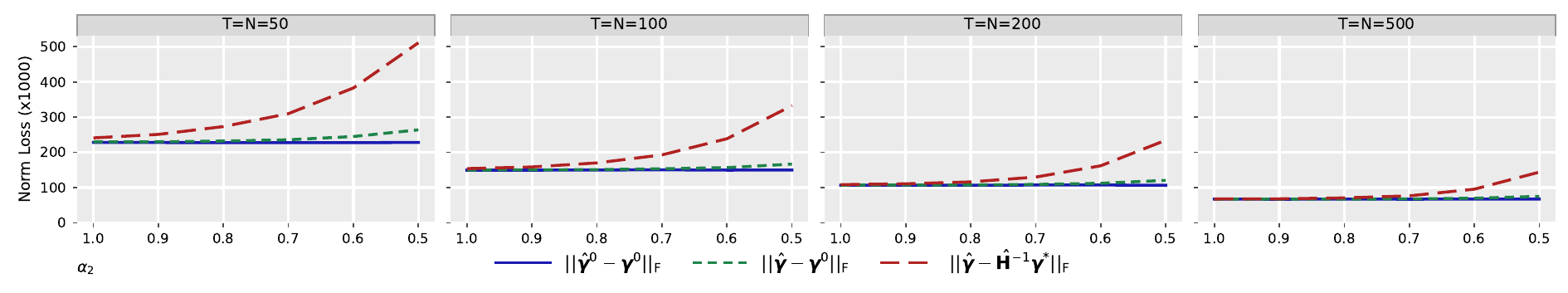}
  \vspace{-5mm}
  \caption{Average norm losses of the regression coefficient estimators:
  $\|\hat{\bgamma}^{0}-\bgamma^{0}\|_{\F}$,
  $\|\hat{\bgamma}-\bgamma^{0}\|_{\F}$, and
  $\|\hat{\bgamma}-\hat{\bH}^{-1}\bgamma^{*}\|_{\F}$.}
  \label{fig:4}
\end{figure}

Figure \ref{fig:4} reports the norm losses
$\|\hat{\bgamma}^{0}-\bgamma^{0}\|_{\F}$,
$\|\hat{\bgamma}-\bgamma^{0}\|_{\F}$, and
$\|\hat{\bgamma}-\hat{\bH}^{-1}\bgamma^{*}\|_{\F}$.
Replacing $\bF^0$ with its consistent estimator $\hat{\bF}$ has little
effect on estimation accuracy. By contrast, centering at the conventional
data-dependent target $\hat{\bH}^{-1}\bgamma^{*}$ leads to substantially
larger losses, especially as the underlying factor model becomes weaker.
All norm losses decrease as the sample size increases.

\subsubsection{Test of linear restrictions}

\begin{figure}[h!]	
	\centering
	\includegraphics[clip,width=14.3cm]{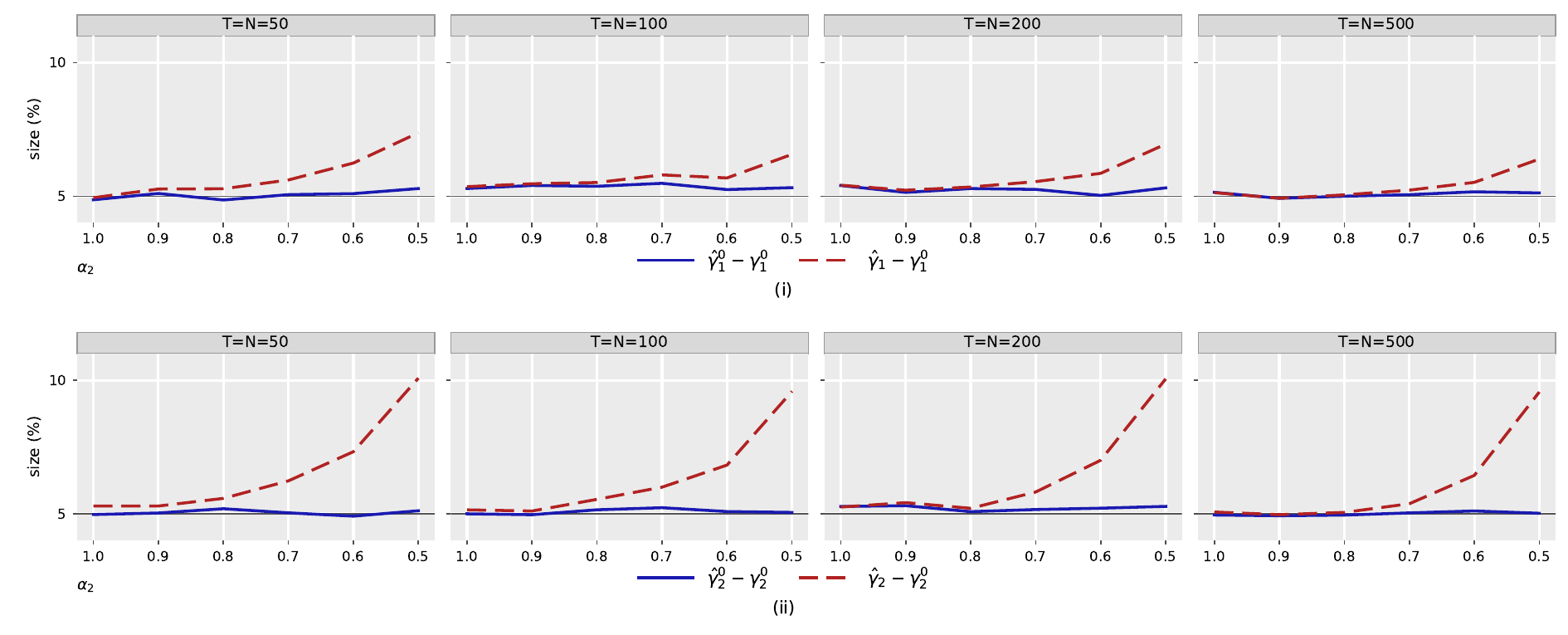}
	\vspace{-5mm}
	\caption{Size of the 5\% level tests for (i) $H_0: \gamma_{1}=\gamma_{1}^0$; (ii) $H_0: \gamma_{2}=\gamma_{2}^0$.}
 \label{fig:6}
\end{figure}

We report the empirical size of the 5\% level tests for
$H_{0}:\gamma_{k}=\gamma_{k}^{0}$ against
$H_{1}:\gamma_{k}\neq\gamma_{k}^{0}$, $k=1,\dots,r$, using the statistic
$z_k(\Delta\gamma_k)
=
\sqrt{T}\Delta\gamma_k/\sigma_{\gamma^{0},kk}^{1/2}$,
where $\Delta\gamma_k=\hat{\gamma}_{k}^{0}-\gamma_{k}^{0}$ or
$\Delta\gamma_k=\hat{\gamma}_{k}-\gamma_{k}^{0}$, and
$\sigma_{\gamma^{0},kk}$ is the $k$th diagonal element of
$\boldsymbol{\Sigma}_{\gamma^{0}}$. The rejection frequencies are computed
using the standard normal critical value, $|z_k(\Delta\gamma_k)|>1.96$.
Figure \ref{fig:6} is consistent with the theory: the tests have empirical
sizes close to the nominal level, except in the weakest model, where the
conditions of Theorem \ref{thm:forecast1} are not satisfied.

The joint chi-square diagnostics, which underlie Wald-type tests for
joint restrictions on \(\bgamma^0\), are reported in
Appendix~\ref{app:far-joint-normal}. They show that the approximation relative
to the fixed pseudo-true coefficient vector is accurate except in the weakest
design, whereas centering at the conventional data-dependent target can lead to
severe over-rejection.

The corresponding Monte Carlo evidence for the confidence interval for the
conditional mean \(y_{T+h|T}\) and the prediction interval for the actual value
\(y_{T+h}\) is reported in Appendix~\ref{app:forecast-interval-mc}.

\section{Empirical application: oil-factor confidence intervals}
\label{subsec:application}
The fixed-target theory developed in this paper has several empirical implications. First, it justifies interpreting PC-normalized factors through inference on their loadings; see Appendix~\ref{subsec:addapplication} and \citet{UY2019inference}. Second, it provides a population interpretation for coefficients in factor-augmented regressions; see \citet{jiang2024Mw}. Here we focus on a third implication: constructing confidence intervals for a population factor obtained from economically motivated restrictions. 
This application is particularly useful in weak-factor settings, where conventional PC theory based on data-dependent rotations does not, in general, deliver confidence intervals for a fixed population factor.

In macroeconomics, factor models are often combined with economically motivated restrictions to give factors a structural interpretation. One of the leading examples is the oil-market factor of \citet{StockWatson2016}, which captures an oil-price shock that moves oil prices directly, with macroeconomic variables responding subsequently. 
Suppose that the first $p$ variables, $\{x_{t,1},\dots,x_{t,p}\}$, are oil price indices. To make the first factor the oil factor, restrict the corresponding loading block to
$\bB_{1}^{\dag} = (\bc, \bzero_{p\times (r-1)} )$, 
where the weights $\bc=(c_1,\dots,c_p)'$ are chosen by the researcher. 

Following \citet{StockWatson2016}, we first estimate the PC-normalized population factors and loadings, $\bF^{0}$ and $\bB^{0}$, and then impose the oil-factor restriction on the first $p$ rows of the loading matrix. In the application below, $pr<r^{2}$, so the restriction $\bB_{1}^{\dag}$ does not identify all factors and loadings. 
Define the restricted loading matrix as $\bB^{\dag}=[\bB_{1}^{\dag\prime},\bB_{2}^{0\prime}]^{\prime}$,
where $\bB_{2}^{0}$ consists of the last $N-p$ rows of $\bB^{0}$. 
The corresponding population factor is $\bF^{\dag}=\bF^{0}\bG$, 
where $\bG = \bB^{0\prime}\bB^{\dag}(\bB^{\dag\prime}\bB^{\dag})^{-1}$. 
The first column of $\bF^{\dag}$ is the oil factor. 
In general, due to overidentification,
$\bH\neq\bG\neq\mathbf{I}_{r}$ and $\bF^{\dag}\bB^{\dag\prime}\neq\bF^{0}\bB^{0\prime}$. 
The associated estimator is $\hat{\bF}^{\dag}=\hat{\bF}\hat{\bG}$, where $\hat{\bG}=\hat{\bB}^{\prime}\hat{\bB}^{\dag}(\hat{\bB}^{\dag\prime}\hat{\bB}^{\dag})^{-1}$ and 
$\hat{\bB}^{\dag}=[\bB_{1}^{\dag\prime},\hat{\bB}_{2}^{\prime}]^{\prime}$. 
Then 
$
\frac{1}{\sqrt{T}}\,
\big\|\hat \bF^{\dagger}-\bF^{\dagger}-(\hat \bF-\bF^{0})\bG\big\|_{\F}
=O_p\big(N^{(\alpha_1-\alpha_r)/2} \tilde{\Delta}_{NT}\big)=o_p(1),
$
with $\tilde{\Delta}_{NT}=N^{-\alpha_r}+N^{1-\alpha_r}/T$ as defined in (B.2) in Appendix.

This expansion shows why fixed-target asymptotics are needed for inference on the oil factor. 
Confidence intervals for the population oil factor $f^{\dag}_{1,t}$ can be constructed from the asymptotic distribution of $\bD\bN^{1/2}(\hat{\bff}_{t}-\bff_{t}^{0})$: writing $\bg_1$ for the first
column of $\bG$, we have
$
\hat f^{\dagger}_{1,t}-f^{\dagger}_{1,t}
=\bg_1'\big(\hat\bff_t-\bff_t^{0}\big)+O_p\!\big(\tilde{\Delta}_{NT}\big)$, with $
\lambda_1^{1/2}\tilde{\Delta}_{NT}=o_p(1).
$
By contrast, an approximation based on a data-dependent rotation, such as
$\hat{\bF}-\bF^{0}\tilde{\bH}$, would imply
$\frac{1}{\sqrt{T}}(\hat{\bF}^{\dag}-\bF^{0}\tilde{\bH}\bG )= \frac{1}{\sqrt{T}}(\hat{\bF}-\bF^{0}\tilde{\bH})\bG+o_{p}(1)$.
The resulting interval would be centered on
$\bF^{0}\tilde{\bH}\bG$, a data-dependent linear combination of $\bF^{0}$, 
rather than on the fixed population factor $\bF^{\dag}$.

We use the balanced panel dataset of \citet{StockWatson2016}, which contains
135 quarterly U.S. macroeconomic and oil-market series over 1987q4--2013q4,
including four oil price indices $(p=4)$. All series are transformed as in
\citet{StockWatson2016} and then standardized. 
The ED estimator of $r$ by \citet{Onatski2010} gives $\hat r=6$.
Using the procedure in Appendix~\ref{subsec:addapplication} with FDR controlled below $0.1$, the loading-based factor-strength estimates $\hat{\alpha}_k$, $k\in[r]$, range from $0.81$ to $0.96$. Following \citet{StockWatson2016}, let $c_\ell=\lambda_1^{1/2}/\sigma_\ell$ and in practice we replace it by the plug-in estimator $\hat{\lambda}_{1}^{1/2}/s_{\ell}$ for $\ell=1,\dots,p$, where $\hat{\lambda}_{1}$ is the largest eigenvalue of the sample covariance matrix and $s_{\ell}$ is the sample standard deviation of the $\ell$-th oil price index before standardization.

\begin{figure}[h!]	
  \centering
  \includegraphics[width=\textwidth]{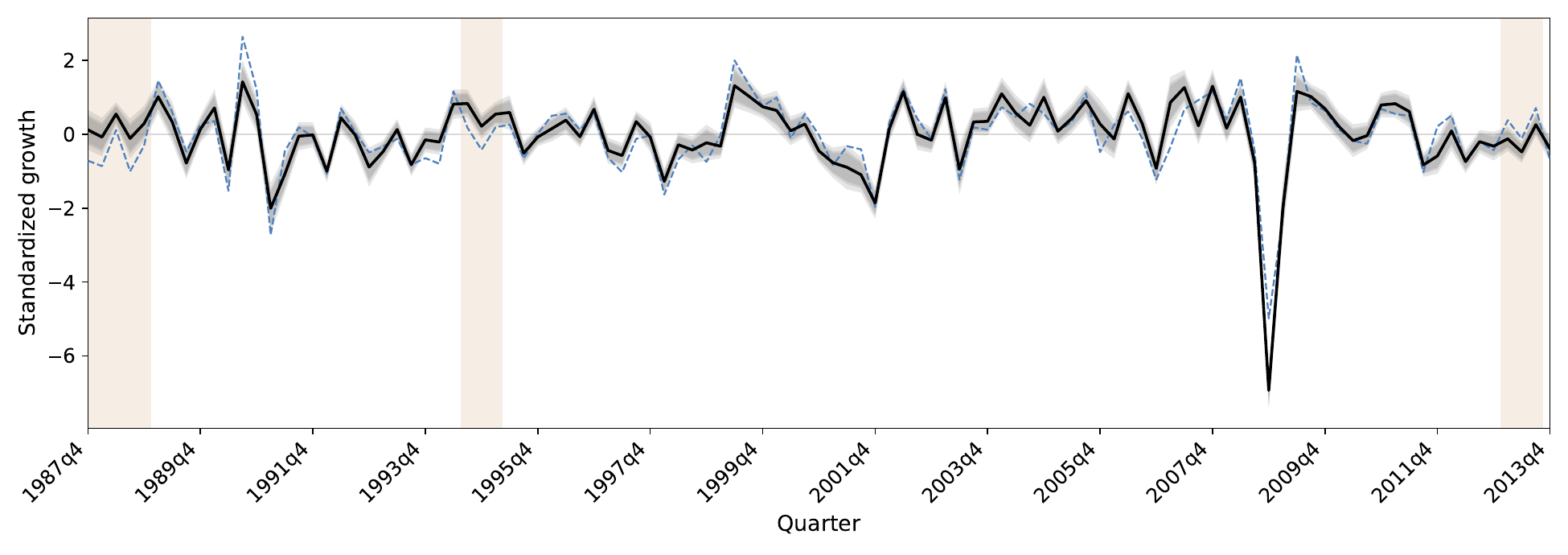}\\[-2pt]
  \includegraphics[width=0.95\textwidth]{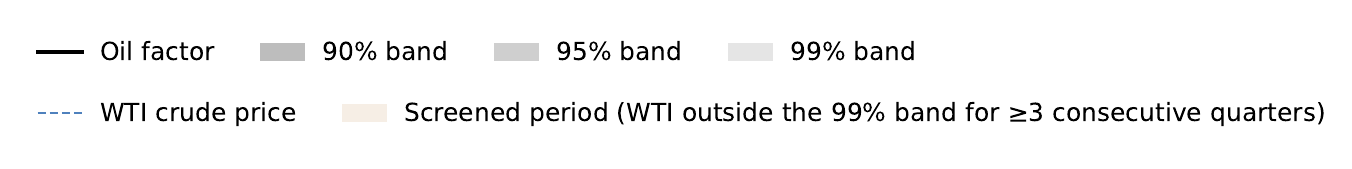}\\[-2pt]
\caption{Confidence bands for the oil factor $f^{\dagger}_{1,t}$ and WTI crude-price growth}
  \label{fig:sw2016}
\end{figure}

Figure~\ref{fig:sw2016} plots the estimated oil factor $\widehat{f}^{\dagger}_{1,t}$ as a
black solid line, together with its $90\%$, $95\%$ and $99\%$ confidence bands \textit{for the population oil factor} $f^{\dag}_{1,t}$, shown by the nested gray shading. WTI crude-price growth is plotted as the blue dashed line for reference. 
Three shaded vertical spans descriptively mark periods in which WTI lies outside the $99\%$ band for at least three consecutive quarters. 
All series are quarterly log-growth rates over 1987q4--2013q4.
The narrow bands illustrate that the restricted population factor is estimated with substantial precision.

Standardized WTI growth co-moves closely with the oil factor, but it often lies outside the factor's confidence band. This is expected: the factor pools four oil-price indices and is therefore estimated precisely, whereas any single index contains a larger idiosyncratic component. 
Thus, the confidence band should be interpreted as uncertainty about the population oil factor, not as a prediction band for individual oil-price indices.

The persistent deviations are economically interpretable. 
The oil factor tracks the common movement in the four oil-price indicators, whereas WTI is the landlocked Cushing benchmark. 
Hence deviations of WTI from the factor reflect WTI-specific basis movements and lead--lag timing among the indicators. 
The clearest case is 2013q1--2013q3, which is naturally interpreted against the post-2010 Cushing bottleneck documented by \citet{BorensteinKellogg2014}; the subsequent narrowing is consistent with new Cushing takeaway capacity coming online \citep{EIAspread2013}.
The earlier episodes, 1987q4--1988q4 and 1994q3--1995q1, predate the post-2010 WTI basis episode and are therefore better viewed as composition or timing divergences among the four oil-price measures.
 

\section{Conclusion}\label{sec:con}

In the literature, including \cite{Bai2003} and \cite{BaiNg2023}, the PC estimators have typically been analyzed relative to $(\bF^{*}\hat{\bH}, \bB^{*}\hat{\bH}^{-1\prime})$, where $(\bF^{*},\bB^{*})$ is a starting representation in model \eqref{model:AFM*} and $\hat{\bH}$ is defined in \eqref{Hhat}. Since these ``rotated parameters'' depend on the PC estimators through $\hat{\bH}$, such results do not establish consistency for fixed population parameters. This raises a natural question: what fixed population parameters do the PC estimators estimate?

To answer this question, we have established several theoretical results. First, under quite general conditions, we prove that, for any starting representation $(\bF^{*},\bB^{*})$ of the common component, there exists a rotation matrix $\bH$, unique up to column signs, that maps the original model \eqref{model:AFM*} to the pseudo-true model \eqref{model:AFM0} satisfying the PC normalization \eqref{pc1} (Theorem \ref{thm:rotation}). Although $\bH$ depends on the chosen starting representation, the resulting pseudo-true representation is invariant to that choice, up to simultaneous column sign changes. We then establish the asymptotic theory of the PC estimators relative to the pseudo-true parameters $(\bF^{0},\bB^{0})$ (Theorems \ref{0} and \ref{thm dist}).

Our asymptotic theory accommodates WF models in which the $r$ largest eigenvalues of $T^{-1}\bB^{*}\bF^{*\prime}\bF^{*}\bB^{*\prime}$ diverge at possibly different rates $N^{\alpha_k}$, with $\alpha_1\geq\dots\geq\alpha_r>0$. Within this framework, Theorem \ref{thm dist} provides asymptotic normality for each pseudo-true factor and loading, thereby supporting componentwise fixed-parameter inference. We also extend the analysis to factor-augmented regressions. Under the conventional approximation $\hat{\bF}\approx\bF^{*}\hat{\bH}$, the coefficient vector associated with $\hat{\bF}$ depends on the data-dependent rotation $\hat{\bH}^{-1}$; under our approximation $\hat{\bF}\approx\bF^{0}$, it has a fixed population counterpart. We establish consistency and asymptotic normality of the least squares estimator under WF models, thereby providing a fixed-parameter interpretation of the usual $t$-tests on the coefficients associated with individual estimated PC factors.

We have also conducted extensive finite-sample experiments across different factor-strength regimes. The results show that approximation by the pseudo-true parameters generally compares favorably with the data-dependent alternatives, and that the corresponding normal approximations and $t$-tests have accurate finite-sample behavior. Similar conclusions hold for factor-augmented regressions.

Finally, we comment on the additional restrictions required to identify a particular starting representation $(\bF^{*},\bB^{*})$ as structural, a problem distinct from the fixed-target problem studied in this paper. As discussed in \cite{UY2019}, to directly identify a particular structural loading matrix $\bB^{*}$, $r^2$ (or more) constraints should be imposed along with the cross-sectional ordering of $x_{t,i}$. Such restrictions must be supplied exogenously. One way to identify structural parameters is to impose constraints guided by economic and financial theory; see discussions in \cite{StockWatson2016}. Another is to use exogenous shocks that result in structural breaks in the statistical model. Recently, \cite{YamamotoHara2022} proposed a method to identify factor-augmented regression models using changes in unconditional shock variances. Extending our approach to incorporate these methods would be a useful direction for future research.

\section*{Acknowledgment}
We are grateful to Silvia Gon\c{c}alves, Naoko Hara, Kazuhiko Hayakawa, and Yohei Yamamoto for helpful discussions and useful comments. 
This work was supported by JSPS KAKENHI (grant numbers 23H00804, 24K16343, 25K00625, 25K05036 and 25H00544).

\section*{Data Availability Statement}
The data used in the empirical illustration are publicly available in the replication materials for Stock and Watson (2016), available from 
\href{https://www.princeton.edu/~mwatson/publi.html}{Mark W. Watson’s website}.

\section*{Declaration of Generative AI Use}
During the preparation of this manuscript, the authors used ChatGPT (OpenAI) to improve English expression, readability, and the organization of the presentation. The authors reviewed and edited all generated content and take full responsibility for the content of the manuscript.


\bibliographystyle{chicago}
\bibliography{references_wfr}

\begin{thebibliography}{}

\bibitem[\protect\citeauthoryear{Ahn and Horenstein}{Ahn and Horenstein}{2013}]{AhnHorenstein2013}
Ahn, S.~C. and A.~R. Horenstein (2013).
\newblock Eigenvalue ratio test for the number of factors.
\newblock {\em Econometrica\/}~{\em 81}, 1203--1227.

\bibitem[\protect\citeauthoryear{Bai}{Bai}{2003}]{Bai2003}
Bai, J. (2003).
\newblock Inferential theory for factor models of large dimensions.
\newblock {\em Econometrica\/}~{\em 71}, 135--171.

\bibitem[\protect\citeauthoryear{Bai and Ng}{Bai and Ng}{2002}]{BaiNg2002}
Bai, J. and S.~Ng (2002).
\newblock Determining the number of factors in approximate factor models.
\newblock {\em Econometrica\/}~{\em 70}, 191--221.

\bibitem[\protect\citeauthoryear{Bai and Ng}{Bai and Ng}{2006}]{BaiNg2006}
Bai, J. and S.~Ng (2006).
\newblock Confidence intervals for diffusion index forecasts and inference with factor-augmented regressions.
\newblock {\em Econometrica\/}~{\em 74}, 1133--1150.

\bibitem[\protect\citeauthoryear{Bai and Ng}{Bai and Ng}{2008}]{BaiNg2008}
Bai, J. and S.~Ng (2008).
\newblock Forecasting economic time series using targeted predictors.
\newblock {\em Journal of Econometrics\/}~{\em 146\/}(2), 304--317.

\bibitem[\protect\citeauthoryear{Bai and Ng}{Bai and Ng}{2009}]{BaiNg2009}
Bai, J. and S.~Ng (2009).
\newblock Boosting diffusion indices.
\newblock {\em Journal of Applied Econometrics\/}~{\em 24\/}(4), 607--629.

\bibitem[\protect\citeauthoryear{Bai and Ng}{Bai and Ng}{2013}]{BaiNg2013}
Bai, J. and S.~Ng (2013).
\newblock Principal components estimation and identification of static factors.
\newblock {\em Journal of Econometrics\/}~{\em 176}, 18--29.

\bibitem[\protect\citeauthoryear{Bai and Ng}{Bai and Ng}{2023}]{BaiNg2023}
Bai, J. and S.~Ng (2023).
\newblock Approximate factor models with weaker loadings.
\newblock {\em Journal of Econometrics\/}~{\em 235\/}(2), 1893--1916.

\bibitem[\protect\citeauthoryear{Bailey, Kapetanios, and Pesaran}{Bailey et~al.}{2016}]{BaileyEtAl2016}
Bailey, N., G.~Kapetanios, and M.~H. Pesaran (2016).
\newblock Exponent of cross-sectional dependence: Estimation and inference.
\newblock {\em Journal of Applied Econometrics\/}~{\em 31}, 929--960.

\bibitem[\protect\citeauthoryear{Borenstein and Kellogg}{Borenstein and Kellogg}{2014}]{BorensteinKellogg2014}
Borenstein, S. and R.~Kellogg (2014).
\newblock The incidence of an oil glut: Who benefits from cheap crude oil in the midwest?
\newblock {\em The Energy Journal\/}~{\em 35\/}(1), 15--34.

\bibitem[\protect\citeauthoryear{Chamberlain and Rothschild}{Chamberlain and Rothschild}{1983}]{CR1983}
Chamberlain, G. and M.~Rothschild (1983).
\newblock Arbitrage, factor structure and mean-variance analysis in large asset markets.
\newblock {\em Econometrica\/}~{\em 51}, 1281--1304.

\bibitem[\protect\citeauthoryear{Cheng and Hansen}{Cheng and Hansen}{2015}]{ChengHansen2015}
Cheng, X. and B.~E. Hansen (2015).
\newblock Forecasting with factor-augmented regression: A frequentist model averaging approach.
\newblock {\em Journal of Econometrics\/}~{\em 186\/}(2), 280--293.

\bibitem[\protect\citeauthoryear{Choi and Yuan}{Choi and Yuan}{2025}]{ChoiYuan2024}
Choi, J. and M.~Yuan (2025).
\newblock High dimensional factor analysis with weak factors.
\newblock {\em Journal of Econometrics\/}~{\em 252}, 106086.

\bibitem[\protect\citeauthoryear{Connor and Korajczyk}{Connor and Korajczyk}{1986}]{ConnorKorajczyk1986}
Connor, G. and R.~A. Korajczyk (1986).
\newblock Performance measurement with the arbitrage pricing theory: A new framework for analysis.
\newblock {\em Journal of Financial Economics\/}~{\em 15}, 373--394.

\bibitem[\protect\citeauthoryear{Connor and Korajczyk}{Connor and Korajczyk}{1993}]{ConnorKorajczyk1993}
Connor, G. and R.~A. Korajczyk (1993).
\newblock A test for the number of factors in an approximate factor model.
\newblock {\em Journal of Finance\/}~{\em 48}, 1263--1291.

\bibitem[\protect\citeauthoryear{De~Mol, Giannone, and Reichlin}{De~Mol et~al.}{2008}]{DeMol2008}
De~Mol, C., D.~Giannone, and L.~Reichlin (2008).
\newblock Forecasting using a large number of predictors: Is {B}ayesian shrinkage a valid alternative to principal components?
\newblock {\em Journal of Econometrics\/}~{\em 146}, 318--328.

\bibitem[\protect\citeauthoryear{Fan, Liao, and Mincheva}{Fan et~al.}{2013}]{FanEtAl2013}
Fan, J., Y.~Liao, and M.~Mincheva (2013).
\newblock Large covariance estimation by thresholding principal orthogonal complements.
\newblock {\em Journal of the Royal Statistical Society Series B\/}~{\em 75}, 603--680.

\bibitem[\protect\citeauthoryear{Fan, Yan, and Zheng}{Fan et~al.}{2024}]{FanEtAl2024}
Fan, J., Y.~Yan, and Y.~Zheng (2024).
\newblock When can weak latent factors be statistically inferred?
\newblock {\em arXiv:2407.03616\/}.

\bibitem[\protect\citeauthoryear{Freyaldenhoven}{Freyaldenhoven}{2022}]{Freyaldenhoven21JoE}
Freyaldenhoven (2022).
\newblock Factor models with local factors - determining the number of relevant factors.
\newblock {\em Journal of Econometrics\/}~{\em 229}, 80--102.

\bibitem[\protect\citeauthoryear{Gonçalves and Perron}{Gonçalves and Perron}{2014}]{GoncalvesPerron2014}
Gonçalves, S. and B.~Perron (2014).
\newblock Bootstrapping factor-augmented regression models.
\newblock {\em Journal of Econometrics\/}~{\em 182\/}(1), 156--173.

\bibitem[\protect\citeauthoryear{Javanmard and Javadi}{Javanmard and Javadi}{2019}]{JJ2019}
Javanmard, A. and H.~Javadi (2019).
\newblock False discovery rate control via debiased lasso.
\newblock {\em Electronic Journal of Statistics\/}~{\em 13}, 1212--1253.

\bibitem[\protect\citeauthoryear{Jiang, Uematsu, and Yamagata}{Jiang et~al.}{2026}]{jiang2024Mw}
Jiang, P., Y.~Uematsu, and T.~Yamagata (2026).
\newblock Bias correction in factor-augmented regression models with weak factors.
\newblock {\em Journal of Business and Economic Statistics\/}.

\bibitem[\protect\citeauthoryear{Ludvigson and Ng}{Ludvigson and Ng}{2009}]{LudvigsonNg2009}
Ludvigson, C.~S. and S.~Ng (2009).
\newblock Macro factors in bond risk premia.
\newblock {\em Review of Financial Studies\/}~{\em 22}, 5027--5067.

\bibitem[\protect\citeauthoryear{Newey and West}{Newey and West}{1987}]{NeweyWest1987}
Newey, W.~K. and K.~D. West (1987).
\newblock A simple, positive semi-definite, heteroskedasticity and autocorrelation consistent covariance matrix.
\newblock {\em Econometrica\/}~{\em 55}, 703--708.

\bibitem[\protect\citeauthoryear{Onatski}{Onatski}{2010}]{Onatski2010}
Onatski, A. (2010).
\newblock Determining the number of factors from empirical distribution of eigenvalues.
\newblock {\em Review of Economics and Statistics\/}~{\em 92}, 1004--1016.

\bibitem[\protect\citeauthoryear{Stock and Watson}{Stock and Watson}{2016}]{StockWatson2016}
Stock, J. and M.~Watson (2016).
\newblock Chapter 8 - dynamic factor models, factor-augmented vector autoregressions, and structural vector autoregressions in macroeconomics.
\newblock Volume~2 of {\em Handbook of Macroeconomics}, pp.\  415--525. Elsevier.

\bibitem[\protect\citeauthoryear{Stock and Watson}{Stock and Watson}{2002a}]{StockWatson2002JASA}
Stock, J.~H. and M.~W. Watson (2002a).
\newblock Forecasting using principal components from a large number of predictors.
\newblock {\em Journal of the American Statistical Association\/}~{\em 97}, 1167--1179.

\bibitem[\protect\citeauthoryear{Stock and Watson}{Stock and Watson}{2002b}]{StockWatson2002JBES}
Stock, J.~H. and M.~W. Watson (2002b).
\newblock Macroeconomic forecasting using diffusion indexes.
\newblock {\em Journal of Business \& Economic Statistics\/}~{\em 30}, 147--162.

\bibitem[\protect\citeauthoryear{Uematsu and Yamagata}{Uematsu and Yamagata}{2023a}]{UY2019}
Uematsu, Y. and T.~Yamagata (2023a).
\newblock Estimation of sparsity-induced weak factor models.
\newblock {\em Journal of Business \& Economic Statistics\/}~{\em 41}, 213--227.

\bibitem[\protect\citeauthoryear{Uematsu and Yamagata}{Uematsu and Yamagata}{2023b}]{UY2019inference}
Uematsu, Y. and T.~Yamagata (2023b).
\newblock Inference in sparsity-induced weak factor models.
\newblock {\em Journal of Business \& Economic Statistics\/}~{\em 41}, 126--139.

\bibitem[\protect\citeauthoryear{{U.S. Energy Information Administration}}{{U.S. Energy Information Administration}}{2013}]{EIAspread2013}
{U.S. Energy Information Administration} (2013, August).
\newblock Spread narrows between {Brent} and {WTI} crude oil benchmark prices.
\newblock Today in Energy.
\newblock \url{https://www.eia.gov/todayinenergy/detail.php?id=12391}.

\bibitem[\protect\citeauthoryear{Wang and Fan}{Wang and Fan}{2017}]{WangFan2017}
Wang, W. and J.~Fan (2017).
\newblock {Asymptotics of empirical eigenstructure for high dimensional spiked covariance}.
\newblock {\em The Annals of Statistics\/}~{\em 45\/}(3), 1342 -- 1374.

\bibitem[\protect\citeauthoryear{Wei and Zhang}{Wei and Zhang}{2026}]{WeiZhang2026}
Wei, J. and Y.~Zhang (2026).
\newblock Can principal component analysis preserve the sparsity in factor loadings?
\newblock {\em Econometric Theory\/}, 1–24.

\bibitem[\protect\citeauthoryear{Yamamoto and Hara}{Yamamoto and Hara}{2022}]{YamamotoHara2022}
Yamamoto, Y. and N.~Hara (2022).
\newblock Identifying factor-augmented vector autoregression models via changes in shock variances.
\newblock {\em Journal of Applied Econometrics\/}~{\em 37\/}(4), 722--745.

\end{thebibliography}


\newpage
\appendix
\setcounter{page}{1}
\setcounter{section}{0}

\phantomsection\label{supplementary-material}
\bigskip

\begin{center}

{\large\bf SUPPLEMENTARY MATERIAL}

\end{center}

\section{Proofs of the Main Results}
\label{subsec:main results}
\setcounter{lem}{0}
\renewcommand{\thelem}{A.\arabic{lem}}
\renewcommand{\theequation}{A.\arabic{equation}}
\setcounter{equation}{0}


\subsection{Proof of Lemma \ref{lem:UV}}

\begin{proof}
	The eigen-decomposition of ${\bB^*}'\bB^*\left(T^{-1}{\bF^*}'\bF^*\right)$ yields
	\begin{align*}
	\bU\bV=\bLambda=\bLambda'=\bV'\bU'=\bV\bU.
	\end{align*}
	Since $\bU$, $\bV$, and $\bLambda$ are invertible under Assumption \ref{cond:eigen}, we obtain
	\begin{align}\label{UV}
	\bLambda^{-1}\bU=\bV^{-1}=\bU\bLambda^{-1}
	~~~\text{and}~~~
	\bV\bLambda^{-1}=\bU^{-1}=\bLambda^{-1}\bV.
	\end{align}
	From the first equation in \eqref{UV}, we have $(\bLambda^{-1}\bU)_{ij}=(\bU\bLambda^{-1})_{ij}$ for any $i,j\in\{1,\dots,r:i\not=j\}$, which entails $(1/\lambda_i-1/\lambda_j)u_{ij}=0$. Because $\lambda_i$'s are distinct and bounded away from zero by Assumption \ref{cond:eigen}, it must reduce $u_{ij}=0$; that is, $\bU$ is a diagonal matrix. In the same way, $\bV$ is also a diagonal matrix from the second equation in \eqref{UV}. This completes the proof. 
\end{proof}

\subsection{Proof of Theorem \ref{thm:rotation}}
\begin{proof}
For $\bF^0=\bF^*\bH$ and $\bB^0=\bB^*{\bH^{-1}}'$ with $\bH=\bP\bV^{-1/2}$, Lemma \ref{lem:UV} gives
	\begin{align*}
	T^{-1}{\bF^0}'\bF^0 
	&= \bV^{-1/2}\bP'\left(T^{-1}{\bF^*}'\bF^* \right)\bP\bV^{-1/2}
	=\bV^{-1/2}\bV\bV^{-1/2}
	=\bI
	\end{align*}
	and 
	\begin{align*}
	{\bB^0}'\bB^0
	= \bV^{1/2}\bP^{-1}{\bB^*}'\bB^*{\bP^{-1}}'\bV^{1/2}
	= \bV^{1/2}\bU\bV^{1/2}
	= \bU\bV
	= \bLambda~~(\text{diagonal}).
	\end{align*}
	Therefore \eqref{pc1} holds.
	
Finally, we prove the uniqueness of $\bH$. Suppose there exists another rotation matrix $\bar{\bH}$ such that $T^{-1}\bar{\bH}'{\bF^*}'\bF^*\bar{\bH}=\bI$ and $\bar{\bH}^{-1}{\bB^*}'\bB^* {\bar{\bH}^{-1\prime}}$ become diagonal (with elements ordered decreasingly). Then we must have
	\begin{align*}
	{\bar{\bH}^{-1\prime}}\bar{\bH}^{-1} = T^{-1}{\bF^*}'\bF^*
	= {\bP^{-1}}'\bV\bP^{-1} 
	= {\bH'}^{-1}{\bH}^{-1},
	\end{align*}
which is equivalent to $\bar\bH\bar\bH'=\bH\bH'$. Thus it is written as  $\bar{\bH}=\bH\bar{\bQ}$ for some orthogonal matrix $\bar{\bQ}$. We further have
	\begin{align*}
	\bar{\bH}^{-1}{\bB^*}'\bB^* {\bar{\bH}^{-1\prime}}
	= \bar{\bQ}'\bH^{-1}{\bB^*}'\bB^*{\bH^{-1}}'\bar{\bQ}
	= \bar{\bQ}'\bLambda\bar{\bQ},
	\end{align*}
which must be diagonalized. Since $\bLambda$ is a diagonal matrix with distinct elements in descending order and $\bar{\bQ}$ is orthogonal, the only choice of $\bar{\bQ}$ is $\bar{\bQ}=\bar{\bS}$, where $\bar{\bS}$ is a sign matrix. Hence, the construction of $\bH$ is unique up to sign changes. This completes the proof.
\end{proof}

\subsection{Proof of Theorem \ref{thm:invariance}}
\begin{proof}
Assumption~\ref{cond:eigen}(i) implies that
$\bF^{\ast}$, $\bB^{\ast}$, $\bF^{\dagger}$, and
$\bB^{\dagger}$ all have column rank $r$. Let
$\bC=\bF^{\ast}\bB^{\ast\prime}
=\bF^{\dagger}\bB^{\dagger\prime}$.
Since $\bB^{\ast\prime}$ and $\bB^{\dagger\prime}$ have row rank $r$, we have
$\operatorname{col}(\bC)
=\operatorname{col}(\bF^{\ast})
=\operatorname{col}(\bF^{\dagger})$. 
Therefore, there exists an $r\times r$ matrix $\bG$ such that
$\bF^{\dagger}=\bF^{\ast}\bG$. This matrix is unique because
$\bF^{\ast}$ has full column rank; explicitly, $\bG=
(\bF^{\ast\prime}\bF^{\ast})^{-1}
\bF^{\ast\prime}\bF^{\dagger}$. Moreover, $\bG$ is nonsingular. Indeed, if $\bG\bv=\bzero$, then
$\bF^{\dagger}\bv=\bF^{\ast}\bG\bv=\bzero$, and the full column rank
of $\bF^{\dagger}$ implies $\bv=\bzero$.

Substituting $\bF^{\dagger}=\bF^{\ast}\bG$ into
$\bF^{\dagger}\bB^{\dagger\prime}
=\bF^{\ast}\bB^{\ast\prime}$ and using the full column rank of
$\bF^{\ast}$ give
$\bB^{\ast\prime}=\bG\bB^{\dagger\prime}$, and hence
$\bB^{\dagger}=\bB^{\ast}\bG^{\prime-1}$.

Now set $\bK=\bG^{-1}\bH$. Then we obtain $\bF^{\dagger}\bK
=\bF^{\ast}\bH
=\bF^{0}$ and $\bB^{\dagger}\bK^{\prime-1}=
\bB^{\ast}\bH^{\prime-1}=\bB^{0}$. 
Thus, $\bK$ maps $(\bF^{\dagger},\bB^{\dagger})$ to the pseudo-true
parameter $(\bF^{0},\bB^{0})$ satisfying \eqref{pc1}. By the uniqueness
result in Theorem~\ref{thm:rotation}, there exists a diagonal sign
matrix $\bar{\bS}$ such that
$\bH^{\dagger}=\bK\bar{\bS}$. It follows that $\bF^{\dagger}\bH^{\dagger}
=\bF^{0}\bar{\bS}$ and $\bB^{\dagger}\bH^{\dagger\prime-1}=
\bB^{0}\bar{\bS}$. Hence, $\bF^{\dagger}\bH^{\dagger}=\bF^{0}$ and
$\bB^{\dagger}\bH^{\dagger\prime-1}=\bB^{0}$ hold up to simultaneous
column sign changes.
\end{proof}

\subsection{Proof of Lemma \ref{lem fh4bq}}
\begin{proof}(i) Since $\hat{\bLambda}$ and $\hat{\bF}$ are the eigenvalues and eigenvectors of $\bX\bX'/T$, respectively, we have
\begin{align*} 
T^{-1}\bX\bX'\hat{\bF}
= \hat{\bF}\hat{\bLambda}.
\end{align*}
By Theorem \ref{thm:rotation}, we can find a rotation matrix $\bH$ such that the rotated parameters $\bF^0=\bF^* {\bH}$ and $\bB^0=\bB^* {\bH}^{\prime -1}$ satisfy the $r^2$ restrictions,
\begin{align*}
\frac{1}{T}\bF^{0 \prime}\bF^0 = \bI_r
~~\text{ and }~~
\bB^{0 \prime}\bB^0 \in \cD(r).
\end{align*}
Meanwhile, using the model $\mathbf{X}=\bF^0 \bB^{0 \prime}+\mathbf{E}$, $\bX\bX'$ can be expanded as
\begin{align*}\label{XX}
\bX\bX'
= \bF^0{\bB^0}'\bB^0{\bF^0}' 
+ \bF^0{\bB^0}'\bE'
+ \bE\bB^0{\bF^0}'+\bE\bE'.
\end{align*}
Thus, we obtain
\[
\hat{\bF}=\frac{1}{T}\left( \bF^0{\bB^0}'\bB^0{\bF^0}' 
+ \bF^0{\bB^0}'\bE'
+ \bE\bB^0{\bF^0}'+\bE\bE' \right) \hat{\bF}\hat{\bLambda}^{-1}.
\]	
From $\hat{\bB}= \bX^{\prime} \hat{\bF}/T$, we can see that
 \begin{align}
 \nonumber
\hat{\bB}&=\frac{1}{T} \bB^0 {\bF^0}^{\prime} \hat{\bF}+\frac{1}{T} \bE^{\prime} \hat{\bF}, \\\nonumber
{\bB^0}' \hat{\bB}&=\frac{1}{T} {\bB^0}' \bB^0 {\bF^0}' \hat{\bF}+\frac{1}{T} {\bB^0}' \bE^{\prime} \hat{\bF}, \\\nonumber
 {\bB^0}' \hat{\bB} \hat{\bLambda}^{-1}&=\frac{1}{T} {\bB^0}' \bB^0 {\bF^0}' \hat{\bF}\hat{\bLambda}^{-1}+\frac{1}{T} {\bB^0}' \bE^{\prime} \hat{\bF} \hat{\bLambda}^{-1}, \\
 \label{eq:H4-H}
\tilde{\bH}_4&=\tilde{\bH}+\frac{1}{T} {\bB^0}' \bE^{\prime} \hat{\bF} \hat{\bLambda}^{-1},
   \end{align}
where $\tilde{\bH}=\frac{1}{T} {\bB^0}' \bB^0 {\bF^0}' \hat{\bF}\hat{\bLambda}^{-1}$. Then, we have
\begin{eqnarray*}
	\hat{\bF} - \bF^* \hat{\bH}_4=\hat{\bF} - \bF^0 \tilde{\bH}_4 =\left( \frac{1}{T} \bE\bE'\hat{\bF}+\frac{1}{T} \bE\bB^0{\bF^0 }'\hat{\bF} \right) \hat{\bLambda}^{-1}.
\end{eqnarray*}
Lemma \ref{lem:BEF}(iii) yields $ \left\| {\bB^0}'\bE'\right\|_{\F} =O_p(\sqrt{TN^{\alpha_{1}}})$. Thus the inequality $\left\| \bA\bB \right\|_{\F} \le \left\| \bA \right\|_2\left\| \bB \right\|_{\F}$ gives
\begin{eqnarray*}
    \left\| \bE\bE'\hat{\bF} \right\|_{\F}
    \le \lambda_{1} [\bE\bE'] \|\hat{\bF} \|_{\F} = O_p\left( N+T\right) \sqrt{T}.
\end{eqnarray*}
Therefore, we obtain the first result:
\begin{align*}
	& \frac{1}{\sqrt{T}}\left\| 	\hat{\bF} - \bF^* \hat{\bH}_4\right\|_{\F} \\
 &\le 
	\left\|  \frac{1}{T^{3/2}} \bE\bE'\hat{\bF}\hat{\bLambda}^{-1} \right\|_{\F}+\left\|\frac{1}{T^{3/2}} \bE\bB^0{\bF^0 }'\hat{\bF}\hat{\bLambda}^{-1} \right\|_{\F} \\
 &\le  \frac{\left\| \bE\bE'\hat{\bF} \right\|_{\F} }{T^{3/2}} \left\| \hat{\bLambda}^{-1} \right\|_{\F} +\left\|\frac{1}{T^{3/2}} \left( \bE\bB^0 \bN^{-\frac{1}{2}}\right) \left( \bN^{\frac{1}{2}}{\bF^0 }'\hat{\bF}\bN^{-\frac{1}{2}}\right) \bN^{\frac{1}{2}}\hat{\bLambda}^{-1} \right\|_{\F} \\
 & = O_p\left(
 \left( \frac{N}{T}+1\right) N^{-\alpha_r}\right)
+O_p\left(N^{-\frac{1}{2}\alpha_r}\right) \\
 & = O_p\left(\frac{N^{1-\alpha_r}}{T}\right)+O_p\left(N^{-\frac{1}{2}\alpha_r}\right),
\end{align*}
where we have used Lemmas \ref{lem:BEF} and \ref{lemma order of lambda}.\\
(ii) By the definition, 
	\begin{eqnarray*}
		\hat{\bB} &=& \frac{1}{T} \bX'\hat{\bF} 
		=\frac{1}{T}  \bB^0{\bF^0}'\hat{\bF} +\frac{1}{T} \bE'\hat{\bF}
	\end{eqnarray*}
	implies 
	\begin{align}
\hat{\bB}- \bB^0\tilde{\bQ}' = \frac{1}{T} \bE'(\hat{\bF}- \bF^0 \tilde{\bH}_4)+ \frac{1}{T} \bE'\bF^0 \tilde{\bH}_4, \label{B}
	\end{align}	
 where $\tilde{\bQ}'=T^{-1}{\bF^0}'\hat{\bF}$. 
	Using $\tilde{\bH}_4=O_p(1)$ in Lemma \ref{lem fh4bq}(iii), Lemma \ref{lemma f(fh-f0)}(i) and
	\begin{eqnarray*}
		\left\|  \frac{1}{T} \bE'\bF^0 \tilde{\bH}_4\right\|_{\F}  &=&  O_p\left(\frac{1}{T} \sqrt{NT}  \right)= O_p\left(\sqrt{\frac{N}{T}}\right) ,
	\end{eqnarray*}
	we have	
	\begin{align*}
		&\frac{1}{\sqrt{N}} \left\| \hat{\bB} - \bB^0 \tilde{\bQ}' \right\|_{\F}\\
		& = \frac{1}{\sqrt{N}} \left\|\frac{1}{T} \bE'(\hat{\bF}- \bF^0 \tilde{\bH}_4)+ \frac{1}{T} \bE'\bF^0 \tilde{\bH}_4 \right\|_F\\
  &\le \frac{1}{\sqrt{N}} \left[  O_p\left(\sqrt{\frac{N^{1-\alpha_r}}{T}}\frac{N^{1-\frac{1}{2}\alpha_r}}{T}\right)+ O_p\left(\frac{N^{1-\frac{1}{2}\alpha_r}}{T}\right)+ O_p\left(N^{-\frac{1}{2}\alpha_{r}}\right) \right] +   \frac{1}{\sqrt{N}}O_p\left(\sqrt{\frac{N}{T}} \right)\\
		&= O_p\left(\frac{1}{\sqrt{T}}\right)+O_p\left(N^{-\frac{1}{2}-\frac{1}{2}\alpha_{r}}\right),
	\end{align*}
if $\frac{N^{1-\alpha_r}}{T}\to 0$.

(iii) Throughout the proof, the column signs of \(\hat \bF\) and \(\hat \bB\) are chosen so that the diagonal elements of \(\tilde \bQ=T^{-1}\hat \bF^\prime \bF^0\) are nonnegative. Without this convention, the limits in Lemma 2(iii) hold up to a diagonal sign matrix. Theorem \ref{thm:rotation} is crucial to construct the rotation matrices with ``tilde" and derive their boundedness.  
Let $\tilde{\bH}_2= ({\bF^0}'\bF^0)^{-1}{\bF^0}'\hat{\bF}={\bF^0}'\hat{\bF}/T=\tilde{\bQ}'$. 
Since both \(T^{-1}{\bF^0}'\bF^0=\bI_r\) and
\(T^{-1}\hat{\bF}'\hat{\bF}=\bI_r\) hold, the matrices
\(\bF^0/\sqrt T\) and \(\hat{\bF}/\sqrt T\) have orthonormal columns.
Therefore, $||\tilde{\bH}_2||_2\le || \frac{1}{\sqrt{T}} \hat{\bF}||_2 || \frac{1}{\sqrt{T}}\bF^0||_2 =O_p(1)$.  We now establish a rough bound: \[ \|\tilde{\bH}_2-\tilde{\bH}_4\|_{\F} = O_p\left( \sqrt{\frac{N^{1-\alpha_r}}{T}} \right). \]
By Lemma \ref{lem:H equality}(ii)
\begin{align}
    \|\tilde{\bH}_4-\tilde{\bH}_2\|_F
    =
   \|T^{-1}{\bF^0}'\bE\hat{\bB}\hat{\bLambda}^{-1}\|_{\F}
   \le \frac{1}{T} \|{\bF^0}'\bE\|_{\F} \|\hat{\bB}\hat{\bLambda}^{-1}\|_{\F}.
\end{align}
Moreover, 
\begin{align} \|\hat{\bB}\hat{\bLambda}^{-1}\|_{\F}^2 &= \operatorname{tr} \left( \hat{\bLambda}^{-1} \hat{\bB}'\hat{\bB} \hat{\bLambda}^{-1} \right) \\ &= \operatorname{tr}(\hat{\bLambda}^{-1}) = \sum_{k=1}^r \hat{\lambda}_k^{-1}. 
\end{align} 
Using \(\hat{\lambda}_k\asymp_p N^{\alpha_k}\) for \(k=1,\ldots,r\) in Lemma \ref{lemma order of lambda}, and noting that \(r\) is fixed, it follows that $ \|\hat{\bB}\hat{\bLambda}^{-1}\|_{\F} = O_p(N^{-\alpha_r/2})$.  Combining the above bounds gives \[ 
\|\tilde{\bH}_2-\tilde{\bH}_4\|_{\F} 
\le \frac{1}{T} O_p(\sqrt{NT}) O_p(N^{-\alpha_r/2}) 
= O_p\left( \sqrt{\frac{N^{1-\alpha_r}}{T}} \right).  \] 
Hence, if \(N^{1-\alpha_r}/T\to 0\), then $\|\tilde{\bH}_2-\tilde{\bH}_4\|_{\F} = o_p(1)$.

Next, we show the probability limit of the rotation matrices assuming that  $\frac{N^{1-\alpha_r}}{T} \to 0$. Given Proposition \ref{prop:Q-I dec}, we obtain $\tilde{\bQ} \CP \bI_r$. The preliminary bound established above gives $ \|\tilde{\bQ}'-\tilde{\bH}_4\|_{\F} = \|\tilde{\bH}_2-\tilde{\bH}_4\|_{\F} = o_p(1)$. Thus, we complete the proof.
\end{proof}

\subsection{Proof of Theorem \ref{0}}

\begin{proof}(i) 
By the definition of $\hat{\bB}$ and expanding $\bX$, we obtain
\begin{align*}
   \hat{\bF}=\frac{1}{T}\bX\bX'\hat{\bF}\hat{\bLambda}^{-1}=\bX\hat{\bB}\hat{\bLambda}^{-1}
=\bF^0{\bB^0}'\hat{\bB}\hat{\bLambda}^{-1}+\bE\hat{\bB}\hat{\bLambda}^{-1}. 
\end{align*}
Then,
$ 
\frac{1}{T} {\bF^0}'\hat{\bF}=\tilde{\bH}_4+ \frac{1}{T} {\bF^0}'\bE\hat{\bB}\hat{\bLambda}^{-1} . 
$ Under the conditions $\frac{N^{1-\alpha_r}}{T}\to 0$, applying Lemmas \ref{lem fh4bq}(i), \ref{lem:FEB0}(i), and Proposition \ref{prop:Q-I dec}, we show the first result:
\begin{align*}
  & \frac{1}{\sqrt{T}} \left\|  \hat{\bF} - \bF^0\right\|_{\F} \\
	&\le  \frac{1}{\sqrt{T}} \left\|  \hat{\bF} - \bF^0\tilde{\bH}_2\right\|_{\F} +\frac{1}{\sqrt{T}} \left\| \bF^0(\tilde{\bH}_2-\bI_r) \right\|_{\F} \\
 & = \frac{1}{\sqrt{T}} \left\|\frac{1}{T} \bE\bE'\hat{\bF}\hat{\bLambda}^{-1} +\frac{1}{T} \bE\bB^0{\bF^0 }'\hat{\bF} \hat{\bLambda}^{-1}- \frac{1}{T}\bF^0 {\bF^0}'\bE\hat{\bB}\hat{\bLambda}^{-1} \right\|_{\F} + O_p\left( \tilde{\Delta}_{NT}\right) \\ 
 &\le O_p\left(\left( \frac{N}{T} +1\right) N^{ - \alpha_{r}}  \right)+O_p\left(N^{ -\frac{1}{2}\alpha_{r}}\right) + O_p\left(\Delta_1\right)+O_p\left(\tilde{\Delta}_{NT} \right) \\
 &= O_p\left(\left( \frac{N}{T} +1\right) N^{ - \alpha_{r}} \right) +O_p\left(N^{ -\frac{1}{2}\alpha_{r}}\right) 
+ O_p\left(\sqrt{\frac{N^{1-\alpha_r}}{T} }N^{-\alpha_r} +\frac{N^{1-\alpha_r}}{T}+\frac{1}{\sqrt{TN^{\alpha_r}}}\right) \\
 & \quad \quad +O_p\left(N^{-\alpha_r} +\frac{N^{1-\alpha_r}}{T} \right) \\
 & =  O_p\left(\frac{N^{1-\alpha_r}}{T}\right)+ O_p\left( N^{-\frac{1}{2}\alpha_{r}}\right),
\end{align*}
 where $\Delta_{1}$ and $\tilde{\Delta}_{NT}$ are defined by \eqref{Delta_1} and \eqref{eq:TDleta}.

 (ii) 	Since \(\bB^{0\prime}\bB^0=\bLambda\), we have the exact identity
	\[
	\left\|\bB^0(\tilde{\bQ}'-\bI_r)\right\|_{\F}^2
	=
	\operatorname{tr}\{(\tilde{\bQ}'-\bI_r)'\bB^{0\prime}\bB^0(\tilde{\bQ}'-\bI_r)\}
	=
	\left\|\bLambda^{1/2}(\tilde{\bQ}'-\bI_r)\right\|_{\F}^2.
	\] 
	Recall that
	\[
	\hat{\bB}-\bB^0 = (\hat{\bB}-\bB^0\tilde{\bQ}') + \bB^0(\tilde{\bQ}'-\bI_r).
	\]
	Therefore,
	\[
	\frac{1}{\sqrt N}\|\hat{\bB}-\bB^0\|_{\F} \le \frac{1}{\sqrt N} \|\hat{\bB}-\bB^0\tilde{\bQ}'\|_{\F} + \frac{1}{\sqrt N} \|\bB^0(\tilde{\bQ}'-\bI_r)\|_{\F}.
	\]
	By Lemma \ref{lem fh4bq}(ii),
	\[
	\frac{1}{\sqrt N} \|\hat{\bB}-\bB^0\tilde{\bQ}'\|_{\F} = O_p\left(\frac{1}{\sqrt T}\right) + O_p\left(N^{-\frac{1}{2}-\frac{1}{2}\alpha_r}\right).
	\]
	By Proposition \ref{prop:Q-I dec}(iii),
	\[
	\frac{1}{\sqrt N} \|\bB^0(\tilde{\bQ}'-\bI_r)\|_{\F}  =\frac{1}{\sqrt N}
	\left\|\bLambda^{1/2}(\tilde{\bQ}'-\bI_r)\right\|_{\F}=o_p\left(\frac{1}{\sqrt T}\right) + O_p\left(N^{-\frac{1}{2}-\frac{1}{2}\alpha_r}\right).
	\]
	Combining the two bounds gives
	\[
	\frac{1}{\sqrt N}\|\hat{\bB}-\bB^0\|_{\F} = O_p\left(\frac{1}{\sqrt T}\right) + O_p\left(N^{-\frac{1}{2}-\frac{1}{2}\alpha_r}\right).
	\]
    
(iii)
Applying Lemmas \ref{lem fh4bq} and \ref{lem:H equality}(vi), we obtain
\begin{align*}
\frac{1}{\sqrt{NT}}\|\hat{\bC}-\bC^*\|_{\F}
&\le
\frac{1}{\sqrt{NT}}
\left\|
(\hat{\bF}-\bF^0\tilde{\bH}_4)\hat{\bB}'
\right\|_{\F}  \\
&\quad+
\frac{1}{\sqrt{NT}}
\left\|
\bF^0\tilde{\bH}_4
\left(
\hat{\bB}
-\bB^0\tilde{\bH}_2
+
\bB^0\tilde{\bH}_2
-
\bB^0\tilde{\bH}_4^{'-1}
\right)'
\right\|_{\F}.
\end{align*}
We bound the two terms separately. First, since \(\hat{\bB}'\hat{\bB}=\hat{\bLambda}\),
\[
\|(\hat{\bF}-\bF^0\tilde{\bH}_4)\hat{\bB}'\|_{\F}^2
=
\operatorname{tr}((\hat{\bF}-\bF^0\tilde{\bH}_4)\hat{\bB}'\hat{\bB}(\hat{\bF}-\bF^0\tilde{\bH}_4)')
=
\|(\hat{\bF}-\bF^0\tilde{\bH}_4)\hat{\bLambda}^{1/2}\|_{\F}^2.
\]
Hence
\[
\frac{1}{\sqrt{NT}}
\left\|
(\hat{\bF}-\bF^0\tilde{\bH}_4)\hat{\bB}'
\right\|_{\F}
=
\frac{1}{\sqrt{NT}}
\left\|
(\hat{\bF}-\bF^0\tilde{\bH}_4)\hat{\bLambda}^{1/2}
\right\|_{\F}.
\]
Using
\[
\hat{\bF}-\bF^0\tilde{\bH}_4
=
\left(
\frac1T\bE\bE'\hat{\bF}
+
\frac1T\bE\bB^0\bF^{0\prime}\hat{\bF}
\right)
\hat{\bLambda}^{-1},
\]
we have
\begin{align}
&
\frac{1}{\sqrt{NT}}
\left\|
(\hat{\bF}-\bF^0\tilde{\bH}_4)\hat{\bB}'
\right\|_{\F}  \notag\\
&\le
\frac{1}{\sqrt{NT}}
\left\|
\frac1T\bE\bE'\hat{\bF}\hat{\bLambda}^{-1/2}
\right\|_{\F}
+
\frac{1}{\sqrt{NT}}
\left\|
\bE\bB^0\tilde{\bQ}'\hat{\bLambda}^{-1/2}
\right\|_{\F},\\
&\le
\frac{1}{\sqrt{NT}}
\frac1T
\|\bE\|_2^2
\|\hat{\bF}\|_{\F}
\|\hat{\bLambda}^{-1/2}\|_2  
+ O_p(N^{-1/2})\notag\\
&=
O_p\left(\frac{N^{(1-\alpha_r)/2}}{T}\right)
+
O_p\left(N^{-1/2-\alpha_r/2}\right)   
+ O_p(N^{-1/2})
\notag\\
&=
O_p\left(\frac{N^{(1-\alpha_r)/2}}{T}\right)
+
O_p(N^{-1/2}).
\label{eq:C-factor-EE}
\end{align}
Next, for the second term, Lemma \ref{lem:H equality}(vi) gives
\begin{align*}
&\frac{1}{\sqrt{NT}}
\left\|
\bF^0\tilde{\bH}_4
\left(
\hat{\bB}
-\bB^0\tilde{\bH}_2
+
\bB^0\tilde{\bH}_2
-
\bB^0\tilde{\bH}_4^{'-1}
\right)'
\right\|_{\F}  \\
&\quad\le
\frac{1}{\sqrt{NT}}
\|\bF^0\tilde{\bH}_4\|_{\F}
\left\|
\frac1T\bE'\hat{\bF}
-
\bB^0
\frac1T
(\hat{\bB}'\bB^0)^{-1}
\hat{\bB}'\bE'\hat{\bF}
\right\|_{\F}  \\
&\quad\le
\frac{1}{\sqrt{NT}}
\|\bF^0\tilde{\bH}_4\|_{\F}
\left\|
\frac1T\bE'\hat{\bF}
\right\|_{\F}.
\end{align*}
Since \(\|\bF^0\tilde{\bH}_4\|_{\F}=O_p(\sqrt T)\) and
\[
\left\|
\frac1T\bE'\hat{\bF}
\right\|_{\F}
=
O_p\left(\sqrt{\frac NT}\right)
+
O_p(N^{-\alpha_r/2}),
\]
we have
\begin{equation}
\frac{1}{\sqrt{NT}}
\|\bF^0\tilde{\bH}_4\|_{\F}
\left\|
\frac1T\bE'\hat{\bF}
\right\|_{\F}
=
O_p\left(\frac1{\sqrt T}\right)
+
O_p(N^{-1/2-\alpha_r/2}).
\label{eq:C-loading-side-final}
\end{equation}

Finally, combining \eqref{eq:C-factor-EE} and \eqref{eq:C-loading-side-final}, and using
\(\frac{N^{(1-\alpha_r)/2}}{T} \lesssim \frac{1}{\sqrt{T}}\), we obtain
\[
\frac{1}{\sqrt{NT}}
\|\hat{\bC}-\bC^*\|_{\F}
=
O_p(T^{-1/2})
+
O_p(N^{-1/2}).
\]
This completes the proof.

\end{proof}

\subsection{Proof of Lemma \ref{lem dist}}

\begin{proof}(i) We multiply $\hat{\bB}(\hat{\bB}'\hat{\bB})^{-1}$ to the both side of $\bX=\bF^0{\bB^{0}}'+\bE $, then, we have,
\begin{eqnarray*}
	\bX\hat{\bB}(\hat{\bB}'\hat{\bB})^{-1}&=&\bF^0{\bB^{0}}'\hat{\bB}(\hat{\bB}'\hat{\bB})^{-1}+\bE \hat{\bB}(\hat{\bB}'\hat{\bB})^{-1}\\
	\hat{\bF}&=& \bF^0\tilde{\bH}_4+\bE\bB^0\tilde{\bQ}'(\hat{\bB}'\hat{\bB})^{-1}+\bE(\hat{\bB}-\bB^0\tilde{\bQ}')(\hat{\bB}'\hat{\bB})^{-1}. 
\end{eqnarray*}
Since $\bF^*\hat{\bH}_4=\bF^0\tilde{\bH}_4$, the $t$-th row of $\hat{\bF}$ is
\begin{eqnarray*}
	\hat{\bff}_t'-  {\bff_t^*}' \hat{\bH}_4=\hat{\bff}_t'-  {\bff_t^0}' \tilde{\bH}_4 =\be_t'\bB^0\tilde{\bQ}'(\hat{\bB}'\hat{\bB})^{-1}+\be_t'(\hat{\bB}-\bB^0\tilde{\bQ}')(\hat{\bB}'\hat{\bB})^{-1}.
\end{eqnarray*}
That is
\begin{eqnarray*}
	\hat{\bff}_t- \hat{\bH}'_4 \bff_t^*=\hat{\bff}_t- \tilde{\bH}'_4 \bff_t^0 &=&(\hat{\bB}'\hat{\bB})^{-1}\tilde{\bQ}{\bB^0}'\be_t+(\hat{\bB}'\hat{\bB})^{-1} (\hat{\bB}-\bB^0\tilde{\bQ}')'\be_t.
\end{eqnarray*}
Then, we have  
\begin{align}
	\nonumber
	&\bD\bN^{\frac{1}{2}}(\hat{\bff}_t -  \hat{\bH}_4' \bff_t^*)\\
 &=\bD\bN^{\frac{1}{2}}(\hat{\bff}_t -  \tilde{\bH}_4' \bff_t^0) \\
	\label{ft-f0,0}
	&= \bD\bN^{\frac{1}{2}}(\hat{\bB}'\hat{\bB})^{-1} \tilde{\bQ} \bB^{0'}\be_t
	+ \bD\bN^{\frac{1}{2}}(\hat{\bB}'\hat{\bB})^{-1}\left( \hat{\bB} - \bB^0\tilde{\bQ}'\right) '\be_t .
\end{align}
We first consider the first term on the right-hand side of the above equation. By the order in Proposition \ref{prop:Q-I dec}(ii), $\bN^{-\frac{1}{2}}\hat{\bB}'\hat{\bB}\bN^{-\frac{1}{2}} - \bN^{-\frac{1}{2}}{\bB^0}'\bB^0\bN^{-\frac{1}{2}}= O_p(\tilde{\Delta}_{NT})=o_p(1)$ if $\frac{N^{1-\alpha_r}}{T}\to 0$. It implies that
\[
\bD^{-1}\bN^{-\frac{1}{2}}\hat{\bB}'\hat{\bB}\bN^{-\frac{1}{2}}\bD^{-1} \CP \bD^{-1}\bN^{-\frac{1}{2}}{\bB^0}'\bB^0\bN^{-\frac{1}{2}}\bD^{-1}= \bI_r. \] Now $\bD^{-1} \bN^{-\frac{1}{2}} \bB^{0'}\be_t \CD N(0, \bGamma_t)$ by Assumption \ref{assumption dist} and $\tilde{\bQ}-\bI_r=O_p(\tilde{\Delta}_{NT})$. The first term on the right-hand side of \eqref{ft-f0,0} is thus asymptotically normal, and we obtain
\[
\bD\bN^{\frac{1}{2}}(\hat{\bB}'\hat{\bB})^{-1} \tilde{\bQ} \bB^{0'}\be_t=\bD\bN^{\frac{1}{2}}(\hat{\bB}'\hat{\bB})^{-1}  \bB^{0'}\be_t(\bI_r +o_p(1)) ,
\]
whose asymptotic distribution is
\[
\bD\bN^{\frac{1}{2}}(\hat{\bB}'\hat{\bB})^{-1}\bN^{\frac{1}{2}}\bD \bD^{-1} \bN^{-\frac{1}{2}}\bB^{0'}\be_t \CD N(\mathbf{0}, \bGamma_t).
\]
Next, the second term on the right-hand side of \eqref{ft-f0,0} is $o_p(1)$ by Lemmas \ref{lemma bh-b0*e} if $ \alpha_r>1 / 2 $ and $\frac{N^{\frac{3}{2}-\alpha_r}}{T} \rightarrow 0
$. Collecting these results, we obtain (i):
\begin{eqnarray*}
	\bD\bN^{\frac{1}{2}}(\hat{\bff}_t -  \hat{\bH}_4' \bff_t^*)\CD N(\mathbf{0}, \bGamma_t),
\end{eqnarray*}
if $ \alpha_r>1 / 2 $ and $\frac{N^{\frac{3}{2}-\alpha_r}}{T} \rightarrow 0
$.\\
(ii) Recall \eqref{B}
\begin{eqnarray*}
	\hat{\bB}-\bB^*\hat{\bQ}'=\hat{\bB}-\bB^0\tilde{\bQ}'= \frac{1}{T} \bE'(\hat{\bF}-\bF^0\tilde{\bH}_4)+ \frac{1}{T} \bE'\bF^0\tilde{\bH}_4.
 \end{eqnarray*}
 The $i$-th row is
\begin{eqnarray*}
	&& \hat{\bb}_i'-   {\bb_i^*}'\hat{\bQ}'=\hat{\bb}_i'-   {\bb_i^0}'\tilde{\bQ}' =
 \frac{1}{T} \be_i'(\hat{\bF}-\bF^0\tilde{\bH}_4)+ \frac{1}{T} \be_i'\bF^0\tilde{\bH}_4 \\
 && \hat{\bb}_i- \hat{\bQ} {\bb_i^*} =\hat{\bb}_i- \tilde{\bQ} {\bb_i^0} =
 \frac{1}{T} (\hat{\bF}-\bF^0\tilde{\bH}_4)'\be_i+ \frac{1}{T} {\tilde{\bH}_4}'{\bF^0}'\be_i.
\end{eqnarray*}
We have
\begin{align}
 \sqrt{T}(\hat{\bb}_i -  \hat{\bQ} \bb_i^*) = \sqrt{T}(\hat{\bb}_i -  \tilde{\bQ} \bb_i^0) 
	\label{bt-b0,0}
	=  \frac{1}{\sqrt{T}} {\tilde{\bH}_4}'{\bF^0}'\be_i + \frac{1}{\sqrt{T}} (\hat{\bF}-\bF^0\tilde{\bH}_4)'\be_i.
\end{align}
By Lemmas \ref{lem:H equality}(ii), \ref{lem:FEB0}(i) and Proposition \ref{prop:Q-I dec}(i),
\begin{eqnarray*}
	  {\tilde{\bH}}'_4-\bI_r = {\tilde{\bH}}'_4-{\tilde{\bH}_2}'+{\tilde{\bH}}_2'-\bI_r =O_p(\Delta_{1} )+O_p(\tilde{\Delta}_{NT} )=o_p(1) ,  \; \text{if} \;  \frac{N^{1-\alpha_r}}{T}\to 0.
\end{eqnarray*}
Under Assumption \ref{assumption dist}, the first term on the right-hand side of \eqref{bt-b0,0} is thus asymptotically normal. That is 
\[
 \frac{1}{\sqrt{T}} {\tilde{\bH}_4}'{\bF^0}'\be_i=\frac{1}{\sqrt{T}}  {\bF^0}'\be_i+({\tilde{\bH}_4}'-\bI_r  )  \frac{1}{\sqrt{T}}  {\bF^0}'\be_i=\frac{1}{\sqrt{T}}  {\bF^0}'\be_i+O_p(\tilde{\Delta}_{NT}),
\]
whose asymptotic distribution is $N(\mathbf{0}, \bPhi_i)$. 

Next, the second term on the right-hand side of \eqref{bt-b0,0} is bounded by
\begin{eqnarray*}
    \left\| \frac{1}{\sqrt{T}}{\be_i}'(\hat{\bF}-\bF^0\tilde{\bH}_4) \right\|_{\F} 
		= O_p\left( \frac{N^{1-\alpha_r}}{\sqrt{T}}\right)+O_p\left(\sqrt{T}N^{-\alpha_r}\right)+O_p\left(N^{ - \frac{1}{2}\alpha_{r}}\right)  ,
\end{eqnarray*}
from Lemma \ref{lemma f(fh-f0)}(iv). It is negligible if $\frac{N^{1-\alpha_r}}{\sqrt{T}} \rightarrow 0$ and $ \frac{\sqrt{T}}{N^{\alpha_r}} \rightarrow 0$ are satisfied. Collecting these results, we obtain (ii)
\begin{align*}
    \sqrt{T}(\hat{\bb}_i-\hat{\bQ}\bb_i^*) \CD N(\mathbf{0}, \bPhi_i).  &
\end{align*}
if $\frac{N^{1-\alpha_{r}}}{\sqrt{T}} \to 0$, and $\frac{\sqrt{T}}{N^{\alpha_r}} \to 0$.
\end{proof}

\subsection{Proof of Theorem \ref{thm dist}}

\begin{proof}(i) We have
\begin{align}
	\nonumber
	\bD\bN^{\frac{1}{2}}(\hat{\bff}_t - \bff_t^0)
 =\bD\bN^{\frac{1}{2}}(\hat{\bff}_t -  \tilde{\bH}_4' \bff_t^0)+\bD\bN^{\frac{1}{2}}( \tilde{\bH}_4' - \bI_r) \bff_t^0.
\end{align}
We only consider the second term on the right-hand side of the above equation.  Lemmas \ref{lem:H equality}(ii), \ref{lem:FEB0}(i), and Proposition \ref{prop:Q-I dec} imply
\begin{align}
	 \bN^{\frac{1}{2}} ({\tilde{\bH}}'_4-\bI_r )=\bN^{\frac{1}{2}} \left( {\tilde{\bH}}'_4-{\tilde{\bH}_2}'+{\tilde{\bH}}_2'-\bI_r \right)
  \le N^{\frac{1}{2}\alpha_1}O_p(\Delta_{1} )+N^{\frac{1}{2}\alpha_1}O_p(\tilde{\Delta}_{NT} )=o_p(1), \label{eq:H4}
\end{align}
if $\frac{N^{\frac{3}{2}-\alpha_r}}{T} \rightarrow 0$, and $ N^{\frac{1}{2}\alpha_1-\alpha_r} \rightarrow 0$.

Using the result in Lemma \ref{lem dist}, we obtain
\begin{eqnarray*}
	\bD\bN^{\frac{1}{2}}(\hat{\bff}_t -   \bff_t^0)\CD N(\mathbf{0}, \bGamma_t),
\end{eqnarray*}
if $\alpha_r>\frac{1}{2}$ and $\frac{N^{\frac{3}{2}-\alpha_r}}{T} \rightarrow 0$.
\\
(ii) We have
\begin{align}
	\nonumber
	\sqrt{T} (\hat{\bb}_i - \bb_i^0)
 =\sqrt{T}(\hat{\bb}_i -  \tilde{\bQ} \bb_i^0)+\sqrt{T}( \tilde{\bQ} - \bI_r) \bb_i^0 .
\end{align}
The second term on the right-hand side of the above equation is bounded by    
\begin{eqnarray*}
    \left\| \sqrt{T}( \tilde{\bQ} - \bI_r) \bb_i^0 \right\|_{\F} 
		=\sqrt{T} O_p\left(
 \tilde{\Delta}_{NT}\right)=o_p(1),
\end{eqnarray*}
if $\sqrt{T} N^{-\alpha_r} \rightarrow 0$ and $\frac{N^{1-\alpha_r}}{\sqrt{T}} \rightarrow 0$ hold. Thus, we have
\begin{align*}
    \sqrt{T}(\hat{\bb}_i-\bb_i^0) \CD N(\mathbf{0}, \bPhi_i).  &
\end{align*} \\
(iii) By the definition of $\hat{c}_{t,i}$, we have
\begin{align}
\nonumber
	&\hat{c}_{t,i} - c_{t,i}^* \\
 &= \hat{\bb}_i'\hat{\bff}_t - {\bb_i^0}'{\bff}^0_t \\\nonumber
  & =\hat{\bb}_i'\hat{\bff}_t - ({\tilde{\bH}_4^{-1}\bb_i^0})'\tilde{\bH}_4^{'}{\bff}^0_t \\
   & =({\tilde{\bH}_4^{-1}\bb_i^0})' ( \hat{\bff}_t-\tilde{\bH}_4'{\bff}^0_t) +(\hat{\bb}_i-\tilde{\bH}_4^{-1}\bb_i^0)'\tilde{\bH}_4'{\bff}_t^0+(\hat{\bb}_i-\tilde{\bH}_4^{-1}\bb_i^0)'(\hat{\bff}_t -\tilde{\bH}_4'{\bff}^0_t). \label{ch-c0}
 \end{align}
The first term on the right-hand side of \eqref{ch-c0} is
\begin{eqnarray*}
	({\tilde{\bH}_4^{-1}\bb_i^0})' (\hat{\bff}_t-\tilde{\bH}_4'{\bff}^0_t) =({\tilde{\bH}_4^{-1}\bb_i^0})' (\hat{\bB}'\hat{\bB})^{-1} \tilde{\bQ} \bB^{0'}\be_t
	+ ({\tilde{\bH}_4^{-1}\bb_i^0})' (\hat{\bB}'\hat{\bB})^{-1}\left( \hat{\bB} - \bB^0\tilde{\bQ}'\right) '\be_t .
\end{eqnarray*}
By Lemma \ref{lem dist}(i), we have, if $ \alpha_r>1 / 2 $ and $\frac{N^{\frac{3}{2}-\alpha_r}}{T} \rightarrow 0
$,
\begin{eqnarray*}
	{\bb_i^0}' \bD\bN^{\frac{1}{2}}(\hat{\bff}_t-\tilde{\bH}_4'{\bff}^0_t) 
 &\CD& N(0, {\bb_i^0}'\bGamma_t\bb_i^0).
\end{eqnarray*}
Next, consider the second term on the right-hand side of \eqref{ch-c0},
\begin{align*}
	& \sqrt{T} (\hat{\bb}_i- \tilde{\bH}_4^{-1}\bb_i^0)'\tilde{\bH}_4'{\bff}_t^0 \\
 &=\sqrt{T} (\hat{\bb}_i - \tilde{\bQ}{\bb_i^0})'\tilde{\bH}_4'{\bff}_t^0
 +\sqrt{T} ( \tilde{\bH}_4^{-1}{\bb_i^0} - \tilde{\bQ}{\bb_i^0})'\tilde{\bH}_4'{\bff}_t^0 
  \\
 &=\sqrt{T} (\hat{\bb}_i - \tilde{\bQ}{\bb_i^0})'\tilde{\bH}_4'{\bff}_t^0+O_p(\sqrt{T}\tilde{\Delta}_{NT}),
\end{align*}
where $\tilde{\bH}_4^{-1}  - \tilde{\bQ}=O_p(\tilde{\Delta}_{NT})$ by Lemmas \ref{lem:H equality}(vi) and \ref{lem:FEB0}(iii). Applying the result in Lemma \ref{lem dist}(ii), if  $\frac{N^{1-\alpha_{r}}}{\sqrt{T}} \to 0$, and $\sqrt{T} N^{-\alpha_r} \rightarrow 0$, we have
\begin{align*}
    \sqrt{T}(\hat{\bb}_i- \tilde{\bH}_4^{-1}\bb_i^0)'\tilde{\bH}_4'{\bff}_t^0 \CD N(\mathbf{0}, {\bff}_t^{0'}\bPhi_i{\bff}_t^0).  &
\end{align*}
The third term on the right-hand side of \eqref{ch-c0} is dominated by the first and the second terms. Note that $\frac{N^{1-\alpha_r}}{\sqrt{T}} =(\frac{N^{\frac{3}{2}-\alpha_r}}{T})^{\frac{1}{2}}N^{\frac{1}{4}-\frac{1}{2}\alpha_r} \to 0$ if $\alpha_r>1/2$ and $\frac{N^{\frac{3}{2}-\alpha_{r}}}{T}  \to 0$ hold. Therefore, we come to that
if $\alpha_r>1/2$, $\frac{N^{\frac{3}{2}-\alpha_{r}}}{T}  \to 0$, and $\sqrt{T} N^{- \alpha_r} \rightarrow 0$, we have
\[
\frac{\hat{c}_{t,i}-c_{t,i}^*}{\sigma_{c(t,i)}}
\CD N(0, 1),\text{ with }\sigma_{c(t,i)}^2=V_{t,i}+U_{t,i},
\]
where $V_{t,i}={\bb_i^0}'\bD^{-1}\bN^{-\frac{1}{2}}\bGamma_t\bD^{-1}\bN^{-\frac{1}{2}}\bb_i^0$, $U_{t,i}=T^{-1}{\bff_t^0}'\bPhi_i\bff_t^0$.
Thus, we complete the proof.
\end{proof}

\begin{rem}[Comparison with data-dependent rotation matrices]
The conventional PC theory often studies the estimators relative to data-dependent rotated parameters, such as \(\hat{\bH}'\bff_t^*\) and \(\hat{\bH}^{-1}\bb_i^*\) in  \cite{Bai2003} for strong factor models, \(\hat{\bH}_4'\bff_t^*\) and \(\hat{\bQ}\bb_i^*\) in \citet[Proposition 7(i)]{BaiNg2023} for weak factor models. We clarify how such a formulation is related to our theory based on pseudo-true parameters. Recall that
 \[ \bff_t^0=\bH'\bff_t^*, \qquad \bb_i^0=\bH^{-1}\bb_i^* . \]
For the data-dependent rotation $\hat \bH$, the following equality holds
\[
   \hat \bH= \bH \tilde \bH,
\]
hence
\[
    \hat \bH' \bff_t^*=\tilde \bH'\bff_t^0,
    \qquad
    \hat \bH^{-1}\bb_i^*=\tilde \bH^{-1}\bb_i^0 .
\]
Therefore, for each \(t\),
\begin{align*}
    \hat \bff_t-\bff_t^0
    &=
   (\hat \bff_t-\hat \bH'\bff_t^*)
    -
    (\bH'\bff_t^*-\hat \bH'\bff_t^*)       \\
    &=
    (\hat \bff_t-\tilde \bH'\bff_t^0)
    +
    (\tilde \bH'-\bI_r)\bff_t^0 ,
\end{align*}
Similarly, 
\begin{align*} 
\hat{\bb}_i-\bb_i^0 
&= (\hat{\bb}_i-\hat{\bH}^{-1}\bb_i^*) + \{\hat{\bH}^{-1}\bb_i^*-\bH^{-1}\bb_i^*\} \\ 
&= (\hat{\bb}_i-\tilde{\bH}^{-1}\bb_i^0) + (\tilde{\bH}^{-1}-\bI_r)\bb_i^0 . 
\end{align*} 
Lemmas \ref{lemma f(fh-f0)}(iii) gives 
\[ 
\tilde{\bH}-\bI_r=O_p(\Delta_{NT}) 
\qquad  \text{or} \quad\bH^{-1}\hat{\bH}=\bI_r+O_p(\Delta_{NT}), 
\]
where $\Delta_{NT}$ is defined by \eqref{Delta} and is \(o_p(1)\) under the conditions $\alpha_1 <2\alpha_r$ and $ \frac{N^{1-\alpha_r}}{T}\to 0$. 
For the rates, it is convenient to use the data-dependent rotations \(\hat{\bH}_4\) and \(\hat{\bQ}\), because \(\hat{\bff}_t-\hat{\bH}_4'\bff_t^*\) and \(\hat{\bb}_i-\hat{\bQ}\bb_i^*\) admit the simplest tractable expansions. These data-dependent rotations are also asymptotically equivalent to the population rotation in the same relative sense. More precisely,
\[ 
\tilde{\bH}_4-\bI_r= O_p(\tilde{\Delta}_{NT}), 
\qquad \tilde{\bQ}-\bI_r=O_p(\tilde{\Delta}_{NT}), 
\]
or equivalently,
\[ 
\bH^{-1}\hat{\bH}_4=\bI_r+O_p(\tilde{\Delta}_{NT}), \qquad \hat{\bQ}\bH=\bI_r+O_p(\tilde{\Delta}_{NT}). 
\]
where $\tilde{\Delta}_{NT}$ is defined by \eqref{eq:TDleta}.
Thus, the estimation error centered at pseudo-true parameters consists of the usual error around a data-dependent rotated target plus an additional asymptotically negligible term. The latter is the cost of replacing the data-dependent target by the population one.
\end{rem}

\subsection{Proof of Theorem \ref{thm:forecast1}}

\begin{proof}
By the definition of $\hat{\bdelta}$, we have the following decomposition:
\begin{align*} 
	&	\sqrt{T} ( \hat{\bdelta}-\bdelta^0)\\
  & =\left( \frac{1}{T} {\hat{\bZ}'\hat{\bZ}}\right)^{-1}\frac{1}{\sqrt{T}} \hat{\bZ}'\bepsilon + \left( \frac{1}{T} {\hat{\bZ}'\hat{\bZ}}\right)^{-1}\frac{1}{\sqrt{T}} \hat{\bZ}'(\bF^0 -\hat{\bF})\bgamma^0  \\
  & =\left( \frac{1}{T} { {\bZ^0}' {\bZ^0}}\right)^{-1}\frac{1}{\sqrt{T}}  {\bZ^0}'\bepsilon+\left( \frac{1}{T} { {\bZ^0}' {\bZ^0}}\right)^{-1}\frac{1}{\sqrt{T}}  (\hat{\bZ}-{\bZ^0})'\bepsilon + \left[\left( \frac{1}{T} {\hat{\bZ}'\hat{\bZ}}\right)^{-1}-   \left( \frac{1}{T} { {\bZ^0}' {\bZ^0}}\right)^{-1}  \right]\frac{1}{\sqrt{T}} \hat{\bZ}'\bepsilon \\
  & \quad \quad +\left( \frac{1}{T}{\hat{\bZ}'\hat{\bZ}}\right)^{-1}\frac{1}{\sqrt{T}} \hat{\bZ}'(\bF^0 -\hat{\bF})\bgamma^0.  \\
\end{align*}
The third term is bounded by $O_p\left( \tilde{\Delta}_{NT} \right)$ because
\begin{align*}
   & \frac{1}{T} {\hat{\bZ}'\hat{\bZ}}\\
   & =\frac{1}{T} {\bZ^0}'\bZ^0 +\frac{1}{T} ({\hat{\bZ}-\bZ^0)'\hat{\bZ}}+\frac{1}{T} {\bZ^0}'(\hat{\bZ}-\bZ^0) \\
  &=\frac{1}{T} {\bZ^0}'\bZ^0  +
  O_p\left( \tilde{\Delta}_{NT} \right) ,
\end{align*}
 by Lemma \ref{lem:forecast1}(i).
The second and fourth terms are bounded by
\begin{align*}
&		 \left\| \left( \frac{1}{T} { {\bZ^0}' {\bZ^0}}\right)^{-1}\frac{1}{\sqrt{T}}(\bZ^0-\hat{\bZ})'\bepsilon \right\|_{\F}
=  O_p\left(\frac{N^{1 - \alpha_{r}}}{\sqrt{T}} \right) + O_p\left( \frac{\sqrt{T}}{N^{\alpha_r}}\right),\\
&\left\| \left( \frac{1}{T}{\hat{\bZ}'\hat{\bZ}}\right)^{-1}\frac{1}{\sqrt{T}} \hat{\bZ}'(\bF^0 -\hat{\bF})\bgamma^0 \right\|_{\F}
	=  O_p\left(\frac{N^{1-\alpha_r}}{\sqrt{T}}\right)+
 + O_p\left( \frac{\sqrt{T}}{N^{\alpha_r}}\right),
  \end{align*}
  where we have used Lemmas \ref{lem:forecast1}(ii) and (i). Collecting these non-dominating terms, 
\begin{align*} 
	&	\sqrt{T} ( \hat{\bdelta}-\bdelta^0)
  = \left( \frac{1}{T} {{\bZ^0}'{\bZ^0}}\right)^{-1}\frac{1}{\sqrt{T}} {\bZ^0}'\bepsilon +   O_p\left(\frac{N^{1 - \alpha_{r}}}{\sqrt{T}} \right) + O_p\left(\sqrt{T} N^{-\alpha_r}\right).
\end{align*}
If $\frac{N^{1-\alpha_r}}{\sqrt{T}} \to  0$, $\sqrt{T}N^{-\alpha_r} \to 0$, and under Assumption \ref{ass: augreg}(iii) 
\[
\left( \frac{1}{T} {{\bZ^0}'{\bZ^0}}\right)^{-1}\frac{1}{\sqrt{T}} {\bZ^0}'\bepsilon
\CD  N(\mathbf{0}, \bSigma_{\delta^0}),
\]
thus, we complete the proof.
\end{proof}


\subsection{Proof of Theorem \ref{thm:forecast2}}

\begin{proof}
We start with the decomposition using the rotation matrix ${\bH}$
\begin{align}
\nonumber
& \hat{y}_{T+h \mid T}-y_{T+h \mid T} \\ \nonumber
& =\hat{\bgamma}^{\prime} \hat{\bff}_T+\hat{\bbeta}^{\prime} \bw_T-\bgamma^{*\prime} \bff_T^*-\bbeta^{\prime} \bw_T \\\nonumber
& =\left(\hat{\bgamma}- {\bH}^{-1} \bgamma^*\right)^{\prime} \hat{\bff}_T+\bgamma^{*\prime}  {\bH} ^{-1\prime}\left(\hat{\bff}_T- {\bH}^\prime  \bff_T^*\right)+(\hat{\bbeta}-\bbeta)^{\prime} \bw_T \\\nonumber
& =\hat{\bz}_T^{\prime}(\hat{\bdelta}-\bdelta^0)+\bgamma^{*\prime}  {\bH} ^{-1\prime}\left(\hat{\bff}_T- {\bH}^\prime  \bff_T^*\right) \\\nonumber
& =T^{-1 / 2} \hat{\bz}_T^{\prime}[\sqrt{T}(\hat{\bdelta}-\bdelta^0)]+ \bgamma^{*\prime} \bH^{-1\prime} \bN^{-1 / 2}\left[\bN^{1 / 2}\left(\hat{\bff}_T-\bH^\prime \bff_T^*\right)\right] \\\nonumber
& =T^{-1 / 2} {\bz}_T^{0\prime}[\sqrt{T}(\hat{\bdelta}-\bdelta^0)]+T^{-1 / 2} (\hat{\bz}_T-\bz^0_T)^{\prime}[\sqrt{T}(\hat{\bdelta}-\bdelta^0)]\\
& \quad \quad + \bgamma^{*\prime} \bH^{-1\prime} \bN^{-1 / 2}\left[\bN^{1 / 2}\left(\hat{\bff}_T- \tilde{\bH}_4^\prime \bff_T^0 \right)\right]+ \bgamma^{*\prime} \bH^{-1\prime} \left(\tilde{\bH}_4^\prime - \bI \right) \bff_T^0 . \label{yh-y}
\end{align}
Consider the first term on the right-hand side of the above equation,
\begin{eqnarray*} 
	&&T^{-1 / 2} {\bz}_T^{0\prime}  \sqrt{T} ( \hat{\bdelta}-\bdelta^0)\\
  &&=T^{-1 / 2} {\bz}_T^{0\prime}\left( \frac{1}{T} {\hat{\bZ}'\hat{\bZ}}\right)^{-1}\frac{1}{\sqrt{T}} \hat{\bZ}'\bepsilon +T^{-1 / 2} {\bz}_T^{0\prime} \left( \frac{1}{T} {\hat{\bZ}'\hat{\bZ}}\right)^{-1}\frac{1}{\sqrt{T}} \hat{\bZ}'(\bF^0 -\hat{\bF}){\bH}^{-1}\bgamma^*  \\
 && =T^{-1 / 2} {\bz}_T^{0\prime} \left( \frac{1}{T} {{\bZ^0}'{\bZ^0}}\right)^{-1}\frac{1}{\sqrt{T}} {\bZ^0}'\bepsilon \\
 &&\quad \quad+  \frac{1}{\sqrt{T}}\left[O_p\left(\frac{N^{1 - \alpha_{r}}}{\sqrt{T}} \right) + O_p\left(\sqrt{T} N^{-\alpha_r}\right)\right],
\end{eqnarray*}
where the term in the bracket is dominated by the first one if $\frac{N^{1 - \alpha_{r}}}{\sqrt{T}} \to 0$ and $\sqrt{T} N^{-\alpha_r} \to 0$.
The second term on the right-hand side of \eqref{yh-y} is dominated by the first one, thus, we ignore it. 

Next, consider the third term in \eqref{yh-y} and Lemma \ref{lemma bh-b0*e} implies
\begin{eqnarray*}
    && \bgamma^{*\prime}  {\bH}^{-1\prime} \bN^{-\frac{1}{2}} \left[\bN^{1 / 2}\left(\hat{\bff}_T-\tilde{\bH}_4^{\prime} \bff_T^0\right)\right] \\
    && = \bgamma^{*\prime}  {\bH}^{-1\prime} \bN^{-\frac{1}{2}} \bD^{-1}\left[\bD\bN^{\frac{1}{2}}(\hat{\bB}'\hat{\bB})^{-1} \tilde{\bQ} \bB^{0'}\be_t
	+ \bD\bN^{\frac{1}{2}}(\hat{\bB}'\hat{\bB})^{-1}\left( \hat{\bB} - \bB^0\tilde{\bQ}'\right) '\be_t
    \right]\\
    && = \bgamma^{*\prime}  {\bH}^{-1\prime} \bN^{-\frac{1}{2}} \bD^{-1}  
  \left[\bD\bN^{\frac{1}{2}}(\hat{\bB}'\hat{\bB})^{-1} \tilde{\bQ} \bB^{0'}\be_t\right]
	+ \bgamma^{*\prime}  {\bH}^{-1\prime} \bN^{-\frac{1}{2}}  \left[ \bN^{\frac{1}{2}}(\hat{\bB}'\hat{\bB})^{-1}\left( \hat{\bB} - \bB^0\tilde{\bQ}'\right) '\be_t\right]
     \\
      && = \bgamma^{*\prime}  {\bH}^{-1\prime} \bN^{-\frac{1}{2}} \bD^{-1}  
  \left[\bD\bN^{\frac{1}{2}}(\hat{\bB}'\hat{\bB})^{-1}   \bB^{0'}\be_t\right]
	+  O_p(\| \bN^{-\frac{1}{2}}\|_F) O_p(\tilde{\Delta}_{NT})  \\
 && \quad \quad + O_p(\| \bN^{-\frac{1}{2}}\|_F)  \left[ O_p(N^{\frac{1}{2}-\alpha_{r}} )+ O_p\left(\frac{N^{\frac{3}{2}-\alpha_{r}}}{T} \right)+ O_p\left(\sqrt{\frac{N^{1-\alpha_{r}}}{T}}\right)\right],
\end{eqnarray*}	
where the last two terms are dominated by the first one if $\frac{N^{\frac{3}{2} - \alpha_{r}}}{T} \to 0$ and $\frac{1}{2}< \alpha_r$. 

The fourth term in \eqref{yh-y} becomes
\begin{eqnarray*}
    &&   \bgamma^{*\prime} \bH^{-1\prime} \left(\tilde{\bH}_4^{\prime} - \bI \right) \bff_T^0 
      =
      \frac{1}{\sqrt{T}} O_p(\sqrt{T}\tilde{\Delta}_{NT}),
\end{eqnarray*}	
which is dominated by the first term in \eqref{yh-y} if $\frac{N^{1 - \alpha_{r}}}{\sqrt{T}} \to 0$ and $\sqrt{T} N^{-\alpha_r} \to 0$. Collecting these terms, if $\sqrt{T} N^{-\alpha_r}\to 0$, $\frac{1}{2}< \alpha_r$, and $\frac{N^{\frac{3}{2}-\alpha_{r}}}{T} \to 0$, we obtain
\[
\frac{\left(\hat{y}_{T+h \mid T}-y_{T+h \mid T}\right)}{{\sigma}_{T+h \mid T}} \stackrel{d}{\longrightarrow} N(0,1),
\]
where ${\sigma}_{T+h \mid T}^2=
T^{-1} {\bz}_T^{0\prime}\bSigma_{\delta^0}{\bz}_T^0
+{\bgamma}^{0 \prime} 
\bD^{-1}\bN^{-1/2}\bGamma_T\bN^{-1/2}\bD^{-1}
{\bgamma}^0$.

\end{proof}

\section{Related Lemmas and their Proofs}
\label{subsec:related lemmas}
\setcounter{lem}{0}
\renewcommand{\thelem}{B.\arabic{lem}}
\setcounter{equation}{0}
\renewcommand{\theequation}{B.\arabic{equation}} 

\begin{lem}
\label{lem:BEF} The following results hold under Assumptions \ref{cond:eigen}--\ref{assumption BaiNgA3}:
\begin{flalign*}
&(i)~ \frac{1}{NT}\be_t' \bE'\bF^0=O_p\left(\frac{1}{\sqrt{NT}}\right)+ O_p\left(\frac{1}{T}\right),&\\
&(ii) ~
\frac{1}{T}\be_i'\bE \bB^0 {\bN}^{-\frac{1}{2}}= O_p\left(\frac{1}{\sqrt{N^{\alpha_r}}}\right)+ O_p\left(\frac{1}{\sqrt{T}}\right), &\\
& (iii) ~
	\bN^{-\frac{1}{2}}{\bB^0}'\bE'=O_p(\sqrt{T}), &\\
	 & (iv) ~ \frac{\bN^{-\frac{1}{2}} \hat{\bF}' \bF^0 \bN^{\frac{1}{2}}}{T}=O_p(1). & 
\end{flalign*}
\end{lem}

\begin{proof}
       (i) We have
    \begin{align*}
       & \frac{1}{NT}\be_t' \bE'\bF^0\\
       &= \frac{1}{NT} \sum_{i=1}^N \sum_{s=1}^Te_{t,i}e_{s,i}\bff_s^{0'}\\
       & =\frac{1}{NT} \sum_{i=1}^N \sum_{s=1}^T\bff_s^{0'}[e_{t,i}e_{s,i}-\mathbb{E}(e_{t,i}e_{s,i})]+\frac{1}{T}  \left( \frac{1}{N} \sum_{s=1}^T\bff_s^{0'} \sum_{i=1}^N \mathbb{E}(e_{t,i}e_{s,i})\right)\\
       & =O_p\left(\frac{1}{\sqrt{NT}}\right)+ O_p\left(\frac{1}{T}\right).
    \end{align*}
    Because for all $i$, $ | \mathbb{E}\left(e_{s, i} e_{t, i}\right)|\leq| \gamma_{s, t}|$ for some $\gamma_{s, t}$ such that $\sum_{t=1}^T\left|\gamma_{s, t}\right| \leq M
$ by Assumption \ref{assumption errors}(iii), and $\|\bff_s^0 \|_2$ is bounded by $\E\|\bff_t^0\|^4 \le M$ in Assumption \ref{assumption BaiNgA3}(i). \\
    (ii) We have
    \begin{align*}
       & \frac{1}{T}\be_i' \bE\bB^0{\bN}^{-\frac{1}{2}}\\
       &= \frac{1}{T} \sum_{j=1}^N \sum_{t=1}^Te_{t,i}e_{t,j}\bb_j^{0'}{\bN}^{-\frac{1}{2}}\\
       & =\frac{1}{T} \sum_{j=1}^N \sum_{t=1}^T\bb_j^{0'}{\bN}^{-\frac{1}{2}}[e_{t,i}e_{t,j}-\mathbb{E}(e_{t,i}e_{t,j})]+\frac{1}{T}   \sum_{t=1}^T  \sum_{j=1}^N \bb_j^{0'}{\bN}^{-\frac{1}{2}}  \mathbb{E}(e_{t,i}e_{t,j}) \\
       & =O_p\left(\frac{1}{\sqrt{T}}\right)+ O_p\left(\frac{1}{\sqrt{N^{\alpha_r}}}\right).
    \end{align*}
    Under weak cross-sectional dependence as in Assumption \ref{assumption errors}(iv),   $(1 / T N) \sum_{t=1}^T \sum_{j=1}^N\left|\mathbb{E}(e_{t,i}e_{t,j})\right| \leq(1 / N) \sum_{j=1}^N \left|\tau_{ i,j}\right| =O\left(N^{-1}\right)$. By $\E\|\bb_i^0\|^4 \le M$ in Assumption \ref{assumption BaiNgA3}(i), $\|\bb_i^0 \|_{2} $ is bounded and $||\bb_j^{0'}\bN^{-\frac{1}{2}}||_F \le||\bb_j^{0'}||_F ||\bN^{-\frac{1}{2}}||_F \le  
\frac{1}{\sqrt{N^{\alpha_r}}} O_p(1)$. 
    The upper bounds in (i) and (ii) are not consistent with \citet[Assumption A3']{BaiNg2023}, because we impose a moment restriction related to $\bb_i^0$ in Assumption \ref{assumption BaiNgA3}(iv).


(iii) Assumption \ref{assumption BaiNgA3}(ii) implies 
 \[
 \frac{\bN^{-\frac{1}{2}}{\bB^0}'\bE' \bE\bB^0 \bN^{-\frac{1}{2}}}{T} =\frac{1}{T} \sum_{t=1}^{T} \left[ \left( \bN^{-\frac{1}{2}}\sum_{i=1}^{N} \bb_i^0e_{t,i}\right)  \left( \bN^{-\frac{1}{2}}\sum_{i=1}^{N} \bb_i^0e_{t,i}\right)'\right]=O_p(1),
\]
where $\bB^0=(\bb_1^0, \cdots, \bb_N^0)'$. Thus,
\[
\left\| \bE\bB^0 \bN^{-\frac{1}{2}}\right\|_{\F}^2=tr(\bN^{-\frac{1}{2}}{\bB^0}'\bE'\bE\bB^0 \bN^{-\frac{1}{2}})=O_p(T),
\]
$\left\|\bE\bB^0 \bN^{-\frac{1}{2}} \right\|_{\F} =O_p( \sqrt{T })$, and $\left\|\bE\bB^0  \right\|_{\F} =O_p(\sqrt{TN^{\alpha_1}})$.

   (iv) A comparable result to \citet[Lemma 1]{BaiNg2023} is
\[
\bN^{ -1} \hat{\bLambda} \CP \bN^{ -1} \bLambda  ,
\]
which can be obtained by Weyl's inequalities in Lemma \ref{lemma order of lambda}.
The equations $\frac{1}{T}\bX\bX' \hat{\bF}=\hat{\bF} \hat{\bLambda}$ and $\hat{\bF}'\hat{\bF}/T =\bI_r$ imply
\[
\frac{1}{T^2} \bN^{-\frac{1}{2}} \hat{\bF}' \bX\bX' \hat{\bF} \bN^{-\frac{1}{2}}  = \bN^{-1}  \hat{\bLambda}.
\]
The dominating term of the left-hand side is 
\[
\frac{\bN^{-\frac{1}{2}} \hat{\bF}' \bF^0 \bN^{\frac{1}{2}}}{T} \bN^{-\frac{1}{2}}{\bB^0}'\bB^0\bN^{-\frac{1}{2}}  \frac{\bN^{\frac{1}{2}}{\bF^0}' \hat{\bF}  \bN^{-\frac{1}{2}}}{T} \CP \bN^{ -1} \bLambda  .
\]
Thus, $\frac{\bN^{-\frac{1}{2}} \hat{\bF}' \bF^0 \bN^{\frac{1}{2}}}{T}=O_p(1)$.
\end{proof}

\begin{lem}
	\label{lemma order of lambda}
	 Suppose that Assumption \ref{assumption errors} holds. If $\frac{N^{1-\alpha_r}}{T} \to 0$, we have $\hat{\lambda}_k=\lambda_k\left[ T^{-1}\bX\bX'\right] \asymp N^{\alpha_{k}}$  with high probability.
\end{lem}

\begin{proof}
Let $\sigma_k[\bA]$ be $k$-th largest singular value of matrix $\bA$. By the definition of WF models, we have
    \begin{eqnarray*}
		\lambda_k=\lambda_k\left[ \frac{\bB^*{\bF^*}'\bF^*{\bB^*}'}{T}\right] =\lambda_k\left[ {\bB^*}'\bB^*\frac{{\bF^*}'\bF^*}{T}\right] \asymp N^{\alpha_{k}}.
	\end{eqnarray*}
By the singular value version of Weyl's inequalities, 
	\begin{eqnarray*}
		\sigma_{k+l-1}[A+B] \le \sigma_k[A]+\sigma_l[B], \quad 1 \le k, l \le \min{(N,T)}.
	\end{eqnarray*}
	We first show the upper bound for $\hat{\lambda}_k$:
\begin{eqnarray*}
	\nonumber
	\sigma_{k}[\bX] &\le& \sigma_k[\bF^*{\bB^*}']+\sigma_1[\bE]  \\ 
	\sqrt{\lambda_k\left[\frac{\bX\bX'}{T} \right]} &\le& \sqrt{\lambda_k\left[ \frac{\bF^*{\bB^*}'\bB^*{\bF^*}'}{T}\right]} +\sqrt{\lambda_1\left[\frac{\bE\bE'}{T} \right]} \\\nonumber
	\lambda_k\left[\frac{\bX\bX'}{T} \right] &\le& \lambda_k\left[ \frac{\bF^*{\bB^*}'\bB^*{\bF^*}'}{T}\right]+\lambda_1\left[\frac{\bE\bE'}{T} \right]+2\sqrt{\lambda_k\left[ \frac{\bF^*{\bB^*}'\bB^*{\bF^*}'}{T}\right]\lambda_1\left[\frac{\bE\bE'}{T} \right]} \\\nonumber
	\lambda_k\left[\frac{\bX\bX'}{T} \right]  &\le& N^{\alpha_k} + \left(\frac{N}{T}+1\right)+\sqrt{N^{\alpha_{k}}\left(\frac{N}{T}+1\right)}.
\end{eqnarray*}
Next, the lower bound for $\hat{\lambda}_k$ becomes 
\begin{eqnarray*}
	\nonumber
	\sigma_{k}[\bX-\bE] &\le& \sigma_k[\bX]+\sigma_1[\bE]  \\
	\nonumber
\sigma_k[\bF^*{\bB^*}'] &\le& \sigma_k[\bX]+\sigma_1[\bE] \\
\nonumber
\sqrt{\lambda_k\left[ \frac{\bF^*{\bB^*}'\bB^*{\bF^*}'}{T}\right]}&\le& \sqrt{	\lambda_k\left[\frac{\bX\bX'}{T} \right] }+\sqrt{\lambda_1\left[\frac{\bE\bE'}{T} \right]}\\
	\sqrt{\lambda_k\left[\frac{\bX\bX'}{T} \right]}  &\ge& \sqrt{\lambda_k\left[ \frac{\bF^*{\bB^*}'\bB^*{\bF^*}'}{T}\right]}- \sqrt{\lambda_1\left[\frac{\bE\bE'}{T} \right]}.
\end{eqnarray*}
Thus, $ \hat{\lambda}_k=\lambda_k\left[ \frac{\bX\bX'}{T}\right] \asymp N^{\alpha_{k}}$. Note that $N^{\alpha_k}$ dominates $\frac{N}{T}$ if $\frac{N^{1-\alpha_r}}{T} \to 0$.\\
\end{proof}


\begin{lem}
	\label{lemma f(fh-f0)} Define 
 \begin{align}
 \label{Delta}
   \Delta_{NT}
	=
	   \frac{N^{1-\alpha_r}}{T}+ N^{\frac{1}{2}\alpha_1- \alpha_r} \frac{N^{1- \alpha_r}}{T} + N^{\frac{1}{2}\alpha_{1}-\frac{3}{2}\alpha_{r}} .
 \end{align}
  Suppose that Assumptions \ref{cond:eigen}--\ref{assumption BaiNgA3} hold. Then, we have
		\begin{flalign*}
			&(i) ~		\left\|  \frac{1}{T} \bE'(\hat{\bF}- \bF^0\tilde{\bH}_4) \right\|_{\F} 
			= O_p \left( \sqrt{\frac{N^{1-\alpha_r}}{T}}\frac{N^{1-\frac{1}{2}\alpha_r}}{T} \right)+ O_p \left( \frac{N^{1-\frac{1}{2}\alpha_r}}{T}  \right)+ O_p \left(  N^{-\frac{1}{2}\alpha_{r}} \right), & \\
  & (ii) ~		\left\|  \frac{1}{T} {\bB^0}'\bE'(\hat{\bF}- \bF^0\tilde{\bH}_4) \right\|_{\F} 
			= O_p \left( N^{\frac{1}{2}\alpha_{1}-\frac{1}{2}\alpha_{r}}\right)+ O_p \left(  N^{\frac{1}{2}\alpha_1-\frac{1}{2}\alpha_r} \frac{N^{1-\frac{1}{2}\alpha_r}}{T} \right), & \\
   &	\quad\quad	\left\|  \frac{1}{T} \bN^{-\frac{1}{2}}{\bB^0}'\bE'(\hat{\bF}- \bF^0\tilde{\bH}_4) \right\|_{\F} 
			= O_p \left( N^{ -\frac{1}{2}\alpha_{r}}\right)+ O_p \left(    \frac{N^{1- \alpha_r}}{T} \right), & \\
		&(iii) ~		\left\| \frac{1}{T}{\bF^0}'(\hat{\bF}-\bF^0\tilde{\bH}) \right\|_{\F} 
		= O_p \left(\Delta_{NT} \right), & \\
 & (iv) ~ \left\| \frac{1}{T}{\be_i}'(\hat{\bF}-\bF^0\tilde{\bH}_4) \right\|_{\F} 
		=  O_p \left( \frac{N^{1-\alpha_r}}{T} \right)+ O_p \left(   N^{-\alpha_r}\right)
   + O_p \left(\frac{1}{\sqrt{TN^{\alpha_r}}}\right) .&  \\
	\end{flalign*}	
\end{lem}

\begin{proof}(i) By the definition of $\hat{\bF}$, we have the following decomposition:
    \begin{align*}
		&	  \left\|\frac{1}{T} \bE'(\hat{\bF}- \bF^0 \tilde{\bH}_4) \right\|_{\F} \\
  & = \left\| 	\frac{1}{T^2} \bE'\bE\bE'\hat{\bF} \hat{\bLambda}^{-1} +	\frac{1}{T^2} \bE'\bE\bB^0 {\bF^0}'\hat{\bF} \hat{\bLambda}^{-1}\right\|_{\F} \\
  & \le \left\| 	\frac{1}{T^2} \bE'\bE\bE'\right\|_F  	\left\|\hat{\bF} \hat{\bLambda}^{-1} \right\|_F +	\left\|\frac{1}{T} \bE'\bE\right\|_F  	\left\|\bB^0 \bN^{-\frac{1}{2}}\right\|_F  	\left\|\bN^{\frac{1}{2}}\frac{1}{T}{\bF^0}'\hat{\bF}\bN^{-\frac{1}{2}}\bN^{\frac{1}{2}} \hat{\bLambda}^{-1}\right\|_{\F} \\
  & =  		\frac{1}{T^2} O_p\left(N^{3/2}+T^{3/2} \right)O_p\left(\sqrt{T}N^{-\alpha_{r}}\right) +O_p\left( \frac{N}{T}+1\right) O_p\left(N^{-\frac{1}{2}\alpha_{r}}\right) \\
  &  = O_p\left(\sqrt{\frac{N^{1-\alpha_r}}{T}}\frac{N^{1-\frac{1}{2}\alpha_r}}{T}\right)
  +O_p\left(\frac{N^{1-\frac{1}{2}\alpha_r}}{T}\right)
  +O_p\left(N^{-\frac{1}{2}\alpha_{r}}\right).
	\end{align*}	
(ii) We have
 \begin{eqnarray*}
		&&	\left\|\frac{1}{T}{\bB^0}'\bE' (\hat{\bF}-\bF^0 \tilde{\bH}_4) \right\|_{\F} \\
  && = \left\| 	\frac{1}{T^2} {\bB^0}'\bE'\bE\bE'\hat{\bF} \hat{\bLambda}^{-1} +	\frac{1}{T^2} {\bB^0}'\bE'\bE\bB^0 {\bF^0}'\hat{\bF} \hat{\bLambda}^{-1}\right\|_{\F}   \\
   && \le \left\| 	\frac{1}{T^2} {\bB^0}'\bE'\right\|_F  	\left\|\bE\bE'\right\|_F  	\left\|\hat{\bF} \hat{\bLambda}^{-1} \right\|_F +	\left\|\bN^{\frac{1}{2}}\frac{1}{T} \bN^{-\frac{1}{2}}{\bB^0}'\bE'\bE\bB^0 \bN^{-\frac{1}{2}}\right\|_F  	\left\|\bN^{\frac{1}{2}}\frac{1}{T}{\bF^0}'\hat{\bF}\bN^{-\frac{1}{2}}\bN^{\frac{1}{2}} \hat{\bLambda}^{-1}\right\|_{\F} \\
  && =  		O_p\left( N^{\frac{1}{2}\alpha_1-\frac{1}{2}\alpha_r}\left( \frac{N^{1-\frac{1}{2}\alpha_r}}{T}+N^{-\frac{1}{2}\alpha_r}\right) \right) 
  +O_p\left(N^{\frac{1}{2}\alpha_1}\right)O_p\left( N^{-\frac{1}{2}\alpha_r}\right) \\
  &&  = O_p\left(  N^{\frac{1}{2}\alpha_1-\frac{1}{2}\alpha_r} \frac{N^{1-\frac{1}{2}\alpha_r}}{T}\right)
  +O_p\left( N^{\frac{1}{2}\alpha_1-\frac{1}{2}\alpha_r}\right).
	\end{eqnarray*}	
 We also obtain
	 \begin{eqnarray*}
		&&	\left\|\frac{1}{T}\bN^{-\frac{1}{2}}{\bB^0}'\bE' (\hat{\bF}-\bF^0 \tilde{\bH}_4) \right\|_{\F} \\
  && = \left\| 	\frac{1}{T^2} \bN^{-\frac{1}{2}}{\bB^0}'\bE'\bE\bE'\hat{\bF} \hat{\bLambda}^{-1} +	\frac{1}{T^2} \bN^{-\frac{1}{2}}{\bB^0}'\bE'\bE\bB^0 {\bF^0}'\hat{\bF} \hat{\bLambda}^{-1}\right\|_{\F}   \\
   && \le \left\| 	\frac{1}{T^2} \bN^{-\frac{1}{2}}{\bB^0}'\bE'\right\|_F  	\left\|\bE\bE'\right\|_F  	\left\|\hat{\bF} \hat{\bLambda}^{-1} \right\|_F +	\left\| \frac{1}{T} \bN^{-\frac{1}{2}}{\bB^0}'\bE'\bE\bB^0 \bN^{-\frac{1}{2}}\right\|_F  	\left\|\bN^{\frac{1}{2}}\frac{1}{T}{\bF^0}'\hat{\bF}\bN^{-\frac{1}{2}}\bN^{\frac{1}{2}} \hat{\bLambda}^{-1}\right\|_{\F} \\
  &&  =  O_p\left( \frac{N^{1- \alpha_r}}{T}\right)+O_p\left(N^{ - \alpha_r} \right)+ O_p\left(N^{ -\frac{1}{2}\alpha_r}\right)\\
  &&  =  O_p\left(\frac{N^{1- \alpha_r}}{T}\right) + O_p\left(N^{ -\frac{1}{2}\alpha_r}\right).
	\end{eqnarray*}	

(iii) We have
	\begin{eqnarray*}
&&	\frac{1}{T}{\bF^0}'(\hat{\bF}-\bF^0\tilde{\bH})  \\
		&&= \frac{1}{T} {\bF^0}'\left( \frac{1}{T} \bE\bE'\hat{\bF}+\frac{1}{T} \bF^0{\bB^0}'\bE'\hat{\bF}+\frac{1}{T} \bE\bB^0{\bF^0}'\hat{\bF} \right) \hat{\bLambda}^{-1} \\
		&&=\frac{1}{T^2} {\bF^0}' \bE\bE'\hat{\bF}\hat{\bLambda}^{-1}+\frac{1}{T}{\bB^0}'\bE'\hat{\bF}\hat{\bLambda}^{-1}+\frac{1}{T^2} {\bF^0}' \bE\bB^0{\bF^0}'\hat{\bF} \hat{\bLambda}^{-1}.
	\end{eqnarray*}
The first term on the right-hand side of the above equation is bounded by
\[
\left\| \frac{1}{T^2} {\bF^0}' \bE\bE'\hat{\bF}\hat{\bLambda}^{-1}\right\|_{\F} \le  O_p\left( \left(\frac{N}{T} +1\right)N^{-\alpha_{r}}\right).
\]
Next, we consider the second term:
\begin{align*}
		&\left\| 	\frac{1}{T}{\bB^0}'\bE'\hat{\bF}\hat{\bLambda}^{-1} \right\|_{\F}    \\
  &\le\left\| 	\frac{1}{T}{\bB^0}'\bE'(\hat{\bF}-\bF^0\tilde{\bH}_4)\hat{\bLambda}^{-1}  \right\|_{\F} + \left\| 	\frac{1}{T}{\bB^0}'\bE'\bF^0\tilde{\bH}_4\hat{\bLambda}^{-1}  \right\|_{\F} \\
  &\le  \left\|	\frac{1}{T}{\bB^0}'\bE'(\hat{\bF}-\bF^0\tilde{\bH}_4) \right\|_{\F} \left\|\hat{\bLambda}^{-1} \right\|_{\F}+ \left\| 	\frac{1}{T}{\bB^0}'\bE'\bF^0 \right\|_{\F} \left\| \tilde{\bH}_4\hat{\bLambda}^{-1}  \right\|_{\F}\\
		&= O_p\left(
 N^{\frac{1}{2}\alpha_1- \alpha_r} \frac{N^{1- \alpha_r}}{T}\right)
 +O_p\left(
N^{\frac{1}{2}\alpha_1-\frac{3}{2}\alpha_r}\right) +O_p\left(
\frac{\sqrt{TN^{\alpha_{1}}}}{T}N^{-\alpha_{r}}\right),
  \end{align*}
  where we have used Lemma \ref{lemma f(fh-f0)} (ii) and Assumption \ref{assumption BaiNgA3}(vi). In addition, by $2\sqrt{ab} \le a+b$, we have $\frac{N^{\frac{1}{2}\alpha_1-\alpha_r}}{\sqrt{T}}\le \frac{N^{\frac{1}{2}\alpha_1-\alpha_r}}{\sqrt{T}}N^{\frac{1}{2}-\frac{1}{4}\alpha_1-\frac{1}{4}\alpha_r} \le \frac{1}{2}\left(\frac{N^{1-\alpha_r}}{T}+N^{\frac{1}{2}\alpha_1-\frac{3}{2}\alpha_r} \right) $.
    
The third term is bounded by
	\begin{eqnarray*}
	&&\left\| 	\frac{1}{T^2} {\bF^0}' \bE\bB^0{\bF^0}'\hat{\bF} \hat{\bLambda}^{-1}\right\|_{\F}   \le \left\| 	\frac{1}{T^2} {\bF^0}' \bE\bB^0 \bN^{-\frac{1}{2}}\right\|_{\F} \left\| \bN^{\frac{1}{2}}{\bF^0}'\hat{\bF}\bN^{-\frac{1}{2}}\bN^{\frac{1}{2}} \hat{\bLambda}^{-1}\right\|_{\F}  
 = O_p\left(\frac{1}{\sqrt{TN^{\alpha_r}}}\right).
\end{eqnarray*}
Collecting terms, we obtain
	\begin{eqnarray*}
	&&\left\| 	\frac{1}{T}{\bF^0}'(\hat{\bF}-\bF^0\tilde{\bH}) \right\|_{\F}   \\
 && \le  O_p\left(\left( \frac{N}{T} +1\right)N^{-\alpha_{r}} \right)
 +O_p\left( N^{\frac{1}{2}\alpha_1- \alpha_r} \frac{N^{1- \alpha_r}}{T}\right)
 +O_p\left(N^{\frac{1}{2}\alpha_1-\frac{3}{2}\alpha_r} \right)    
 \\
	&&=   O_p\left( \frac{N^{1-\alpha_r}}{T}\right)+ O_p\left(N^{\frac{1}{2}\alpha_1- \alpha_r} \frac{N^{1- \alpha_r}}{T}\right) + O_p\left(N^{\frac{1}{2}\alpha_{1}-\frac{3}{2}\alpha_{r}} \right),
\end{eqnarray*}
which is $o_p(1)$ if $\frac{N^{1-\alpha_r}}{T} \to 0$ and $\frac{1}{2}\alpha_1<\alpha_r$.

(iv) We have
 \begin{eqnarray*}
	&& \left\| \frac{1}{T} \be_i'(\hat{\bF}- \bF^0\tilde{\bH}_4)\right\|_{\F}	\\ && = \left\| \frac{1}{T} \left(\frac{1}{T} \be_i'\bE\bE'\hat{\bF}  +\frac{1}{T} \be_i'\bE\bB^0{\bF^0}'\hat{\bF} \right) \hat{\bLambda}^{-1}  \right\|_{\F} \\
		&& \le O_p\left( \left( \frac{N}{T}+1 \right)N^{-\alpha_r} \right)
   + 
   \left\|  \frac{1}{T} \be_i'\bE\bB^0 \bN^{-\frac{1}{2}} \right\|_{\F} \left\|\frac{1}{T}\bN^{\frac{1}{2}}{\bF^0}'\hat{\bF}\bN^{-\frac{1}{2}}\bN^{\frac{1}{2}}\hat{\bLambda}^{-1}
   \right\|_{\F} \\
   && = 
  O_p\left( \left( \frac{N}{T}+1 \right)N^{-\alpha_r} \right)
   + \left[ O_p\left(\frac{1}{\sqrt{N^{\alpha_r}}} \right)+O_p\left(\frac{1}{\sqrt{T}}\right)\right] O_p\left(N^{-\frac{1}{2}\alpha_r}\right) \\
   && =
    O_p\left(\frac{N^{1-\alpha_r}}{T}\right) 
   + O_p\left(N^{-\alpha_r}\right)+ O_p\left(\frac{1}{\sqrt{TN^{\alpha_r}}}\right),
	\end{eqnarray*}
 since $\frac{1}{T} \be_i'\bE\bB^0 \bN^{-\frac{1}{2}} =O_p\left( \frac{1}{\sqrt{N^{\alpha_r}}}\right)+O_p\left( \frac{1}{\sqrt{T}}\right)$ by Lemma \ref{lem:BEF}(ii).
\end{proof}

We next show the bound of $\tilde{\bQ} - \bI_r$ in Proposition \ref{prop:Q-I dec}.

\begin{prop}\label{prop:Q-I dec}
    Suppose that Assumptions \ref{cond:eigen}--\ref{assumption BaiNgA3} hold and  $\frac{N^{1-\alpha_r}}{T} \to  0$ as $N, T \to \infty$,  then we have 
\begin{flalign*}
&(i) ~ \|\tilde{\bQ} -\bI_r \|_{\F}=O_p(\tilde{\Delta}_{NT}) , &\\
&(ii) ~ \|\bN^{-1/2}(\hat{\bLambda}-\bLambda) \bN^{-1/2}\|_{\F} = O_p(\tilde{\Delta}_{NT}), &\\
&(iii) ~ \|\bLambda^{\frac{1}{2}}(\tilde{\bQ}' -\bI_r) \|_{\F}=o_p\left( \sqrt{\frac{N}{T}}\right) +O_p\left( N^{-\alpha_r/2}\right) ,
\end{flalign*}
where 
\begin{align}
\label{eq:TDleta}
    \tilde{\Delta}_{NT} = N^{-\alpha_r}+\frac{N^{1-\alpha_r}}{T} .
\end{align}
\end{prop}

\begin{rem}
An earlier version of this paper defined the rates in  \eqref{Delta} and \eqref{eq:TDleta} as 
$\Delta_{NT} = \frac{N^{1-\alpha_r}}{T} 
+ N^{\frac{1}{2}\alpha_1-\alpha_r}\frac{N^{1-\alpha_r}}{T} 
+ N^{\frac{1}{2}\alpha_1-\frac{3}{2}\alpha_r} 
+ \frac{N^{\frac{1}{2}\alpha_1-\alpha_r}}{\sqrt{T}}$ 
and 
$\tilde{\Delta}_{NT} = \frac{N^{1-\alpha_r}}{T} 
+ N^{-\alpha_r} 
+ \frac{N^{\frac{1}{2}\alpha_1-\alpha_r}}{\sqrt{T}}$. These earlier formulations, which contain additional cross terms, are used in
\citet[Lemma A.1]{jiang2024Mw}, which cites an earlier version of this paper.
The additional terms are dominated by those retained in the present definitions.
Consequently, the earlier and current formulations are equivalent in order, and
all results stated using the earlier expressions remain valid.
\end{rem}

\begin{proof}[Proof of Proposition \ref{prop:Q-I dec}]
Define $\hat\bU:=\hat\bF/\sqrt T$, $\bU^0:=\bF^0/\sqrt T$ and $\bQ:=\bU^{0\prime}\hat{\bU}$. Then, we have equality
$
\tilde{\bQ}:=T^{-1}\hat\bF'\bF^0=\hat\bU'\bU^0 =\bQ'.
$
Let $\hat\bA:=T^{-1}\bX\bX', 
\bA^0:=\bU^0\bLambda\bU^{0\prime},
\bDelta_{\bA}:=\hat\bA-\bA^0$. 
Combining the eigen-equations $\hat{\bA} \hat{\bU}  =  \hat{\bU} \hat{\bLambda}$ and $\bA^0\bU^0  =  \bU^0 {\bLambda}$, and premultiplying by $\bU^{0\prime}$, we obtain 
\begin{align}
    \bU^{0\prime} (\bA^0 + \bDelta _{\bA}) \hat{\bU} & = \bU^{0\prime} \hat{\bU} \hat{\bLambda} \\
   \bLambda\bU^{0\prime} \hat{\bU}+ \bU^{0\prime}  \bDelta _{\bA}\hat{\bU} & = \bU^{0\prime} \hat{\bU} \hat{\bLambda} \\
     \bU^{0\prime} \hat{\bU}+  \bLambda^{-1}\bU^{0\prime}  \bDelta _{\bA}\hat{\bU} & =  \bLambda^{-1} \bU^{0\prime} \hat{\bU} \hat{\bLambda} \\
       \bQ - \bLambda^{-1} \bQ\hat{\bLambda}&= - \bLambda^{-1}\bU^{0\prime}  \bDelta _{\bA}\hat{\bU} .  \label{eq:Q}
\end{align}

Define $\bM:=\bU^{0\prime}\bDelta_{\bA}\hat{\bU}$. For \(k\ne \ell\), the \((k,\ell)\)-th element of \eqref{eq:Q} satisfies
	\[
	\left(1-\frac{\hat\lambda_\ell}{\lambda_k}\right)q_{k\ell}
	=
	-\lambda_k^{-1}m_{k\ell}.
	\]
	Equivalently,
	\begin{equation}\label{eq:qkl-eig}
		q_{k\ell}
		=
		-\frac{m_{k\ell}}{\lambda_k-\hat\lambda_\ell}.
	\end{equation}
	By the relative eigengap condition and Lemma \ref{lemma order of lambda},
	\[
	|\lambda_k-\hat\lambda_\ell|
	\asymp
	\max\{\lambda_k,\lambda_\ell\}
	\asymp
	N^{\max(\alpha_k,\alpha_\ell)}
	\qquad (k\ne \ell).
	\]
	We now expand \(\bM\). Since
	\[
	\bDelta_{\bA}
	=
	T^{-1}\bF^0\bB^{0\prime}\bE'
	+
	T^{-1}\bE\bB^0\bF^{0\prime}
	+
	T^{-1}\bE\bE',
	\]
	 we obtain
	\begin{equation}\label{eq:M-expansion}
		\bM
		=
		\frac1T\bB^{0\prime}\bE'\bF^0 \tilde{\bH}_4
		+
        \frac1T\bB^{0\prime}\bE'(\hat{\bF}-\bF^0 \tilde{\bH}_4)
        +
		\frac1T\bF^{0\prime}\bE\bB^0\bQ
		+
		\frac1{T^2}\bF^{0\prime}\bE\bE'\hat{\bF}.
	\end{equation}
(i) We  first control the off-diagonal part and bound the contribution of the four terms in \eqref{eq:M-expansion} to
	\(q_{k\ell}\), \(k\ne\ell\). Let
	\[
	q_{\text{off}}:=\max_{k\ne \ell}|q_{k\ell}|.
	\]
    First, consider $m^{(1)}_{k\ell}
	:=\left(
	\frac1T\bB^{0\prime}\bE'\bF^0\tilde{\bH}_4
	\right)_{k\ell}$. 
	Since \(\tilde{\bH}_4=O_p(1)\) and Assumption \ref{assumption BaiNgA3} gives
	$
	\left\|
	\bN^{-1/2}\frac1{\sqrt T}\bB^{0\prime}\bE'\bF^0
	\right\|_{\F}
	=
	O_p(1),$
	this term satisfies
	\[
	\left| \left(
	\frac1T\bB^{0\prime}\bE'\bF^0\tilde{\bH}_4
	\right)_{k\ell}\right|
	=
   N^{\alpha_k/2} \left|\left(\bN^{-1/2}
	\frac1T\bB^{0\prime}\bE'\bF^0\tilde{\bH}_4
	\right)_{k\ell}\right|
    =
	O_p\left(\frac{N^{\alpha_k/2}}{\sqrt T}\right).
	\]
	Dividing by the eigengap in \eqref{eq:qkl-eig} gives
	\[
	\frac{
		O_p(N^{\alpha_k/2}/\sqrt T)
	}{
		N^{\max(\alpha_k,\alpha_\ell)}
	}
	=
	O_p\left(\frac{N^{-\alpha_r/2}}{\sqrt T}\right).
	\]
	Indeed, if \(\alpha_k\ge\alpha_\ell\), the exponent is
	\(-\alpha_k/2\le-\alpha_r/2\); if \(\alpha_k<\alpha_\ell\), the exponent is
	\(\alpha_k/2-\alpha_\ell\le-\alpha_\ell/2\le-\alpha_r/2\). Moreover,  by $2\sqrt{ab} \le a+b$, this term is dominated by $\frac{N^{1-\alpha_r}}{T} + \frac{1}{N^{\alpha_r}}$ since $\frac{N^{-\alpha_r/2}}{\sqrt T}\le \sqrt{\frac{N^{1-\alpha_r}}{T}}\frac{1}{\sqrt{N^{\alpha_r}}}\le \frac{1}{2}\left( \frac{N^{1-\alpha_r}}{T} + \frac{1}{N^{\alpha_r}}\right)$.
    
	For the second term $m^{(2)}_{k\ell}
	:=\left(
	\frac1T\bB^{0\prime}\bE'
	(\hat{\bF}-\bF^0\tilde{\bH}_4)
	\right)_{k\ell}$, Lemma \ref{lemma f(fh-f0)}(ii) implies
	\[
	\left| \left(
	\frac1T\bB^{0\prime}\bE'
	(\hat{\bF}-\bF^0\tilde{\bH}_4)
	\right)_{k\ell} \right|
	\le 
    N^{\alpha_k /2}\left\|
	\frac1T\bN^{-1/2}\bB^{0\prime}\bE'
	(\hat{\bF}-\bF^0\tilde{\bH}_4)
	\right\|_{\F}
    =
	O_p\left(
	N^{\alpha_k/2}
	\left\{
	N^{-\alpha_r/2}
	+
	\frac{N^{1-\alpha_r}}{T}
	\right\}
	\right).
	\]
	After division by \(N^{\max(\alpha_k,\alpha_\ell)}\), the second term is bounded by
	\[
	O_p(N^{-\alpha_r})
	+
	O_p\left(\frac{N^{1-\frac{3}{2}\alpha_r}}{T}\right).
	\]
	Next, consider
	\begin{align*}
	    m^{(3)}_{k\ell}
   & =
	\left(\frac1T\bF^{0\prime}\bE\bB^0\bN^{-1/2}\bN^{1/2}\bQ\right)_{k\ell}\\
	&=
	   \left(\frac1T\bF^{0\prime}\bE\bB^0\bN^{-1/2}\right)_{k\ell}N^{\alpha_\ell/2}q_{\ell\ell}
	+
	\sum_{j\ne\ell}   \left(\frac1T\bF^{0\prime}\bE\bB^0\bN^{-1/2}\right)_{kj}N^{\alpha_j/2}q_{j\ell}.
	\end{align*}
	Since \(\|\bQ_{\cdot\ell}\|_2\le1\), we have \(|q_{\ell\ell}|\le1\). Thus the first term is
	\begin{align}
\label{eq:m_2}\left(\frac1T\bF^{0\prime}\bE\bB^0\bN^{-1/2}\right)_{k\ell}N^{\alpha_\ell/2}q_{\ell\ell}
	=
	O_p\left(\frac{N^{\alpha_\ell/2}}{\sqrt T}\right).
	\end{align}

	After division by \(N^{\max(\alpha_k,\alpha_\ell)}\), it is bounded by
	$
	O_p\left(\frac{N^{-\alpha_r/2}}{\sqrt T}\right).
	$

	It remains to handle the sum over \(j\ne\ell\):
	\[
	\left|
\sum_{j\ne\ell}\left(\frac1T\bF^{0\prime}\bE\bB^0\bN^{-1/2}\right)_{kj}N^{\alpha_j/2}q_{j\ell}
	\right|
	\le
	O_p\left(\frac{N^{\alpha_1/2}}{\sqrt T}\right)q_{\text{off}}.
	\]
	After division by the eigengap, this term is bounded by
	$
	O_p\left(\frac{N^{\alpha_1/2-\alpha_r}}{\sqrt T}\right)q_{\text{off}}$.	Since \(N^{1-\alpha_r}/T\to0\) and \(\alpha_1\le1\),
	$
	\frac{N^{\alpha_1/2-\alpha_r}}{\sqrt T}
	=
	o(1)$.
	Thus , we obtain $ 
		O_p\left(\frac{N^{-\alpha_r/2}}{\sqrt T}\right)
		+
		o_p(1)q_{\text{off}}$.
	
	We next consider $	m^{(4)}_{k\ell}
	:=
	\left(
	\frac1{T^2}\bF^{0\prime}\bE\bE'\hat{\bF}
	\right)_{k\ell}$. We can see that
	\[
	| m^{(4)}_{k\ell}|
	\le
	\frac1{T^2}\|\bF^0\|_{\F}\|\bE\|_2^2\|\hat{\bF}\|_{\F}
	=
	O_p\left(1+\frac NT\right).
	\]
	Dividing by the eigengap gives
	\begin{equation}\label{eq:m3-contribution}
		O_p(N^{-\alpha_r})
		+
		O_p\left(\frac{N^{1-\alpha_r}}{T}\right).
	\end{equation}
	
	Combining the above bounds, we obtain, uniformly over \(k\ne\ell\),
	\[
	|q_{k\ell}|
	\le
	O_p(\tilde{\Delta}_{NT})	+
	o_p(1)q_{\text{off}}.
	\]
    Thus, $q_{\text{off}}=O_p(\tilde{\Delta}_{NT})$.
	Since \(r\) is fixed, this implies
	\begin{equation}\label{eq:off-Q-unweighted}
		\|\operatorname{off}(\bQ)\|_{\F}
		=
		O_p(\tilde{\Delta}_{NT}),
	\end{equation}
		where \(\operatorname{off}(\bQ)\) denotes the matrix obtained from \(\bQ\) by setting its diagonal entries equal to zero. It remains to control the diagonal part. Let $\bU_\perp\in\mathbb R^{T\times(T-r)}$ be an orthonormal basis of the orthogonal
complement of $\mathrm{span}(\bU^0)$ such that
\[
[\,\bU^0,\;\bU_\perp\,]\in\mathbb R^{T\times T}
\quad\text{is orthogonal},\qquad
\bU_\perp'\bU^0=\bzero,\quad \bU_\perp'\bU_\perp=\bI_{T-r}.
\]
Define the associated matrix
$\bS :=\bU_\perp'\hat\bU$.
Using the decomposition $\bU^0 \bU^{0\prime} + \bU_\perp \bU_\perp'=\bI_T$, we can write 
\[
\hat\bU = (\bU^0 \bU^{0\prime} + \bU_\perp \bU_\perp') \hat{\bU} = \bU^{0}\bQ+\bU_\perp\bS.
\]
Hence,
\begin{align}
\label{eq:SS}
    \hat{\bU}'\hat{\bU}
=
\bQ'\bQ+\bS'\bS
=
\bI_r.
\end{align}
Taking the \(k\)-th diagonal element of \eqref{eq:SS} gives
	\[
	q_{kk}^2+\sum_{j\ne k}q_{jk}^2+\|\bS_{\cdot k}\|_2^2=1.
	\]
	Since \(q_{kk}\ge0\) by the sign convention,
	\[
	1-q_{kk}
	=
	\frac{1-q_{kk}^2}{1+q_{kk}}
	=
	\frac{\sum_{j\ne k}q_{jk}^2+\|\bS_{\cdot k}\|_2^2}{1+q_{kk}}
	\le
	\sum_{j\ne k}q_{jk}^2+\|\bS_{\cdot k}\|_2^2.
	\]
	Thus
	\begin{equation}\label{eq:diag-basic-Q}
		|1-q_{kk}|
		\le
		\sum_{j\ne k}q_{jk}^2+\|\bS_{\cdot k}\|_2^2.
	\end{equation}
We next focus on the term $\bS$. From the eigenvector equation,
\[
(\bA^0+\bDelta_{\bA})\hat\bU = \hat\bU \hat\bLambda,
\]
premultiplying by $\bU_\perp'$ and the identity $\bU_\perp' \bA^0=\mathbf{0}$ give
\begin{align}
\label{eq:S1}
    \bU_\perp'\bDelta_{\bA}\,\hat\bU &= (\bU_\perp'\hat\bU)\hat\bLambda = \bS\hat\bLambda,\\
 \bS &= \bU_\perp'\bDelta_{\bA}\,\hat\bU\,\hat\bLambda^{-1}.
\end{align}
Substituting the expression of $\bDelta_{\bA}$ gives
\[
\bS
=
\frac{1}{\sqrt T}\bU_\perp'\bE\bB^0\bQ\hat\bLambda^{-1}
+
\frac{1}{T^{3/2}}\bU_\perp'\bE\bE'\hat\bF\hat\bLambda^{-1}.
\]
Since $\bU_\perp'\bU_\perp=I_{T-r}$, we have
\begin{align}
\label{eq:S2}
  \| \bS \|_{\F} 
 &  \le  \left\|\frac{1}{\sqrt{T}} \bU_\perp'\bE\bB^0 \bQ \hat\bLambda^{-1} \right\|_{\F}
   +
\left\|\frac{1}{T^{3/2}} \bU_\perp'\bE\bE'\hat{\bF} \hat\bLambda^{-1}\right\|_{\F}\\
 &  \le  \| \bU_\perp'\|_2  \left\|\frac{1}{\sqrt{T}}(\bE\bB^0 \bN^{-\frac{1}{2}})(\bN^{\frac{1}{2}}\bQ\bN^{-\frac{1}{2}})\bN^{\frac{1}{2}} \hat\bLambda^{-1}\right\|_{\F}
   +
 \| \bU_\perp'\|_2 \left\|\frac{1}{T^{3/2}}\bE\bE'\hat{\bF} \hat\bLambda^{-1}\right\|_{\F}\\
&=
O_p(N^{-\frac{1}{2}\alpha_r})
+
O_p\left(\frac{N^{1-\alpha_r}}{T} +N^{-\alpha_r} \right),
\end{align}
where we used $ \| \bU_\perp'\|_2=1$, $\|T^{-1/2}\bE\bB^0\bN^{-1/2}\|_{\F}=O_p(1)$, and
$\|\bN^{1/2}\bQ\bN^{-1/2}\|_{\F}=O_p(1)$.
Consequently,
\begin{align}
\| \bS'\bS\|_{\F} \le \|\bS\|_F^2
=
&
O_p(N^{-\alpha_r})
+ o_p\left(\frac{N^{1-\alpha_r}}{T}  \right)  .
\label{eq:SS_lead}
\end{align}

	Using \eqref{eq:off-Q-unweighted} and \eqref{eq:diag-basic-Q}, uniformly over \(k\),
	\[
	|1-q_{kk}|
	\le
	O_p(\tilde{\Delta}_{NT}^2)
	+
	O_p(N^{-\alpha_r})
	+
	o_p\left(\frac{N^{1-\alpha_r}}{T}\right).
	\]
	Since \(\tilde{\Delta}_{NT}\to0\),
	\[
	\max_{1\le k\le r}|1-q_{kk}|
	=
	O_p(\tilde{\Delta}_{NT}).
	\]
	Because \(r\) is fixed,
	\begin{equation}\label{eq:diag-Q-unweighted}
		\|\operatorname{Diag}(\bQ)-\bI_r\|_{\F}
		=
		O_p(\tilde{\Delta}_{NT}).
	\end{equation}
	
	Finally,
	\[
	\bQ-\bI_r
	=
	\operatorname{off}(\bQ)
	+
	\{\operatorname{Diag}(\bQ)-\bI_r\}.
	\]
    Because 
	\(\operatorname{off}(\bQ)\) has zero diagonal entries and
	\(\operatorname{Diag}(\bQ)-\bI_r\) is diagonal, the two matrices are orthogonal under the Frobenius inner product, and
	\[
	\|\bQ-\bI_r\|_{\F}^2
	=
	\|\operatorname{off}(\bQ)\|_{\F}^2
	+
	\|\operatorname{Diag}(\bQ)-\bI_r\|_{\F}^2.
	\]
	Combining \eqref{eq:off-Q-unweighted} and \eqref{eq:diag-Q-unweighted} gives
	\[
	\|\bQ-\bI_r\|_{\F}
	=
	O_p(\tilde{\Delta}_{NT}).
	\]

(ii) The bounds used in the proof of part (i) also imply
	\begin{equation}\label{eq:Lambda-inv-M-bound}
		\|\bLambda^{-1}\bM\|_{\F}
		=
		O_p(\tilde{\Delta}_{NT}).
	\end{equation}
Moreover, (i) gives
\[
q_{kk}=1+O_p(\tilde\Delta_{NT})=1+o_p(1).
\]
Taking the \((k,k)\)-th element of \eqref{eq:Q} yields
\[
q_{kk}\left(1-\frac{\hat\lambda_k}{\lambda_k}\right)
=-(\bLambda^{-1}\bU^{0\prime}\bDelta_{\bA}\hat \bU)_{kk}.
\]
Therefore,
\[
\left|
\frac{\hat\lambda_k}{\lambda_k}-1
\right|
=
\left|\frac{(\bLambda^{-1}\bU^{0\prime}\bDelta_{\bA}\hat \bU)_{kk}}{q_{kk}}\right|
=
O_p(\tilde\Delta_{NT}),
\]
uniformly over \(k=1,\ldots,r\). Since
\(\lambda_k=d_k^2N^{\alpha_k}\), the desired result follows.

(iii)	For \(k \ne \ell\), \eqref{eq:Q} implies that the \((k,\ell)\)-th element satisfies
	\[
	\left(1-\frac{\hat\lambda_\ell}{\lambda_k}\right)q_{k\ell} = -\lambda_k^{-1} \left(\bU^{0\prime}\bDelta_{\bA}\hat{\bU}\right)_{k\ell}.
	\]
    By Lemma \ref{lemma order of lambda}, 
$\hat\lambda_\ell=O_p(N^{\alpha_\ell})$ and $\lambda_\ell=O_p(N^{\alpha_\ell})$ 
for \(\ell=1,\ldots,r\). Hence, under Assumption 3, for \(k\ne\ell\),
\[
\frac{\lambda_k}{\lambda_k-\hat\lambda_\ell}
=
O_p(1)
\quad \text{or} \quad o_p(1),
\]
depending on the relative orders of \(\lambda_k\) and \(\lambda_\ell\). Hence, for the off-diagonal part,
	\[
	\left\| \bLambda^{1/2}\operatorname{off}(\bQ) \right\|_{\F} = O_p\left( \left\| \bLambda^{-1/2}\bU^{0\prime}\bDelta_{\bA}\hat{\bU} \right\|_{\F} \right),
	\]
	where \(\operatorname{off}(\bQ)\) denotes the matrix obtained from \(\bQ\) by setting its diagonal entries equal to zero.
	
	Now expanding \(\bDelta_{\bA}\), we have
	\begin{align}
		\left\| \bLambda^{-1/2}\bU^{0\prime}\bDelta_{\bA}\hat{\bU} \right\|_{\F}
		&\le
		\left\| \bLambda^{-1/2}\frac{1}{T}\bB^{0\prime}\bE'\bF^0\tilde{\bH}_4 \right\|_{\F} 
	+
		\left\| \bLambda^{-1/2}\frac{1}{T}\bB^{0\prime}\bE'(\hat{\bF}-\bF^0\tilde{\bH}_4) \right\|_{\F}\\
		&\quad+
		\left\| \bLambda^{-1/2}\frac{1}{T}\bF^{0\prime}\bE\bB^0\bQ \right\|_{\F} +
		\left\| \bLambda^{-1/2}\frac{1}{T^2}\bF^{0\prime}\bE\bE'\hat{\bF} \right\|_{\F}\\
		&=
		O_p\left(\frac{1}{\sqrt T}\right)
		+
		O_p(N^{-\alpha_r/2})
		+
		O_p\left(\frac{N^{1-\alpha_r}}{T}\right) \\
		&\quad+
		O_p\left( \frac{N^{(\alpha_1-\alpha_r)/2}}{\sqrt T} \right)
		+
		O_p\left(\frac{N^{1-\alpha_r/2}}{T}\right) ,
	\end{align}
	where we have used \(\bLambda=\bD^2\bN\), \(\tilde{\bH}_4=O_p(1)\), Lemma \ref{lemma f(fh-f0)}(ii), and 
	\[
	\left\| \bN^{-1/2}\frac{1}{\sqrt T}\bB^{0\prime}\bE'\bF^0 \right\|_{\F} = O_p(1)
	\]
in Assumption 4(vi). In addition, since \(\alpha_1-\alpha_r\le1\) and \(N^{1-\alpha_r}/T\to0\), we have
	\[
	\frac{N^{(\alpha_1-\alpha_r)/2}}{\sqrt T}
	\lesssim
	\sqrt{\frac NT},
	\quad
	\frac{N^{1-\alpha_r/2}}{T}
	=
	\sqrt{\frac NT}\sqrt{\frac{N^{1-\alpha_r}}{T}}
	=
	o\left(\sqrt{\frac NT}\right).
	\]
	Thus
	\begin{align}
		\label{eq:off-Q-final}
	\left\|
	\bLambda^{1/2}\operatorname{off}(\bQ)
	\right\|_{\F}
	=
	o_p\left(\sqrt{\frac NT}\right)
	+
	O_p\left(N^{-\alpha_r/2}\right).
	\end{align}
			It remains to control the diagonal part:
	\[
	1-q_{kk}
	\le
	\sum_{j\ne k}q_{jk}^2+\|\bS_{\cdot k}\|_2^2.
	\]
Premultiplied by $\lambda_k^{1/2}$, that is
	\begin{equation}\label{eq:diag-q-basic}
		\lambda_k^{1/2}|1-q_{kk}|
		\le
		\lambda_k^{1/2}\sum_{j\ne k}q_{jk}^2
		+
		\lambda_k^{1/2}\|\bS_{\cdot k}\|_2^2.
	\end{equation}
	
	We first bound \(\bS_{\cdot k}\). For each \(k\),
	\[
	\bS_{\cdot k}
	=
	\hat\lambda_k^{-1}
	\bU_\perp'\bDelta_{\bA}\hat{\bU}_k.
	\]
	Using the expansion of \(\bDelta_{\bA}\) and \(\bU_\perp'\bF^0=0\), the term
	\(T^{-1}\bF^0\bB^{0\prime}\bE'\) drops out, and hence
	\[
	\bS_{\cdot k}
	=
	\hat\lambda_k^{-1}
	\left[
	\frac1{\sqrt T}\bU_\perp'\bE\bB^0\bQ_{\cdot k}
	+
	\frac1T\bU_\perp'\bE\bE'\hat{\bU}_k
	\right],
	\]
	because \(\bF^{0\prime}\hat{\bU}_k=\sqrt T\,\bQ_{\cdot k}\).
	Write
	\[
	\bB^0\bQ_{\cdot k}
	=
	\bB_k^0 q_{kk}
	+
	\sum_{j\ne k}\bB_j^0q_{jk}.
	\]
	Since \(\|\bQ_{\cdot k}\|_2\le1\), we have \(|q_{kk}|\le1\). Moreover,
	\[
	\|\bE\bB_j^0N^{-\alpha_j/2}\|_2
	=
	O_p(\sqrt{T}),
	\qquad j=1,\ldots,r,
	\]
	because
	\[
	\|\bE\bB^0\bN^{-1/2}\|_{\F}=O_p(\sqrt T).
	\] 
	Therefore,
	\begin{align}
		\left\|
		\frac1{\sqrt T}\bU_\perp'\bE\bB^0\bQ_{\cdot k}
		\right\|_2
		&\le
		O_p(N^{\alpha_k/2})
		+
		\sum_{j\ne k}O_p(N^{\alpha_j/2})|q_{jk}|. 
		\label{eq:S-first-column}
	\end{align}
	By \eqref{eq:off-Q-final},
	\[
	\sum_{j\ne k}N^{\alpha_j/2}|q_{jk}|
	=
	o_p\left(\sqrt{\frac NT}\right)
	+
	O_p\left(N^{-\alpha_r/2}\right).
	\]
	Thus
	\[
	\left\|
	\frac1{\sqrt T}\bU_\perp'\bE\bB^0\bQ_{\cdot k}
	\right\|_2
	=
	O_p(N^{\alpha_k/2})+
	o_p\left(\sqrt{\frac NT}\right)
	+
	O_p\left(N^{-\alpha_r/2}\right).
	\]
	Since \(\hat\lambda_k\asymp_p N^{\alpha_k}\),  \(N^{1-\alpha_r}/T\to0\) and \(\alpha_r\le\alpha_k\), we have
	\[
	\hat\lambda_k^{-1}
	\left\|
	\frac1{\sqrt T}\bU_\perp'\bE\bB^0\bQ_{\cdot k}
	\right\|_2
	=
	O_p(N^{-\alpha_k/2}).
	\]
	For the remaining term,
	\[
	\begin{aligned}
		\left\|
		\hat\lambda_k^{-1}
		\frac1T\bU_\perp'\bE\bE'\hat{\bU}_k
		\right\|_2
		&\le
		\hat\lambda_k^{-1}\frac1T\|\bE\|_2^2\|\hat{\bU}_k\|_2 \\
		&=
		O_p(N^{-\alpha_k})
		+
		O_p\left(\frac{N^{1-\alpha_k}}{T}\right).
	\end{aligned}
	\]
	Consequently,
	\begin{equation}\label{eq:S-column-final}
		\|\bS_{\cdot k}\|_2
		=
		O_p(N^{-\alpha_k/2})
		+
		O_p\left(\frac{N^{1-\alpha_k}}{T}\right).
	\end{equation}
	
	Now consider the two terms on the right-hand side of \eqref{eq:diag-q-basic}. 
	For the term involving \(\bS_{\cdot k}\), \eqref{eq:S-column-final} gives
	\begin{align}
		\lambda_k^{1/2}\|\bS_{\cdot k}\|_2^2
		&=
		O_p(N^{-\alpha_k/2})
		+
		o_p\left(\frac{N^{1-\alpha_k/2}}{T}\right).
		\label{eq:S-diag-contribution}
	\end{align}
		Next consider the off-diagonal squared term in \eqref{eq:diag-q-basic}. For \(j\ne k\),
	the \((j,k)\)-th element of the eigen-equation gives
	\[
	q_{jk}
	=
	-\frac{	1}{\lambda_j-\hat\lambda_k} \left[
	\bU^{0\prime}\bDelta_{\bA}\hat{\bU}
	\right]_{jk}.
	\]
	By the relative eigengap argument,
	\[
	|\lambda_j-\hat\lambda_k|
	\asymp_p
	N^{\max(\alpha_j,\alpha_k)},
	\qquad j\ne k.
	\]
	Hence
	\[
	\lambda_k^{1/2}q_{jk}^2
	\le
	O_p(N^{-\alpha_k/2})
	\left(\lambda_j^{-1/2} \left[
	\bU^{0\prime}\bDelta_{\bA}\hat{\bU}
	\right]_{jk} \right)^2 .
	\]
	Summing over \(j\ne k\) and using
	\[
	\|\bLambda^{-1/2}\bU^{0\prime}\bDelta_{\bA}\hat{\bU}\|_{\F}
	=
o_p\left(\sqrt{\frac NT}\right)
+
O_p\left(N^{-\alpha_r/2}\right).
	\]
	and $\sqrt{\frac NT}
	+
	N^{-\alpha_r/2} \lesssim N^{\alpha_k/2}$,
	we obtain
	\begin{align}
		\lambda_k^{1/2}\sum_{j\ne k}q_{jk}^2
		&=
		o_p\left(\sqrt{\frac NT}\right)
		+
		O_p(N^{-\alpha_r/2}).
		\label{eq:q-off-square-contribution}
	\end{align}
	
	Combining \eqref{eq:diag-q-basic}, \eqref{eq:S-diag-contribution}, and
	\eqref{eq:q-off-square-contribution}, we have, uniformly over \(k=1,\ldots,r\),
	\[
	\lambda_k^{1/2}|1-q_{kk}|
	=
	o_p\left(\sqrt{\frac NT}\right)
	+
	O_p(N^{-\alpha_r/2}).
	\]
	Since \(r\) is fixed, this implies
	\begin{equation}\label{eq:diag-Q-final}
		\left\|
		\bLambda^{1/2}
		\{\operatorname{Diag}(\bQ)-\bI_r\}
		\right\|_{\F}
		=
		o_p\left(\sqrt{\frac NT}\right)
		+
		O_p(N^{-\alpha_r/2}).
	\end{equation}
	Finally, decompose
	\[
	\bQ-\bI_r
	=
	\operatorname{off}(\bQ)
	+
	\{\operatorname{Diag}(\bQ)-\bI_r\},
	\]
	then,
	\[
	\left\|
	\bLambda^{1/2}(\bQ-\bI_r)
	\right\|_{\F}^2
	=
	\left\|
	\bLambda^{1/2}\operatorname{off}(\bQ)
	\right\|_{\F}^2
	+
	\left\|
	\bLambda^{1/2}\{\operatorname{Diag}(\bQ)-\bI_r\}
	\right\|_{\F}^2.
	\]
	Using \eqref{eq:off-Q-final} and \eqref{eq:diag-Q-final}, we obtain
	\[
	\left\|
	\bLambda^{1/2}(\bQ-\bI_r)
	\right\|_{\F}
	=
	o_p\left(\sqrt{\frac NT}\right)
	+
	O_p\left(N^{-\alpha_r/2}\right).
	\]
	Since \(\bQ=\tilde{\bQ}'\), the proof is complete.
\end{proof}

For later use, define
\begin{align*}
   &  \tilde{\bH}_1
    :=
    \bLambda(\hat{\bB}'\bB^0)^{-1},
    \qquad
    \tilde{\bH}_2
    :=
    T^{-1}{\bF^0}'\hat{\bF}
    =
    \tilde{\bQ}', \\
  &  \tilde{\bH}_3
    :=
    \left(T^{-1}\hat{\bF}'\bF^0\right)^{-1},
    \qquad
    \tilde{\bH}_4
    :=
    {\bB^0}'\hat{\bB}\hat{\bLambda}^{-1},\\
  &     \tilde{\bH}
    :=
    {\bB^0}'\bB^0
    \left(T^{-1}{\bF^0}'\hat{\bF}\right)
    \hat{\bLambda}^{-1}. 
\end{align*}

\begin{lem} Suppose that Assumptions \ref{cond:eigen}--\ref{assumption BaiNgA3} hold. Then, we have,
\label{lem:H equality}
\begin{flalign*}
     &(i)~{\tilde{\bH}_4} -\tilde{\bH}=\frac{1}{T} {\bB^0}' \bE^{\prime} \hat{\bF} \hat{\bLambda}^{-1}, &\\
    & (ii)~ {\tilde{\bH}_2} -\tilde{\bH}_4= \frac{1}{T}{\bF}^{0'}\bE \hat{\bB} \hat{\bLambda}^{-1}, &\\
  & (iii)~  {\tilde{\bH}_3}  -\tilde{\bH}_1= \left( \frac{1}{T} \hat{\bF}'\bF^0\right)^{-1}\frac{1}{T}  \hat{\bF}^{'}\bE {\bB}^0 ({\hat{\bB}}'{\bB^0})^{-1} , &\\
   & (iv) ~{\tilde{\bH}_3}  -\tilde{\bH}_4= \left( \frac{1}{T} \hat{\bF}'\bF^0\right)^{-1}\frac{1}{T}  \hat{\bF}^{'}\bE \hat{\bB} \hat{\bLambda}^{-1},  &\\
   &(v) ~{\tilde{\bH}_1}^{'-1} -\tilde{\bH}_2= \frac{1}{T}\bLambda^{-1}{\bB}^{0'}\bE'\hat{\bF},  &\\
   & (vi)~ {\tilde{\bH}_4}^{'-1} -\tilde{\bH}_2= \frac{1}{T}({\hat{\bB}}'{\bB^0})^{-1}  \hat{\bB}'  \bE'\hat{\bF}. &
\end{flalign*}
\end{lem}

\begin{proof}(i) We proceed by left multiplying ${\bB^0}'$ and right multiplying $\hat{\bLambda}^{-1}$ to the first equation, to get
    \begin{eqnarray*}
\hat{\bB}&=&\frac{1}{T} \bX^{\prime} \hat{\bF}=\frac{1}{T} \bB^0{\bF^0}'  \hat{\bF}+\frac{1}{T} \bE^{\prime} \hat{\bF} \\
{\bB^0}' \hat{\bB}&=&\frac{1}{T} {\bB^0}' \bB^0 {\bF^0}' \hat{\bF}+\frac{1}{T} {\bB^0}' \bE^{\prime} \hat{\bF} \\
 {\bB^0}' \hat{\bB} \hat{\bLambda}^{-1}&=&\frac{1}{T} {\bB^0}' \bB^0 {\bF^0}' \hat{\bF}\hat{\bLambda}^{-1}+\frac{1}{T} {\bB^0}' \bE^{\prime} \hat{\bF} \hat{\bLambda}^{-1} \\
\tilde{\bH}_4&=&\tilde{\bH}+\frac{1}{T} {\bB^0}' \bE^{\prime} \hat{\bF} \hat{\bLambda}^{-1}.
\end{eqnarray*}

(ii) By the definition of $\hat{\bB}$ and expanding $\bX$, we obtain
\begin{eqnarray*}
\hat{\bF}&=& \frac{1}{T}\bX\bX'\hat{\bF}\hat{\bLambda}^{-1}=\bX  \hat{\bB}\hat{\bLambda}^{-1} =\bF^0 {\bB^0}'\hat{\bB}\hat{\bLambda}^{-1}+\bE \hat{\bB}\hat{\bLambda}^{-1} \\
\frac{1}{T} {\bF^0}'\hat{\bF} &=& {\bB^0}'\hat{\bB}\hat{\bLambda}^{-1}+\frac{1}{T} {\bF^0}'\bE \hat{\bB}\hat{\bLambda}^{-1} .
\end{eqnarray*}

(iii) ${\bB^0}' \hat{\bB}$ in (i) implies
 \begin{eqnarray*}
\hat{\bB}'\bB^0 &=&\frac{1}{T}\hat{\bF}' {\bF^0}{\bB^0}' \bB^0 +\frac{1}{T} \hat{\bF}'\bE{\bB^0}  \\
 \tilde{\bH}_3 &=&\tilde{\bH}_1+ \left( \frac{1}{T}\hat{\bF}' {\bF^0}\right)^{-1} \frac{1}{T} \hat{\bF}'\bE{\bB^0} \left(\hat{\bB}'\bB^0 \right)^{-1} .
\end{eqnarray*}

(iv) Multiplying $\hat{\bB}'$ to the first equation in (i),
 \begin{eqnarray*}
{\hat{\bB}}' \hat{\bB}&=&\frac{1}{T} {\hat{\bB}}' \bB^0 {\bF^0}' \hat{\bF}+\frac{1}{T} {\hat{\bB}}' \bE^{\prime} \hat{\bF} \\
{\hat{\bB}}' \hat{\bB}&=&\frac{1}{T}\hat{\bF}'\bF^0 {\bB^0}' \hat{\bB} +\frac{1}{T} \hat{\bF}' \bE{\hat{\bB}}  \\
  \tilde{\bH}_3 &=&\tilde{\bH}_4+ \left( \frac{1}{T}\hat{\bF}' {\bF^0}\right)^{-1}\frac{1}{T} \hat{\bF}' \bE{\hat{\bB}}\left({\hat{\bB}}' \hat{\bB}\right)^{-1}.
\end{eqnarray*}

(v) Post-multiplying $\bB^0$ to the transpose of the first equation in (i),
\begin{eqnarray*}
\hat{\bB}'\bB^0 &=&\frac{1}{T}\hat{\bF}' {\bF^0}{\bB^0}' \bB^0 +\frac{1}{T} \hat{\bF}'\bE{\bB^0} \\
  \tilde{\bH}_1^{-1}&=&\tilde{\bH}'_2+\frac{1}{T} \hat{\bF}'\bE{\bB^0} \left({\bB^0}' \bB^0\right)^{-1}\\
 \tilde{\bH}_1^{'-1}&=&\tilde{\bH}_2 + \left({\bB^0}' \bB^0\right)^{-1} \frac{1}{T} {\bB^0}' \bE^{\prime} \hat{\bF}  .
\end{eqnarray*}

(vi) The first equation in (iv) implies
\begin{eqnarray*}
\tilde{\bH}_4^{'-1}&=&\tilde{\bH}_2 + \left(\hat{\bB}' \bB^0\right)^{-1} \frac{1}{T} \hat{\bB}' \bE^{\prime} \hat{\bF}  .
\end{eqnarray*}
\end{proof}

\begin{lem} Suppose that Assumptions \ref{cond:eigen}--\ref{assumption BaiNgA3} hold. If $\frac{N^{1-\alpha_r}}{T}\to 0$, then, we have,
    \begin{flalign*}
  &(i)~ \left\| \frac{1}{T}  {\bF^0}'\bE\hat{\bB}({\hat{\bB}}'\hat{\bB})^{-1}\right\|_{\F} = O_p(\Delta_1)  ,&\\
   & (ii)~ \left\| \frac{1}{T}  \hat{\bF}^{'}\bE {\bB}^0 ({\hat{\bB}}'{\bB^0})^{-1}\right\|_{\F} = O_p(\Delta_2) ,&\\
    &(iii) ~\left\| \frac{1}{T}  \hat{\bF}^{'}\bE \hat{\bB} ({\hat{\bB}}'\hat{\bB})^{-1}\right\|_{\F} = O_p(\Delta_1+\Delta_2) .&
 \end{flalign*}
where 
\begin{flalign}
    & \Delta_1= \sqrt{\frac{N^{1-\alpha_r}}{T} }
  N^{-\alpha_r} +\frac{N^{1-\alpha_r}}{T}+\frac{1}{\sqrt{TN^{\alpha_r}}} ,\label{Delta_1} &\\
   &\Delta_2=\frac{N^{1-\alpha_r}}{T}N^{-\frac{1}{2} \alpha_r}+N^{- \alpha_r}+\frac{1}{\sqrt{T N^{\alpha_r}}} .\label{Delta_2}&
\end{flalign}
\label{lem:FEB0}
\end{lem}

\begin{proof}
(i) 
Using the results in the proof of Lemma \ref{lem fh4bq}, we show that
    \begin{align*}
& \left\|\frac{1}{T} \bF^{0^{\prime}} \bE \hat{\bB}(\hat{\bB} \hat{\bB})^{-1}\right\|_{\F} \\
&\le \left\|\frac{1}{T} {\bF^0}' \bE\left(\hat{\bB}-\bB^0 \tilde{\bQ}^{\prime}\right)\left(\hat{\bB}^{\prime} \hat{\bB}\right)^{-1}\right\|_{\F}+\left\|\frac{1}{T} {\bF^0}' \bE  \bB^0 \tilde{\bQ}^{\prime}\left(\hat{\bB}^{\prime} \hat{\bB}\right)^{-1}\right\|_{\F} \\
& \le  \frac{1}{T} O_p(\sqrt{N T}) O_p\left(  
\sqrt{\frac{N^{1-\alpha_r}}{T}} \frac{N^{1-\frac{1}{2} \alpha_r}}{T}+\frac{N^{1-\frac{1}{2} \alpha_r}}{T}+N^{-\frac{1}{2} \alpha_r}+\sqrt{\frac{N}{T}}
 \right)O_p( N^{-\alpha_r}) +O_p\left(
\frac{1}{\sqrt{TN^{\alpha_r}}} \right)\\
& =  
O_p\left( \sqrt{\frac{N^{1-\alpha_r}}{T} }N^{-\alpha_r}\right) 
+O_p\left( \frac{N^{1-\alpha_r}}{T}\right)
+O_p\left( \frac{1}{\sqrt{TN^{\alpha_r}}}\right)   .
\end{align*}

(ii) Lemma \ref{lemma f(fh-f0)}(ii) implies
\begin{align*}
& \left\|\frac{1}{T} \hat{\bF}^{\prime} \bE \bB^0\left(\hat{\bB}^{\prime} \bB^0\right)^{-1}\right\|_{\F} \\
&\le
\left\|\frac{1}{T}(\hat{\bF}-\bF^0 \tilde{\bH}_4)^{\prime} \bE \bB^0\left(\hat{\bB}' \bB^0\right)^{-1}\right\|_{\F}+\left\|\frac{1}{T} \tilde{\bH}_4^{\prime}\left({\bF^0}^{\prime} \bE \bB^0\right)\left(\hat{\bB}' \bB^0\right)^{-1}\right\|_{\F} \\ 
& \le  O_p\left(  \left( \frac{N^{1-\alpha_r}}{T}+N^{-\frac{1}{2}\alpha_r}\right)    N^{-\frac{1}{2}\alpha_r}\right)
+O_p\left(\frac{\sqrt{T N^{\alpha_1}}}{T} N^{-\alpha_r} \right)\\
& = O_p\left(\frac{N^{1-\alpha_r}}{T}N^{-\frac{1}{2} \alpha_r}\right)
+O_p\left(N^{-\alpha_r}\right)
+O_p\left(\frac{1}{\sqrt{T N^{\alpha_r}}}\right),
\end{align*}
where we have used
\[
    \left\|
    \bN^{1/2}(\hat{\bB}'\bB^0)^{-1}
    \right\|_{\F}
    =
    \left\|
    \bN^{1/2}\bLambda^{-1}\tilde{\bH}_1  
    \right\|_{\F}
    =
    O_p(N^{-\alpha_r/2}),
\]
since $\tilde{\bH}_1^{-1}-\tilde{\bH}_2=o_p(1)$ and $\tilde{\bH}_1=O_p(1)$.

(iii) Lemma \ref{lemma f(fh-f0)}(i) implies
\begin{align*}
& \left\|\frac{1}{T} \hat{\bF}^{\prime} \bE \hat{\bB} \left(\hat{\bB}^{\prime} \hat{\bB}\right)^{-1}\right\|_{\F} \\
& \le \left\|\frac{1}{T}(\hat{\bF}-\bF^0 \tilde{\bH}_4)^{\prime} \bE \hat{\bB} \left(\hat{\bB}' \hat{\bB}\right)^{-1}\right\|_{\F}+\left\|\frac{1}{T} \tilde{\bH}_4^{\prime}\left({\bF^0}^{\prime} \bE \hat{\bB} \right)\left(\hat{\bB}' \hat{\bB}\right)^{-1}\right\|_{\F} \\ 
& = \left\|\frac{1}{T}(\hat{\bF}-\bF^0 \tilde{\bH}_4)^{\prime} \bE (\hat{\bB}-\bB^0\tilde{\bQ}) \left(\hat{\bB}' \hat{\bB} \right)^{-1}+ \frac{1}{T}(\hat{\bF}-\bF^0 \tilde{\bH}_4)^{\prime} \bE \bB^0\tilde{\bQ}  \left(\hat{\bB}' \hat{\bB} \right)^{-1}\right\|_{\F}\\
&\quad \quad  +\left\|\frac{1}{T} \tilde{\bH}_4^{\prime}\left({\bF^0}^{\prime} \bE \hat{\bB} \right)\left(\hat{\bB}' \hat{\bB}\right)^{-1}\right\|_{\F} \\
& \le O_p( \Delta_1)+O_p(\Delta_2),
\end{align*}
where $\Delta_1+\Delta_2 \lesssim \tilde{\Delta}_{NT}$ if $\frac{N^{1-\alpha_r}}{T}\to 0$.
\end{proof}
 
\begin{lem}
	\label{lemma bh-b0*e}
	Suppose that Assumptions \ref{cond:eigen}--\ref{assumption BaiNgA3} hold. If $\frac{N^{1-\alpha_r}}{T}\to 0$, then, we have
\begin{eqnarray*}
    && \left\| \bN^{\frac{1}{2}}(\hat{\bB}'\hat{\bB})^{-1} \left(  \hat{\bB} - \bB^0 \tilde{\bQ}'\right)' \be_t\right\|_{\F}   =
O_p(N^{\frac{1}{2}-\alpha_{r}} )+ O_p\left(\frac{N^{\frac{3}{2}-\alpha_{r}}}{T} \right)+ O_p\left(\sqrt{\frac{N^{1-\alpha_{r}}}{T}}\right).
\end{eqnarray*}	
\end{lem}
\begin{proof}
\eqref{B} implies
\begin{align*}
	(\hat{\bB}-\bB^0\tilde{\bQ}')'\be_t 
	&= \frac{1}{T} (\hat{\bF}-\bF^0 \tilde{\bH}_4)' \bE\be_t +{\tilde{\bH}
 _4}' \frac{1}{T} {\bF^0}'\bE \be_t \\
	&= \frac{1}{T} \hat{\bLambda}^{-1}  \left(\frac{1}{T} \bE'\bE\bE'\hat{\bF}   +\frac{1}{T} \bE'\bE\bB^0{\bF^0}'\hat{\bF} \right)'\be_t +{\tilde{\bH}_4}' \frac{1}{T} {\bF^0}'\bE \be_t . 
\end{align*}
Consider the first term on the right-hand side of the above equation:
\begin{eqnarray*}
	&&\left\| \bN^{\frac{1}{2}}(\hat{\bB}'\hat{\bB})^{-1}	\frac{1}{T^2} \hat{\bLambda}^{-1}{ \hat{\bF}}'\bE\bE'\bE\be_t\right\|_{\F}\\
	&&= \left\|	\bN^{\frac{1}{2}}(\hat{\bB}'\hat{\bB})^{-1}\bN^{\frac{1}{2}} \bN^{-\frac{1}{2}}	\frac{1}{T^2} \hat{\bLambda}^{-1}{ \hat{\bF}}'\bE\bE'\bE\be_t \right\|_{\F}\\
	&&\le 	\left\| \bN^{\frac{1}{2}}(\hat{\bB}'\hat{\bB})^{-1}\bN^{\frac{1}{2}} \bN^{-\frac{1}{2}}	\frac{1}{T^2} \hat{\bLambda}^{-1}\right\|_{\F} \left\| { \hat{\bF}}'\right\|_{\F} \left\| \bE\right\| ^3_{2} \left\| \be_t\right\|_{\F}  \\
	&&=   O_p\left(\frac{1}{T^2} N^{-\frac{3}{2}\alpha_{r}} \sqrt{T} (N^{3/2}+T^{3/2}) \sqrt{N} \right)\\
	&& = O_p\left(
 \frac{N^{\frac{3}{2}-\alpha_r}}{T} \sqrt{\frac{N^{1-\alpha_{r}}}{T}} \right)+O_p\left(
N^{\frac{1}{2}-\frac{3}{2}\alpha_{r}}\right).
\end{eqnarray*}

Next, consider the upper bound of the second term.
\begin{eqnarray*}
	&&\left\|\bN^{\frac{1}{2}}(\hat{\bB}'\hat{\bB})^{-1}	\frac{1}{T^2} \hat{\bLambda}^{-1}{ \hat{\bF}}'\bF^0{\bB^0}'\bE'\bE\be_t \right\|_{\F} \\
	&&\le  	\left\| \bN^{\frac{1}{2}}(\hat{\bB}'\hat{\bB})^{-1}\bN^{\frac{1}{2}} 	\frac{1}{T} \hat{\bLambda}^{-1}\right\|_{\F} 
	\left\| \bN^{-\frac{1}{2}}\frac{1}{T}{\hat{\bF}}'\bF^0 \bN^{\frac{1}{2}}\right\|_{\F} \left\|  \bN^{-\frac{1}{2}}\bB^0 \right\|_{\F} \lambda_{1}[  \bE'\bE] \left\| \be_t \right\|_{\F} 
	\\
	&&=  O_p\left( \frac{1}{T} N^{-\alpha_{r}}   (N+T)\sqrt{N} \right)\\
	&&=  O_p\left(
\frac{N^{\frac{3}{2}-\alpha_{r}}}{T} \right)+O_p\left(
N^{\frac{1}{2}- \alpha_{r}}\right) .
\end{eqnarray*}
Consider the third term and Lemma \ref{lem:BEF}(i) implies
\begin{eqnarray*}
	&&\left\|\bN^{\frac{1}{2}}(\hat{\bB}'\hat{\bB})^{-1}	 {\tilde{\bH}_4}' \frac{1}{T} {\bF^0}'\bE \be_t \right\|_{\F}\\
	&&\le  	\left\| \bN^{\frac{1}{2}}(\hat{\bB}'\hat{\bB})^{-1} {\tilde{\bH}_4}' \right\|_{\F} 
	\left\| 	  \frac{1}{T} {\bF^0}'\bE \be_t   \right\|_{\F} 
	\\
	&&= O_p\left(  N^{-\frac{1}{2} \alpha_{r}} N\left(\frac{1}{\sqrt{NT}}+\frac{1}{T}\right) \right) \\
	&&=   O_p\left(\sqrt{\frac{N^{1-\alpha_{r}}}{T}} \right)+ O_p\left(\frac{N^{1-\frac{1}{2}\alpha_{r}}}{T}\right),
\end{eqnarray*}
where $\frac{N^{1-\frac{1}{2}\alpha_{r}}}{T}=\frac{N^{\frac{3}{2}-\alpha_{r}}}{T}N^{\frac{1}{2}\alpha_r-\frac{1}{2}}=o(1)$. 
If $\frac{N^{1-\alpha_r}}{T}\to 0$ holds, collecting these terms completes the proof.
\end{proof}

\begin{lem}\label{lem:forecast1}
Suppose that Assumptions \ref{cond:eigen}--\ref{ass: augreg} hold. If $ \frac{N^{1-\alpha_r}}{T} \to 0$, then, we have
\begin{flalign*}
		&(i)	~\frac{1}{T}	 \| \hat{\bZ}'(\bF^0-\hat{\bF}) \|_{\F}
		=  O_p\left(\tilde{\Delta}_{NT} \right)  , &\\
		&(ii)~\frac{1}{\sqrt{T}}	 \| (\hat{\bF}-\bF^0)'\bepsilon \|_{\F}
		=   O_p\left(\sqrt{T} \tilde{\Delta}_{NT}\right)  .&
	\end{flalign*}
\end{lem}
\begin{proof}(i) Because $\hat{\bZ}$ includes $\hat{\bF}$ and $\bW$, we have
\begin{eqnarray*}
	\frac{1}{T}	 \left\| \hat{\bF}'(\bF^0-\hat{\bF}) \right\|_{\F}
		&=&
 	 \left\| \tilde{\bQ}-\bI_r \right\|_{\F} =
  O_p\left(\tilde{\Delta}_{NT}\right) ,
\end{eqnarray*}
\begin{eqnarray*}
	&&\frac{1}{T}	 \left\| \bW'(\bF^0-\hat{\bF}) \right\|_{\F}\\
	&&	=
  \frac{1}{T}	 \left\| \bW'(\bF^0-\bF^0\tilde{\bH}_4+\bF^0\tilde{\bH}_4-\hat{\bF}) \right\|_{\F}  \\
 & &\le \frac{1}{T}	 \left\| \bW'\bF^0(\bI_r- \tilde{\bH}_4)  \right\|_{\F} + \frac{1}{T} \left\|\bW'(\bF^0\tilde{\bH}_4-\hat{\bF}) \right\|_{\F} \\
 && = 
O_p\left( \tilde{\Delta}_{NT} \right)
+ O_p\left(\frac{N^{1- \alpha_{r}}}{T}\right)
+O_p\left(\frac{1}{N^{\alpha_r}}  \right)
+O_p\left(   \frac{1}{\sqrt{TN^{\alpha_r}}} \right).
\end{eqnarray*}
The final inequality is from
\begin{eqnarray*}
	&& \frac{1}{T}\left\| \bW'(\bF^0\tilde{\bH}_4-\hat{\bF}) \right\|_{\F} \\
 &&= \frac{1}{T}\left\|\bW'\left( \frac{1}{T} \bE\bE'\hat{\bF}+\frac{1}{T} \bE\bB^0{\bF^0}'\hat{\bF} \right) \hat{\bLambda}^{-1} \right\|_{\F} \\
 && \le 
 \frac{1}{T^2}\left\|\bW'  \bE\bE'\hat{\bF} \hat{\bLambda}^{-1} \right\|_{\F}  + \frac{1}{T^2}\left\| (\bW'\bE\bB^0){\bF^0}'\hat{\bF} \hat{\bLambda}^{-1} \right\|_{\F} \\
 && =O_p
 \left(  \frac{N^{1- \alpha_{r}}}{T}+\frac{1}{N^{\alpha_r}} \right)
 +O_p\left(
   \frac{1}{\sqrt{TN^{\alpha_r}}}\right) \\
 && =
  O_p\left(\frac{N^{1- \alpha_{r}}}{T}\right)
  +O_p\left(\frac{1}{N^{\alpha_r}}  \right)
  +O_p\left(   \frac{1}{\sqrt{TN^{\alpha_r}}}\right),
\end{eqnarray*}
which is less than $O_p(\tilde{\Delta}_{NT})$.

(ii) Using Lemmas \ref{lem:H equality}(ii), \ref{lem:FEB0}(i), and Proposition \ref{prop:Q-I dec}
\begin{eqnarray*}
    && \left\|\tilde{\bH}_4-\bI_r\right\|_{\F}\\
    &&=\left\|\tilde{\bH}_4-\tilde{\bH}_2+\tilde{\bH}_2-\bI_r\right\|_{\F}\\
    && = \left\|- \frac{1}{T}{\bF}^{0'}\bE \hat{\bB} \bLambda^{-1}+\tilde{\bH}_2-\bI_r\right\|_{\F} \\
    &&\le O_p(\Delta_1)+O_p(\tilde{\Delta}_{NT}) = O_p( \tilde{\Delta}_{NT}),
\end{eqnarray*}
if $\frac{N^{1-\alpha_r}}{T} \to 0$. Then,
\begin{eqnarray*}
	&&\frac{1}{\sqrt{T}}	 \left\| (\hat{\bF}-\bF^0)'\bepsilon \right\|_{\F} 
 \\
 && \le
 \frac{1}{\sqrt{T}}	 \left\| (\hat{\bF}-\bF^0\tilde{\bH}_4)'\bepsilon \right\|_{\F} + \frac{1}{\sqrt{T}}	\left\| (\tilde{\bH}_4-\bI_r)'{\bF^0}'\bepsilon \right\|_{\F} \\
 && \le  \frac{1}{T^{3/2}}\left\|\hat{\bLambda}^{-1} {\hat{\bF}}' \bE\bE' \bepsilon\right\|_{\F} + \frac{1}{T^{3/2}}\left\| \hat{\bLambda}^{-1}{\hat{\bF}}'\bF^0 {\bB^0}'{\bE}'\bepsilon \right\|_{\F}  + \frac{1}{\sqrt{T}}	\left\| (\tilde{\bH}_4-\bI_r)'{\bF^0}'\bepsilon \right\|_{\F}\\
 && \le 
 O_p\left( \frac{N^{1- \alpha_{r}}}{\sqrt{T}}+\frac{\sqrt{T}}{N^{\alpha_r}} \right)
 +O_p\left(  \frac{\sqrt{T}}{N^{\alpha_r}} +\frac{1}{\sqrt{N^{\alpha_r}}}\right)   
   + O_p(\tilde{\Delta}_{NT}) \\
   && =  O_p\left(\frac{N^{1- \alpha_{r}}}{\sqrt{T}}\right)
   +O_p\left(\frac{\sqrt{T}}{N^{\alpha_r}} \right) ,
\end{eqnarray*}
since $\frac{1}{T}\bN^{-\frac{1}{2}}{\bB^0}'{\bE}'\bepsilon =O_p\left(\frac{1}{\sqrt{T}}+\frac{1}{\sqrt{N^{\alpha_r}}}\right)$, $\frac{N^{1-\alpha_r}}{T} \to 0$ and $2\sqrt{ab} \le a+b$.
\end{proof}

\begin{lem}
\label{lem:delta34}
Suppose that Assumptions \ref{cond:eigen}--\ref{ass: augreg} hold. 
If $\frac{N^{1-\alpha_r}}{\sqrt{T}} \to  0$, and $\sqrt{T}N^{-\alpha_r} \to 0$ as $N, T \to \infty$, we have 
\begin{eqnarray*}
		\nonumber
		&&\sqrt{T} ( \hat{\bdelta}-\bdelta_3) \CD N(\mathbf{0}, \bSigma_{\delta^0}) , \\
  &&\sqrt{T} ( \hat{\bdelta}-\bdelta_4) \CD N(\mathbf{0}, \bSigma_{\delta^0}) ,
\end{eqnarray*}
where $\bgamma_j=\tilde{\bH}_j^{-1} \bgamma^0 , \; \bdelta_j =({\bgamma}_j ',\bbeta')', \; j=3,4$, $\bSigma_{\delta^0}=\bSigma_{z^0}^{-1}\bSigma_{z^0\epsilon}\bSigma_{z^0}^{-1}$.
\end{lem}

\begin{proof}
The augmented regression model is rewritten as 
 \begin{eqnarray*}
          y_{t+h}&=&{\bgamma^*}'\bff_t^*+\bbeta' \bw_t+\epsilon_{t+h}\\
          &=&{\bgamma^0}'\tilde{{\bH}}_j^{'-1}\hat{\bff}_t+\bbeta' \bw_t+\epsilon_{t+h}+{\bgamma^0}'\tilde{{\bH}}_j^{'-1}\tilde{{\bH}}_j^{\prime}\bff_t^0-{\bgamma^0}'\tilde{{\bH}}_j^{'-1}\hat{\bff}_t. 
      \end{eqnarray*}

(i) We start with $j=3$.
\begin{eqnarray*}
   && \frac{1}{\sqrt{T}} \hat{\bF}'(\bF^0\tilde{\bH}_3 -\hat{\bF})\tilde{\bH}_3^{-1}\bgamma^0=0 \\
   && \frac{1}{\sqrt{T}} \bW'(\bF^0\tilde{\bH}_3 -\hat{\bF})\tilde{\bH}_3^{-1}\bgamma^0 \\
   &&= \frac{1}{\sqrt{T}} (\bW'\bF^0  -\bW'\hat{\bF}\tilde{\bH}_3^{-1})\bgamma^0
   \\
   &&= \frac{1}{\sqrt{T}} (\bW'\bF^0  -\bW'\bX\hat{\bB}\hat{\bLambda}^{-1}\tilde{\bQ})\bgamma^0\\
   &&   = \frac{1}{\sqrt{T}} \bW'\bF^0(  \bI_r- {\bB^0}'\hat{\bB}\hat{\bLambda}^{-1}\tilde{\bQ})-\frac{1}{\sqrt{T}}( \bW'\bE'\hat{\bB}\hat{\bLambda}^{-1}\tilde{\bQ})\bgamma^0 \\
   && = \sqrt{T}O_p(\tilde{\Delta}_{NT}).
\end{eqnarray*}
By the proof of Theorem \ref{thm:forecast1}, we have
\begin{eqnarray*} 
	&&	\sqrt{T} ( \hat{\bdelta}-\bdelta_3)
  \\
  &&=\left( \frac{1}{T} {\hat{\bZ}'\hat{\bZ}}\right)^{-1}\frac{1}{\sqrt{T}} \hat{\bZ}'\bepsilon + \left( \frac{1}{T} {\hat{\bZ}'\hat{\bZ}}\right)^{-1}\frac{1}{\sqrt{T}} \hat{\bZ}'(\bF^0\tilde{\bH}_3 -\hat{\bF})\tilde{\bH}_3^{-1}\bgamma^0  \\
  && =\left( \frac{1}{T} { {\bZ^0}' {\bZ^0}}\right)^{-1}\frac{1}{\sqrt{T}}  {\bZ^0}'\bepsilon+\left( \frac{1}{T} { {\bZ^0}' {\bZ^0}}\right)^{-1}\frac{1}{\sqrt{T}}  (\hat{\bZ}-{\bZ^0})'\bepsilon + \left[\left( \frac{1}{T} {\hat{\bZ}'\hat{\bZ}}\right)^{-1}-   \left( \frac{1}{T} { {\bZ^0}' {\bZ^0}}\right)^{-1}  \right]\frac{1}{\sqrt{T}} \hat{\bZ}'\bepsilon \\
  && \quad \quad +\left( \frac{1}{T}{\hat{\bZ}'\hat{\bZ}}\right)^{-1}\frac{1}{\sqrt{T}} \hat{\bZ}'(\bF^0\tilde{\bH}_3 -\hat{\bF})\tilde{\bH}_3^{-1}\bgamma^0 \\
 && = \left( \frac{1}{T} {{\bZ^0}'{\bZ^0}}\right)^{-1}\frac{1}{\sqrt{T}} {\bZ^0}'\bepsilon +   O_p\left(\frac{N^{1 - \alpha_{r}}}{\sqrt{T}} \right) + O_p\left(\sqrt{T} N^{-\alpha_r}\right) .
\end{eqnarray*}
Thus,
$
\sqrt{T} ( \hat{\bdelta}-\bdelta_3) \CD  N(\mathbf{0}, \bSigma_{\delta^0})
$
if $\frac{N^{1-\alpha_r}}{\sqrt{T}} \to  0$, and $\sqrt{T}N^{-\alpha_r} \to 0$.

(ii) Let $j=4$, 
\begin{eqnarray*}
   && \frac{1}{\sqrt{T}} \hat{\bF}'(\bF^0\tilde{\bH}_4 -\hat{\bF})\tilde{\bH}_4^{-1}\bgamma^0\\
   &&= \sqrt{T} \left(\frac{1}{T}\hat{\bF}'\bF^0  -\tilde{\bH}_4^{-1}\right) \bgamma^0\\
   &&= \sqrt{T} \left[\frac{1}{T}\hat{\bF}'\bE\hat{\bB}(\hat{\bB}'\bB^0)^{-1}  \right] \bgamma^0 \\
   &&= \sqrt{T}O_p(\Delta_{1}+\Delta_2),\\
   && \frac{1}{\sqrt{T}} \bW'(\bF^0\tilde{\bH}_4 -\hat{\bF})\tilde{\bH}_4^{-1}\bgamma^0 \\
   &&=  \sqrt{T} \left(-\frac{1}{T^2} \bW'\bE\bB^0{\bF^0}'\hat{\bF}\hat{\bLambda}^{-1}  -\frac{1}{T^2}\bW'\bE\bE'\hat{\bF}\hat{\bLambda}^{-1}\right)\tilde{\bH}_4^{-1}\bgamma^0
   \\
   && = O_p\left(N^{-\frac{1}{2} \alpha_{r} }\right)+O_p\left(\frac{N^{1 - \alpha_{r}}}{\sqrt{T}}  + \sqrt{T} N^{-\alpha_r}\right) .
\end{eqnarray*}
By the proof of Theorem \ref{thm:forecast1}, we have
\begin{eqnarray*} 
	&&	\sqrt{T} ( \hat{\bdelta}-\bdelta_4)
  \\&& =\left( \frac{1}{T} {\hat{\bZ}'\hat{\bZ}}\right)^{-1}\frac{1}{\sqrt{T}} \hat{\bZ}'\bepsilon + \left( \frac{1}{T} {\hat{\bZ}'\hat{\bZ}}\right)^{-1}\frac{1}{\sqrt{T}} \hat{\bZ}'(\bF^0\tilde{\bH}_4 -\hat{\bF})\tilde{\bH}_4^{-1}\bgamma^0  \\
  && =\left( \frac{1}{T} { {\bZ^0}' {\bZ^0}}\right)^{-1}\frac{1}{\sqrt{T}}  {\bZ^0}'\bepsilon+\left( \frac{1}{T} { {\bZ^0}' {\bZ^0}}\right)^{-1}\frac{1}{\sqrt{T}}  (\hat{\bZ}-{\bZ^0})'\bepsilon + \left[\left( \frac{1}{T} {\hat{\bZ}'\hat{\bZ}}\right)^{-1}-   \left( \frac{1}{T} { {\bZ^0}' {\bZ^0}}\right)^{-1}  \right]\frac{1}{\sqrt{T}} \hat{\bZ}'\bepsilon \\
  && \quad \quad +\left( \frac{1}{T}{\hat{\bZ}'\hat{\bZ}}\right)^{-1}\frac{1}{\sqrt{T}} \hat{\bZ}'(\bF^0\tilde{\bH}_4 -\hat{\bF})\tilde{\bH}_4^{-1}\bgamma^0 \\
 && = \left( \frac{1}{T} {{\bZ^0}'{\bZ^0}}\right)^{-1}\frac{1}{\sqrt{T}} {\bZ^0}'\bepsilon +   O_p\left(\frac{N^{1 - \alpha_{r}}}{\sqrt{T}} \right) + O_p\left(\sqrt{T} N^{-\alpha_r}\right).
\end{eqnarray*}
Thus,
$
\sqrt{T} ( \hat{\bdelta}-\bdelta_4) \CD  N(\mathbf{0}, \bSigma_{\delta^0})
$
if $\frac{N^{1-\alpha_r}}{\sqrt{T}} \to  0$, and $\sqrt{T}N^{-\alpha_r} \to 0$.

\end{proof}

\begin{lem}\label{lem:yh-y34}
Under the assumptions of Theorem \ref{thm:forecast2}, if we use rotation matrices $\tilde{\bH}_3 $ or $\tilde{\bH}_4 $, we obtain the same results as in Theorem \ref{thm:forecast2}.
\[
\frac{\left(\hat{y}_{T+h \mid T}-y_{T+h \mid T}\right)}{{\sigma}_{T+h \mid T}} \stackrel{d}{\longrightarrow} N(0,1),
\]
where ${\sigma}_{T+h \mid T}^2=
T^{-1} {\bz}_T^{0\prime}\bSigma_{\delta^0}{\bz}_T^0
+{\bgamma}^{0 \prime} 
\bD^{-1}\bN^{-1/2}\bGamma_T\bN^{-1/2}\bD^{-1}
{\bgamma}^0$.
\end{lem}

\begin{proof}
Expand the term
\begin{eqnarray*}
&& \hat{y}_{T+h \mid T}-y_{T+h \mid T} \\
&& =\hat{\bgamma}^{\prime} \hat{\bff}_T+\hat{\bbeta}^{\prime} \bw_T-\bgamma^{*\prime} \bff_T^*-\bbeta^{\prime} \bw_T \\
&& =\left(\hat{\bgamma}-\hat{\bH}_j^{-1} \bgamma^*\right)^{\prime} \hat{\bff}_T+\bgamma^{*\prime} \hat{\bH}_j^{-1\prime}\left(\hat{\bff}_T-\hat{\bH}_j^{\prime} \bff_T^*\right)+(\hat{\bbeta}-\bbeta)^{\prime} \bw_T \\
&& =\hat{\bz}_T^{\prime}(\hat{\bdelta}-\bdelta_j)+\bgamma^{*\prime} \hat{\bH}_j^{-1}\left(\hat{\bff}_T-\hat{\bH}_j^{\prime} \bff_T\right) \\
&& =T^{-1 / 2} \hat{\bz}_T^{\prime}[\sqrt{T}(\hat{\bdelta}-\bdelta_j)]+ \bgamma^{*\prime} \hat{\bH}_j^{-1\prime} \bN^{-1 / 2}\left[\bN^{1 / 2}\left(\hat{\bff}_T-\tilde{\bH}_j^{\prime} \bff^0_T\right)\right] \\
&& =T^{-1 / 2} \bz_T^{0\prime}[\sqrt{T}(\hat{\bdelta}-\bdelta_j)]+T^{-1 / 2} (\hat{\bz}_T-\bz^0_T)^{\prime}[\sqrt{T}(\hat{\bdelta}-\bdelta_j)]\\
&& \quad \quad + \bgamma^{*\prime} \hat{\bH}_j^{-1\prime} \bN^{-1 / 2}\left[\bN^{1 / 2}\left(\hat{\bff}_T-\tilde{\bH}_4^{\prime} \bff^0_T\right)\right]+ \bgamma^{*\prime} \hat{\bH}_j^{-1\prime} \bN^{-1 / 2}\left[\bN^{1 / 2}\left(\tilde{\bH}_4^{\prime} -\tilde{\bH}_j^{\prime} \right) \bff^0_T\right] .
\end{eqnarray*}
(i) We start with $j=3$. The second term on the righat-hand side of the above equation is dominated by the first one, thus, we next focus on the first, third and fourth terms.
\begin{eqnarray*}
    && T^{-1 / 2} \bz_T^{0\prime}\sqrt{T}(\hat{\bdelta}-\bdelta_3) \\
    && =T^{-1 / 2} {\bz}_T^{0\prime} \left( \frac{1}{T} {{\bZ^0}'{\bZ^0}}\right)^{-1}\frac{1}{\sqrt{T}} {\bZ^0}'\bepsilon  \\
    && \quad \quad + T^{-1 / 2}  \left[ O_p\left(\frac{N^{1 - \alpha_{r}}}{\sqrt{T}} \right) + O_p\left(\sqrt{T} N^{-\alpha_r}\right)\right] .
    \end{eqnarray*}	
    \begin{eqnarray*}
    && \bgamma^{*\prime}  \hat{\bH}_3^{-1\prime} \bN^{-\frac{1}{2}} \left[\bN^{1 / 2}\left(\hat{\bff}_T-\tilde{\bH}_4^{\prime} \bff_T^0\right)\right] \\
      && = \bgamma^{*\prime}  \hat{\bH}_3^{-1\prime} \bN^{-\frac{1}{2}} \bD^{-1}  
  \left[\bD\bN^{\frac{1}{2}}(\hat{\bB}'\hat{\bB})^{-1}   \bB^{0'}\be_t\right]
	+  O(\| \bN^{-\frac{1}{2}}\|_F) O_p(\tilde{\Delta}_{NT})  \\
 && \quad \quad + O(\| \bN^{-\frac{1}{2}}\|_{\bF})  \left[ O_p(N^{\frac{1}{2}-\alpha_{r}} )+ O_p\left(\frac{N^{\frac{3}{2}-\alpha_{r}}}{T} \right)+ O_p\left(\sqrt{\frac{N^{1-\alpha_{r}}}{T}}\right)\right].\\
    &&   \bgamma^{*\prime} \hat{\bH}_3^{-1\prime} \left(\tilde{\bH}_3^{\prime} - \tilde{\bH}_4^{\prime} \right) \bff_T^0 
      =
      \frac{1}{\sqrt{T}} O_p( \sqrt{T} \tilde{\Delta}_{NT}), \\
\end{eqnarray*}	
where
\begin{eqnarray*}
    && \sqrt{T}\tilde{\Delta}_{NT}
	=
	   \frac{N^{1 -\alpha_r}}{\sqrt{T}} +\sqrt{T} N^{ -\alpha_{r}} .
\end{eqnarray*}
Thus, if $\sqrt{T} N^{-\alpha_r}\to 0$, $\frac{1}{2}< \alpha_r$, and $\frac{N^{\frac{3}{2}-\alpha_{r}}}{T} \to 0$, we have
\[
\frac{\left(\hat{y}_{T+h \mid T}-y_{T+h \mid T}\right)}{{\sigma}_{T+h \mid T}} \stackrel{d}{\longrightarrow} N(0,1),
\]
where ${\sigma}_{T+h \mid T}^2=
T^{-1} {\bz}_T^{0\prime}\bSigma_{\delta^0}{\bz}_T^0
+{\bgamma}^{0 \prime} 
\bD^{-1}\bN^{-1/2}\bGamma_T\bN^{-1/2}\bD^{-1}
{\bgamma}^0$.

(ii)
Let $j=4$, 
\begin{eqnarray*}
    && T^{-1 / 2} \bz_T^{0\prime}\sqrt{T}(\hat{\bdelta}-\bdelta_4) \\
   && =T^{-1 / 2} {\bz}_T^{0\prime} \left( \frac{1}{T} {{\bZ^0}'{\bZ^0}}\right)^{-1}\frac{1}{\sqrt{T}} {\bZ^0}'\bepsilon  \\
    && \quad \quad + T^{-1 / 2}  \left[ O_p\left(\frac{N^{1 - \alpha_{r}}}{\sqrt{T}} \right) + O_p\left(\sqrt{T} N^{-\alpha_r}\right)\right] .\\
    && \bgamma^{*\prime}  \hat{\bH}_4^{-1\prime} \bN^{-\frac{1}{2}} \left[\bN^{1 / 2}\left(\hat{\bff}_T-\tilde{\bH}_4^{\prime} \bff_T^0\right)\right] \\
      && = \bgamma^{*\prime}  \hat{\bH}_4^{-1\prime} \bN^{-\frac{1}{2}} \bD^{-1}  
  \left[\bD\bN^{\frac{1}{2}}(\hat{\bB}'\hat{\bB})^{-1}   \bB^{0'}\be_t\right]
	+  O_p(\| \bN^{-\frac{1}{2}}\|_F) O_p(\tilde{\Delta}_{NT})  \\
 && \quad \quad + O_p(\| \bN^{-\frac{1}{2}}\|_F)  \left[ O_p(N^{\frac{1}{2}-\alpha_{r}} )+ O_p\left(\frac{N^{\frac{3}{2}-\alpha_{r}}}{T} \right)+ O_p\left(\sqrt{\frac{N^{1-\alpha_{r}}}{T}}\right)\right].
\end{eqnarray*}	
Thus, if $\sqrt{T} N^{-\alpha_r}\to 0$, $\frac{1}{2}< \alpha_r$, and $\frac{N^{\frac{3}{2}-\alpha_{r}}}{T} \to 0$, we have
\[
\frac{\left(\hat{y}_{T+h \mid T}-y_{T+h \mid T}\right)}{{\sigma}_{T+h \mid T}} \stackrel{d}{\longrightarrow} N(0,1),
\]
where ${\sigma}_{T+h \mid T}^2=
T^{-1} {\bz}_T^{0\prime}\bSigma_{\delta^0}{\bz}_T^0
+{\bgamma}^{0 \prime} 
\bD^{-1}\bN^{-1/2}\bGamma_T\bN^{-1/2}\bD^{-1}
{\bgamma}^0$.
\end{proof}

\section{Additional Experimental Results: Sparse loadings}
\label{subsec:addMC}

\setcounter{table}{0}
\renewcommand{\thetable}{C.\arabic{table}}
In this section, we report experimental results in the case of sparsity-induced weak factor models. The results are summarized in Tables \ref{table:sparseNT50}-\ref{table:sparseNT500}. The data generating process is identical to the one used for the results reported in Section \ref{sec: MC}, except that the loadings are sparse; see design 2 in Section \ref{mc: fac des} for more details.
The results are qualitatively the same as those for the non-sparse
designs in Section~\ref{sec: MC}.


\renewcommand{\arraystretch}{.9}
\begin{table}[htb!]
\caption{Experimental results with sparse factor loadings for $N=T=50$}
\centering\small
\begin{tabular}{lrrrrrr}
\hline
$(\alpha_{1},\alpha_{2})$ & (1.0,1.0) & (1.0,0.9) & (1.0,0.8) & (0.9,0.7) & (0.8,0.6) & (0.7,0.5) \\ \hline
\multicolumn{7}{c}{Factor models} \\ \hline
\multicolumn{7}{l}{Estimation Accuracy (x1000)} \\
\ \ $T^{-\frac{1}{2}}||\mathbf{\hat{F}}-\mathbf{F}^{0}||_{\text{F}}$ & 110.3 & 128.8 & 149.9 & 183.6 & 212.9 & 261.9 \\
\ \ $T^{-\frac{1}{2}}||\mathbf{\hat{F}}-\mathbf{F}^{\ast}\mathbf{\hat{H}}_{4}||_{\text{F}}$ & 111.7 & 131.2 & 153.2 & 188.1 & 218.3 & 269.3 \\
\ \ $T^{-\frac{1}{2}}||\mathbf{\hat{F}}-\mathbf{F}^{\ast}\mathbf{\hat{H}}||_{\text{F}}$ & 119.0 & 142.1 & 169.4 & 208.8 & 242.5 & 300.3 \\
\ \ $T^{-\frac{1}{2}}||\mathbf{\hat{F}}-\mathbf{F}^{\ast}\mathbf{\hat{R}}_{\mathrm{Fan},F}||_{\text{F}}$ & 109.3 & 128.2 & 149.4 & 183.1 & 212.1 & 261.0 \\
\multicolumn{7}{l}{} \\[-7pt]
\ \ $N^{-\frac{1}{2}}||\mathbf{\hat{B}}-\mathbf{B}^{0}||_{\text{F}}$ & 142.5 & 142.3 & 142.6 & 143.3 & 144.3 & 145.9 \\
\ \ $N^{-\frac{1}{2}}||\mathbf{\hat{B}}-\mathbf{B}^{\ast}\mathbf{\hat{Q}}^{\prime}||_{\text{F}}$ & 142.3 & 142.7 & 143.4 & 144.4 & 145.6 & 147.7 \\
\ \ $N^{-\frac{1}{2}}||\mathbf{\hat{B}}-\mathbf{B}^{\ast}\mathbf{\hat{H}}^{\prime-1}||_{\text{F}}$ & 145.1 & 147.3 & 150.4 & 152.8 & 155.2 & 160.3 \\
\ \ $N^{-\frac{1}{2}}||\mathbf{\hat{B}}-\mathbf{B}^{\ast}\mathbf{\hat{R}}_{\mathrm{Fan},B}||_{\text{F}}$ & 142.0 & 142.3 & 142.9 & 143.5 & 144.5 & 146.1 \\
\multicolumn{7}{l}{} \\[-7pt]
\ \ $(NT)^{-\frac{1}{2}}||\mathbf{\hat{C}}-\mathbf{C}^{\ast}||_{\text{F}}$ & 198.6 & 198.8 & 199.4 & 200.1 & 201.0 & 202.6 \\

\multicolumn{7}{l}{Normal Approximation (rejection frequencies of 5\% $\chi_{r}^{2}$ tests)} \\
\ \ $\mathbf{\hat{f}}_{t}-\mathbf{f}_{t}^{0}$ & 5.0 & 5.0 & 5.0 & 5.1 & 5.1 & 5.0 \\
\ \ $\mathbf{\hat{f}}_{t}-\mathbf{\hat{H}}_{4}^{\prime}\mathbf{f}_{t}^{\ast}$ & 4.9 & 5.0 & 5.1 & 5.1 & 5.2 & 5.2 \\
\ \ $\mathbf{\hat{f}}_{t}-\mathbf{\hat{H}}^{\prime}\mathbf{f}_{t}^{\ast}$ & 5.4 & 5.4 & 5.6 & 5.7 & 5.8 & 6.1 \\
\ \ $\mathbf{\hat{f}}_{t}-\mathbf{\hat{R}}_{\mathrm{Fan},F}^{\prime}\mathbf{f}_{t}^{\ast}$ & 4.7 & 4.8 & 4.8 & 5.0 & 5.0 & 4.9 \\
\multicolumn{7}{l}{} \\[-7pt]
\ \ $\mathbf{\hat{b}}_{i}-\mathbf{b}_{i}^{0}$ & 5.5 & 5.3 & 5.4 & 5.5 & 6.0 & 6.6 \\
\ \ $\mathbf{\hat{b}}_{i}-\mathbf{\hat{Q}b}_{i}^{\ast}$ & 5.4 & 5.3 & 5.7 & 6.1 & 7.0 & 8.7 \\
\ \ $\mathbf{\hat{b}}_{i}-\mathbf{\hat{H}}^{-1}\mathbf{b}_{i}^{\ast}$ & 5.9 & 7.3 & 9.7 & 13.1 & 18.1 & 30.9 \\
\ \ $\mathbf{\hat{b}}_{i}-\mathbf{\hat{R}}_{\mathrm{Fan},B}'\mathbf{b}_{i}^{\ast}$ & 5.4 & 5.2 & 5.4 & 5.5 & 6.0 & 6.5 \\
\multicolumn{7}{l}{Test for linear restrictions (rejection frequencies of 5\% $N(0,1)$ tests)} \\
\ \ $H_{0}:f_{t,1}=f_{t,1}^{0}$ & 5.3 & 5.1 & 4.9 & 5.1 & 5.1 & 5.0 \\
\ \ $H_{0}:f_{t,2}=f_{t,2}^{0}$ & 4.8 & 4.9 & 4.7 & 4.9 & 4.8 & 4.9 \\
\multicolumn{7}{l}{} \\[-7pt]
\ \ $H_{0}:b_{i,1}=b_{i,1}^{0}$ & 5.3 & 5.4 & 5.3 & 5.4 & 5.6 & 6.0 \\
\ \ $H_{0}:b_{i,2}=b_{i,2}^{0}$ & 5.1 & 5.2 & 5.4 & 5.5 & 5.8 & 6.4 \\
\multicolumn{7}{l}{} \\[-7pt]
\ \ $H_{0}:c_{t,i}=c_{t,i}^{\ast}$ & 4.8 & 5.1 & 5.0 & 5.1 & 5.3 & 5.5 \\
\hline
\multicolumn{7}{c}{Factor Augmented Regressions} \\ \hline
\multicolumn{7}{l}{Estimation Accuracy (x1000)} \\
\ \ $||\boldsymbol{\hat{\gamma}}^{0}-\boldsymbol{\gamma}^{0}||_{\text{F}}$ & 227.1 & 227.7 & 228.4 & 228.4 & 227.7 & 226.8 \\
\ \ $||\boldsymbol{\hat{\gamma}}-\boldsymbol{\gamma}^{0}||_{\text{F}}$ & 227.7 & 228.8 & 229.8 & 231.5 & 233.8 & 238.6 \\
\ \ $||\boldsymbol{\hat{\gamma}}-\mathbf{\hat{H}}^{-1}\boldsymbol{\gamma}^{\ast}||_{\text{F}}$ & 233.1 & 239.2 & 248.5 & 265.6 & 286.1 & 338.0 \\
\multicolumn{7}{l}{Normal Approximation (rejection frequencies of 5\% $\chi_{r}^{2}$ tests)} \\
\ \ $\boldsymbol{\hat{\gamma}}^{0}-\boldsymbol{\gamma}^{0}$ & 4.9 & 5.0 & 4.9 & 5.0 & 5.0 & 4.9 \\
\ \ $\boldsymbol{\hat{\gamma}}-\boldsymbol{\gamma}^{0}$ & 4.9 & 5.0 & 5.1 & 5.5 & 5.8 & 6.2 \\
\ \ $\boldsymbol{\hat{\gamma}}-\mathbf{\hat{H}}^{-1}\boldsymbol{\gamma}^{\ast}$ & 5.7 & 6.6 & 7.9 & 10.7 & 14.3 & 24.1 \\
\multicolumn{7}{l}{Test for linear restrictions (rejection frequencies of 5\% $N(0,1)$ tests)} \\
\ \ $\hat{\gamma}_{1}^{0}-\gamma_{1}^{0}$ & 4.8 & 5.0 & 5.0 & 5.1 & 5.0 & 4.9 \\
\ \ $\hat{\gamma}_{1}-\gamma_{1}^{0}$ & 4.9 & 5.1 & 5.1 & 5.4 & 5.5 & 5.8 \\
\multicolumn{7}{l}{} \\[-7pt]
\ \ $\hat{\gamma}_{2}^{0}-\gamma_{2}^{0}$ & 4.9 & 4.8 & 5.1 & 5.1 & 5.0 & 4.9 \\
\ \ $\hat{\gamma}_{2}-\gamma_{2}^{0}$ & 5.0 & 5.1 & 5.2 & 5.5 & 5.8 & 6.5 \\
\multicolumn{7}{l}{Empirical coverage} \\
\ \ \ 95\% $CI(y_{T+h|T})$ & 94.9 & 94.4 & 94.0 & 93.1 & 92.4 & 90.8 \\
\ \ \ 95\% $CI(y_{T+h})$ & 94.9 & 94.9 & 95.1 & 95.1 & 95.0 & 94.7 \\
\multicolumn{7}{l}{} \\[-7pt]
\ \ \ 90\% $CI(y_{T+h|T})$ & 89.9 & 89.4 & 88.8 & 87.7 & 86.9 & 84.8 \\
\ \ \ 90\% $CI(y_{T+h})$ & 90.0 & 89.9 & 90.1 & 90.1 & 90.0 & 89.7 \\
\hline
\end{tabular}
\label{table:sparseNT50}
\end{table}

\clearpage

\renewcommand{\arraystretch}{.9}
\begin{table}[htb!]
\caption{Experimental results with sparse factor loadings for $N=T=100$}
\centering\small
\begin{tabular}{lrrrrrr}
\hline
$(\alpha_{1},\alpha_{2})$ & (1.0,1.0) & (1.0,0.9) & (1.0,0.8) & (0.9,0.7) & (0.8,0.6) & (0.7,0.5) \\ \hline
\multicolumn{7}{c}{Factor models} \\ \hline
\multicolumn{7}{l}{Estimation Accuracy (x1000)} \\
\ \ $T^{-\frac{1}{2}}||\mathbf{\hat{F}}-\mathbf{F}^{0}||_{\text{F}}$ & 78.5 & 94.4 & 116.2 & 144.3 & 184.0 & 233.1 \\
\ \ $T^{-\frac{1}{2}}||\mathbf{\hat{F}}-\mathbf{F}^{\ast}\mathbf{\hat{H}}_{4}||_{\text{F}}$ & 79.1 & 95.3 & 117.6 & 146.3 & 187.2 & 237.9 \\
\ \ $T^{-\frac{1}{2}}||\mathbf{\hat{F}}-\mathbf{F}^{\ast}\mathbf{\hat{H}}||_{\text{F}}$ & 81.8 & 99.8 & 125.5 & 156.8 & 202.0 & 259.6 \\
\ \ $T^{-\frac{1}{2}}||\mathbf{\hat{F}}-\mathbf{F}^{\ast}\mathbf{\hat{R}}_{\mathrm{Fan},F}||_{\text{F}}$ & 78.2 & 94.2 & 116.0 & 144.1 & 183.8 & 232.8 \\
\multicolumn{7}{l}{} \\[-7pt]
\ \ $N^{-\frac{1}{2}}||\mathbf{\hat{B}}-\mathbf{B}^{0}||_{\text{F}}$ & 100.4 & 100.3 & 100.5 & 100.9 & 101.6 & 102.7 \\
\ \ $N^{-\frac{1}{2}}||\mathbf{\hat{B}}-\mathbf{B}^{\ast}\mathbf{\hat{Q}}^{\prime}||_{\text{F}}$ & 100.3 & 100.5 & 100.9 & 101.4 & 102.4 & 103.8 \\
\ \ $N^{-\frac{1}{2}}||\mathbf{\hat{B}}-\mathbf{B}^{\ast}\mathbf{\hat{H}}^{\prime-1}||_{\text{F}}$ & 101.3 & 102.3 & 103.9 & 105.0 & 107.1 & 110.5 \\
\ \ $N^{-\frac{1}{2}}||\mathbf{\hat{B}}-\mathbf{B}^{\ast}\mathbf{\hat{R}}_{\mathrm{Fan},B}||_{\text{F}}$ & 100.2 & 100.3 & 100.6 & 101.0 & 101.8 & 102.8 \\
\multicolumn{7}{l}{} \\[-7pt]
\ \ $(NT)^{-\frac{1}{2}}||\mathbf{\hat{C}}-\mathbf{C}^{\ast}||_{\text{F}}$ & 140.9 & 141.1 & 141.3 & 141.7 & 142.4 & 143.5 \\

\multicolumn{7}{l}{Normal Approximation (rejection frequencies of 5\% $\chi_{r}^{2}$ tests)} \\
\ \ $\mathbf{\hat{f}}_{t}-\mathbf{f}_{t}^{0}$ & 4.9 & 5.0 & 4.9 & 4.9 & 5.0 & 5.0 \\
\ \ $\mathbf{\hat{f}}_{t}-\mathbf{\hat{H}}_{4}^{\prime}\mathbf{f}_{t}^{\ast}$ & 4.9 & 5.0 & 5.1 & 5.1 & 5.2 & 5.6 \\
\ \ $\mathbf{\hat{f}}_{t}-\mathbf{\hat{H}}^{\prime}\mathbf{f}_{t}^{\ast}$ & 5.8 & 5.8 & 6.2 & 6.6 & 7.2 & 8.1 \\
\ \ $\mathbf{\hat{f}}_{t}-\mathbf{\hat{R}}_{\mathrm{Fan},F}^{\prime}\mathbf{f}_{t}^{\ast}$ & 4.5 & 4.8 & 4.7 & 4.8 & 4.8 & 4.9 \\
\multicolumn{7}{l}{} \\[-7pt]
\ \ $\mathbf{\hat{b}}_{i}-\mathbf{b}_{i}^{0}$ & 5.2 & 5.1 & 5.3 & 5.5 & 5.7 & 6.5 \\
\ \ $\mathbf{\hat{b}}_{i}-\mathbf{\hat{Q}b}_{i}^{\ast}$ & 5.1 & 5.2 & 5.5 & 6.0 & 6.8 & 9.6 \\
\ \ $\mathbf{\hat{b}}_{i}-\mathbf{\hat{H}}^{-1}\mathbf{b}_{i}^{\ast}$ & 5.5 & 6.2 & 8.1 & 10.7 & 18.0 & 35.1 \\
\ \ $\mathbf{\hat{b}}_{i}-\mathbf{\hat{R}}_{\mathrm{Fan},B}\mathbf{b}_{i}^{\ast}$ & 5.1 & 5.1 & 5.3 & 5.6 & 5.7 & 6.5 \\
\multicolumn{7}{l}{Test for linear restrictions (rejection frequencies of 5\% $N(0,1)$ tests)} \\
\ \ $H_{0}:f_{t,1}=f_{t,1}^{0}$ & 5.5 & 5.1 & 5.1 & 5.0 & 5.1 & 5.0 \\
\ \ $H_{0}:f_{t,2}=f_{t,2}^{0}$ & 4.5 & 4.8 & 4.9 & 4.8 & 5.0 & 5.1 \\
\multicolumn{7}{l}{} \\[-7pt]
\ \ $H_{0}:b_{i,1}=b_{i,1}^{0}$ & 5.1 & 5.1 & 5.3 & 5.4 & 5.2 & 5.5 \\
\ \ $H_{0}:b_{i,2}=b_{i,2}^{0}$ & 5.1 & 5.1 & 5.2 & 5.6 & 5.9 & 6.6 \\
\multicolumn{7}{l}{} \\[-7pt]
\ \ $H_{0}:c_{t,i}=c_{t,i}^{\ast}$ & 4.9 & 4.9 & 5.0 & 5.0 & 4.9 & 5.2 \\
\hline
\multicolumn{7}{c}{Factor Augmented Regressions} \\ \hline
\multicolumn{7}{l}{Estimation Accuracy (x1000)} \\
\ \ $||\boldsymbol{\hat{\gamma}}^{0}-\boldsymbol{\gamma}^{0}||_{\text{F}}$ & 150.2 & 150.0 & 149.4 & 150.1 & 150.1 & 150.0 \\
\ \ $||\boldsymbol{\hat{\gamma}}-\boldsymbol{\gamma}^{0}||_{\text{F}}$ & 150.5 & 150.4 & 149.9 & 151.4 & 153.3 & 157.0 \\
\ \ $||\boldsymbol{\hat{\gamma}}-\mathbf{\hat{H}}^{-1}\boldsymbol{\gamma}^{\ast}||_{\text{F}}$ & 152.4 & 154.5 & 159.2 & 169.0 & 191.5 & 238.3 \\
\multicolumn{7}{l}{Normal Approximation (rejection frequencies of 5\% $\chi_{r}^{2}$ tests)} \\
\ \ $\boldsymbol{\hat{\gamma}}^{0}-\boldsymbol{\gamma}^{0}$ & 5.5 & 5.6 & 5.4 & 5.5 & 5.6 & 5.5 \\
\ \ $\boldsymbol{\hat{\gamma}}-\boldsymbol{\gamma}^{0}$ & 5.5 & 5.8 & 5.5 & 5.6 & 6.0 & 6.8 \\
\ \ $\boldsymbol{\hat{\gamma}}-\mathbf{\hat{H}}^{-1}\boldsymbol{\gamma}^{\ast}$ & 5.9 & 6.6 & 7.4 & 9.8 & 15.8 & 30.6 \\
\multicolumn{7}{l}{Test for linear restrictions (rejection frequencies of 5\% $N(0,1)$ tests)} \\
\ \ $\hat{\gamma}_{1}^{0}-\gamma_{1}^{0}$ & 5.2 & 5.4 & 5.2 & 5.4 & 5.3 & 5.3 \\
\ \ $\hat{\gamma}_{1}-\gamma_{1}^{0}$ & 5.3 & 5.4 & 5.2 & 5.5 & 5.7 & 5.8 \\
\multicolumn{7}{l}{} \\[-7pt]
\ \ $\hat{\gamma}_{2}^{0}-\gamma_{2}^{0}$ & 5.2 & 5.4 & 5.1 & 5.1 & 5.2 & 5.2 \\
\ \ $\hat{\gamma}_{2}-\gamma_{2}^{0}$ & 5.2 & 5.5 & 5.3 & 5.4 & 5.8 & 6.9 \\
\multicolumn{7}{l}{Empirical coverage} \\
\ \ \ 95\% $CI(y_{T+h|T})$ & 95.1 & 95.0 & 94.7 & 94.7 & 94.7 & 94.4 \\
\ \ \ 95\% $CI(y_{T+h})$ & 95.0 & 95.2 & 95.1 & 95.0 & 94.9 & 94.9 \\
\multicolumn{7}{l}{} \\[-7pt]
\ \ \ 90\% $CI(y_{T+h|T})$ & 90.0 & 90.0 & 89.6 & 89.6 & 89.5 & 88.9 \\
\ \ \ 90\% $CI(y_{T+h})$ & 90.1 & 90.0 & 90.1 & 90.0 & 89.8 & 90.0 \\
\hline
\end{tabular}
\label{table:sparseNT100}
\end{table}

\clearpage

\renewcommand{\arraystretch}{.9}
\begin{table}[htb!]
\caption{Experimental results with sparse factor loadings for $N=T=200$}
\centering\small
\begin{tabular}{lrrrrrr}
\hline
$(\alpha_{1},\alpha_{2})$ & (1.0,1.0) & (1.0,0.9) & (1.0,0.8) & (0.9,0.7) & (0.8,0.6) & (0.7,0.5) \\ \hline
\multicolumn{7}{c}{Factor models} \\ \hline
\multicolumn{7}{l}{Estimation Accuracy (x1000)} \\
\ \ $T^{-\frac{1}{2}}||\mathbf{\hat{F}}-\mathbf{F}^{0}||_{\text{F}}$ & 55.7 & 69.4 & 87.7 & 113.4 & 144.5 & 184.6 \\
\ \ $T^{-\frac{1}{2}}||\mathbf{\hat{F}}-\mathbf{F}^{\ast}\mathbf{\hat{H}}_{4}||_{\text{F}}$ & 55.9 & 69.8 & 88.3 & 114.3 & 145.9 & 187.0 \\
\ \ $T^{-\frac{1}{2}}||\mathbf{\hat{F}}-\mathbf{F}^{\ast}\mathbf{\hat{H}}||_{\text{F}}$ & 56.9 & 71.6 & 91.8 & 119.4 & 153.3 & 198.1 \\
\ \ $T^{-\frac{1}{2}}||\mathbf{\hat{F}}-\mathbf{F}^{\ast}\mathbf{\hat{R}}_{\mathrm{Fan},F}||_{\text{F}}$ & 55.6 & 69.4 & 87.7 & 113.3 & 144.4 & 184.5 \\
\multicolumn{7}{l}{} \\[-7pt]
\ \ $N^{-\frac{1}{2}}||\mathbf{\hat{B}}-\mathbf{B}^{0}||_{\text{F}}$ & 70.8 & 70.8 & 70.9 & 71.1 & 71.4 & 71.9 \\
\ \ $N^{-\frac{1}{2}}||\mathbf{\hat{B}}-\mathbf{B}^{\ast}\mathbf{\hat{Q}}^{\prime}||_{\text{F}}$ & 70.8 & 70.9 & 71.1 & 71.3 & 71.7 & 72.4 \\
\ \ $N^{-\frac{1}{2}}||\mathbf{\hat{B}}-\mathbf{B}^{\ast}\mathbf{\hat{H}}^{\prime-1}||_{\text{F}}$ & 71.2 & 71.6 & 72.3 & 72.9 & 73.7 & 75.1 \\
\ \ $N^{-\frac{1}{2}}||\mathbf{\hat{B}}-\mathbf{B}^{\ast}\mathbf{\hat{R}}_{\mathrm{Fan},B}||_{\text{F}}$ & 70.8 & 70.9 & 71.0 & 71.2 & 71.5 & 72.0 \\
\multicolumn{7}{l}{} \\[-7pt]
\ \ $(NT)^{-\frac{1}{2}}||\mathbf{\hat{C}}-\mathbf{C}^{\ast}||_{\text{F}}$ & 99.8 & 99.9 & 100.0 & 100.2 & 100.5 & 100.9 \\

\multicolumn{7}{l}{Normal Approximation (rejection frequencies of 5\% $\chi_{r}^{2}$ tests)} \\
\ \ $\mathbf{\hat{f}}_{t}-\mathbf{f}_{t}^{0}$ & 4.9 & 5.0 & 4.9 & 5.2 & 5.1 & 4.9 \\
\ \ $\mathbf{\hat{f}}_{t}-\mathbf{\hat{H}}_{4}^{\prime}\mathbf{f}_{t}^{\ast}$ & 4.9 & 5.0 & 4.9 & 5.2 & 5.1 & 4.9 \\
\ \ $\mathbf{\hat{f}}_{t}-\mathbf{\hat{H}}^{\prime}\mathbf{f}_{t}^{\ast}$ & 4.9 & 5.0 & 4.9 & 5.2 & 5.2 & 4.9 \\
\ \ $\mathbf{\hat{f}}_{t}-\mathbf{\hat{R}}_{\mathrm{Fan},F}^{\prime}\mathbf{f}_{t}^{\ast}$ & 4.9 & 4.9 & 4.9 & 5.2 & 5.1 & 4.9 \\
\multicolumn{7}{l}{} \\[-7pt]
\ \ $\mathbf{\hat{b}}_{i}-\mathbf{b}_{i}^{0}$ & 4.8 & 5.2 & 5.1 & 5.4 & 5.5 & 6.2 \\
\ \ $\mathbf{\hat{b}}_{i}-\mathbf{\hat{Q}b}_{i}^{\ast}$ & 4.9 & 5.2 & 5.2 & 5.7 & 6.3 & 8.6 \\
\ \ $\mathbf{\hat{b}}_{i}-\mathbf{\hat{H}}^{-1}\mathbf{b}_{i}^{\ast}$ & 5.0 & 5.7 & 6.9 & 9.1 & 14.0 & 26.1 \\
\ \ $\mathbf{\hat{b}}_{i}-\mathbf{\hat{R}}_{\mathrm{Fan},B}'\mathbf{b}_{i}^{\ast}$ & 4.8 & 5.2 & 5.1 & 5.4 & 5.5 & 6.3 \\
\multicolumn{7}{l}{Test for linear restrictions (rejection frequencies of 5\% $N(0,1)$ tests)} \\
\ \ $H_{0}:f_{t,1}=f_{t,1}^{0}$ & 5.0 & 5.0 & 5.0 & 5.0 & 4.9 & 5.1 \\
\ \ $H_{0}:f_{t,2}=f_{t,2}^{0}$ & 4.9 & 4.8 & 4.8 & 5.1 & 5.0 & 4.9 \\
\multicolumn{7}{l}{} \\[-7pt]
\ \ $H_{0}:b_{i,1}=b_{i,1}^{0}$ & 5.0 & 5.2 & 5.1 & 5.0 & 5.1 & 5.2 \\
\ \ $H_{0}:b_{i,2}=b_{i,2}^{0}$ & 4.8 & 4.9 & 5.1 & 5.4 & 5.6 & 6.4 \\
\multicolumn{7}{l}{} \\[-7pt]
\ \ $H_{0}:c_{t,i}=c_{t,i}^{\ast}$ & 5.1 & 5.0 & 5.1 & 5.3 & 5.1 & 5.1 \\
\hline
\multicolumn{7}{c}{Factor Augmented Regressions} \\ \hline
\multicolumn{7}{l}{Estimation Accuracy (x1000)} \\
\ \ $||\boldsymbol{\hat{\gamma}}^{0}-\boldsymbol{\gamma}^{0}||_{\text{F}}$ & 107.0 & 106.8 & 107.3 & 107.3 & 106.2 & 106.8 \\
\ \ $||\boldsymbol{\hat{\gamma}}-\boldsymbol{\gamma}^{0}||_{\text{F}}$ & 107.0 & 106.9 & 107.6 & 107.9 & 107.8 & 111.1 \\
\ \ $||\boldsymbol{\hat{\gamma}}-\mathbf{\hat{H}}^{-1}\boldsymbol{\gamma}^{\ast}||_{\text{F}}$ & 107.7 & 108.5 & 111.9 & 117.3 & 127.7 & 155.4 \\
\multicolumn{7}{l}{Normal Approximation (rejection frequencies of 5\% $\chi_{r}^{2}$ tests)} \\
\ \ $\boldsymbol{\hat{\gamma}}^{0}-\boldsymbol{\gamma}^{0}$ & 5.3 & 5.0 & 5.3 & 5.3 & 5.0 & 5.0 \\
\ \ $\boldsymbol{\hat{\gamma}}-\boldsymbol{\gamma}^{0}$ & 5.4 & 5.1 & 5.3 & 5.4 & 5.5 & 6.1 \\
\ \ $\boldsymbol{\hat{\gamma}}-\mathbf{\hat{H}}^{-1}\boldsymbol{\gamma}^{\ast}$ & 5.5 & 5.5 & 6.7 & 8.6 & 12.2 & 23.4 \\
\multicolumn{7}{l}{Test for linear restrictions (rejection frequencies of 5\% $N(0,1)$ tests)} \\
\ \ $\hat{\gamma}_{1}^{0}-\gamma_{1}^{0}$ & 5.2 & 5.1 & 5.3 & 5.2 & 5.1 & 5.1 \\
\ \ $\hat{\gamma}_{1}-\gamma_{1}^{0}$ & 5.2 & 5.1 & 5.3 & 5.4 & 5.4 & 5.7 \\
\multicolumn{7}{l}{} \\[-7pt]
\ \ $\hat{\gamma}_{2}^{0}-\gamma_{2}^{0}$ & 5.3 & 4.9 & 5.2 & 5.3 & 5.0 & 5.1 \\
\ \ $\hat{\gamma}_{2}-\gamma_{2}^{0}$ & 5.3 & 5.0 & 5.4 & 5.5 & 5.5 & 6.6 \\
\multicolumn{7}{l}{Empirical coverage} \\
\ \ \ 95\% $CI(y_{T+h|T})$ & 94.9 & 94.8 & 94.9 & 94.8 & 94.8 & 94.6 \\
\ \ \ 95\% $CI(y_{T+h})$ & 95.0 & 95.1 & 95.3 & 94.9 & 95.1 & 95.0 \\
\multicolumn{7}{l}{} \\[-7pt]
\ \ \ 90\% $CI(y_{T+h|T})$ & 89.9 & 89.8 & 89.8 & 89.6 & 89.6 & 89.5 \\
\ \ \ 90\% $CI(y_{T+h})$ & 90.2 & 90.1 & 90.3 & 90.1 & 90.2 & 90.0 \\
\hline
\end{tabular}
\label{table:sparseNT200}
\end{table}

\clearpage

\renewcommand{\arraystretch}{.9}
\begin{table}[htb!]
\caption{Experimental results with sparse factor loadings for $N=T=500$}
\centering\small
\begin{tabular}{lrrrrrr}
\hline
$(\alpha_{1},\alpha_{2})$ & (1.0,1.0) & (1.0,0.9) & (1.0,0.8) & (0.9,0.7) & (0.8,0.6) & (0.7,0.5) \\ \hline
\multicolumn{7}{c}{Factor models} \\ \hline
\multicolumn{7}{l}{Estimation Accuracy (x1000)} \\
\ \ $T^{-\frac{1}{2}}||\mathbf{\hat{F}}-\mathbf{F}^{0}||_{\text{F}}$ & 35.3 & 45.8 & 60.5 & 82.8 & 112.9 & 149.8 \\
\ \ $T^{-\frac{1}{2}}||\mathbf{\hat{F}}-\mathbf{F}^{\ast}\mathbf{\hat{H}}_{4}||_{\text{F}}$ & 35.4 & 45.9 & 60.7 & 83.1 & 113.5 & 151.0 \\
\ \ $T^{-\frac{1}{2}}||\mathbf{\hat{F}}-\mathbf{F}^{\ast}\mathbf{\hat{H}}||_{\text{F}}$ & 35.6 & 46.4 & 61.9 & 85.1 & 116.9 & 156.8 \\
\ \ $T^{-\frac{1}{2}}||\mathbf{\hat{F}}-\mathbf{F}^{\ast}\mathbf{\hat{R}}_{\mathrm{Fan},F}||_{\text{F}}$ & 35.3 & 45.7 & 60.5 & 82.8 & 112.9 & 149.8 \\
\multicolumn{7}{l}{} \\[-7pt]
\ \ $N^{-\frac{1}{2}}||\mathbf{\hat{B}}-\mathbf{B}^{0}||_{\text{F}}$ & 44.8 & 44.8 & 44.8 & 44.9 & 45.0 & 45.3 \\
\ \ $N^{-\frac{1}{2}}||\mathbf{\hat{B}}-\mathbf{B}^{\ast}\mathbf{\hat{Q}}^{\prime}||_{\text{F}}$ & 44.8 & 44.8 & 44.8 & 44.9 & 45.1 & 45.4 \\
\ \ $N^{-\frac{1}{2}}||\mathbf{\hat{B}}-\mathbf{B}^{\ast}\mathbf{\hat{H}}^{\prime-1}||_{\text{F}}$ & 44.8 & 45.0 & 45.2 & 45.4 & 45.8 & 46.5 \\
\ \ $N^{-\frac{1}{2}}||\mathbf{\hat{B}}-\mathbf{B}^{\ast}\mathbf{\hat{R}}_{\mathrm{Fan},B}||_{\text{F}}$ & 44.7 & 44.8 & 44.8 & 44.9 & 45.0 & 45.3 \\
\multicolumn{7}{l}{} \\[-7pt]
\ \ $(NT)^{-\frac{1}{2}}||\mathbf{\hat{C}}-\mathbf{C}^{\ast}||_{\text{F}}$ & 63.2 & 63.2 & 63.3 & 63.3 & 63.5 & 63.7 \\

\multicolumn{7}{l}{Normal Approximation (rejection frequencies of 5\% $\chi_{r}^{2}$ tests)} \\
\ \ $\mathbf{\hat{f}}_{t}-\mathbf{f}_{t}^{0}$ & 4.9 & 4.8 & 5.0 & 4.9 & 5.0 & 5.2 \\
\ \ $\mathbf{\hat{f}}_{t}-\mathbf{\hat{H}}_{4}^{\prime}\mathbf{f}_{t}^{\ast}$ & 4.9 & 4.9 & 5.0 & 5.0 & 5.1 & 5.2 \\
\ \ $\mathbf{\hat{f}}_{t}-\mathbf{\hat{H}}^{\prime}\mathbf{f}_{t}^{\ast}$ & 4.9 & 5.0 & 5.1 & 5.1 & 5.3 & 5.4 \\
\ \ $\mathbf{\hat{f}}_{t}-\mathbf{\hat{R}}_{\mathrm{Fan},F}^{\prime}\mathbf{f}_{t}^{\ast}$ & 4.9 & 4.8 & 5.0 & 4.9 & 5.0 & 5.2 \\
\multicolumn{7}{l}{} \\[-7pt]
\ \ $\mathbf{\hat{b}}_{i}-\mathbf{b}_{i}^{0}$ & 5.0 & 5.1 & 5.3 & 5.1 & 5.6 & 6.3 \\
\ \ $\mathbf{\hat{b}}_{i}-\mathbf{\hat{Q}b}_{i}^{\ast}$ & 5.0 & 5.1 & 5.3 & 5.4 & 6.6 & 8.8 \\
\ \ $\mathbf{\hat{b}}_{i}-\mathbf{\hat{H}}^{-1}\mathbf{b}_{i}^{\ast}$ & 5.1 & 5.4 & 6.3 & 7.8 & 12.8 & 26.5 \\
\ \ $\mathbf{\hat{b}}_{i}-\mathbf{\hat{R}}_{\mathrm{Fan},B}\mathbf{b}_{i}^{\ast}$ & 5.0 & 5.1 & 5.3 & 5.1 & 5.6 & 6.3 \\
\multicolumn{7}{l}{Test for linear restrictions (rejection frequencies of 5\% $N(0,1)$ tests)} \\
\ \ $H_{0}:f_{t,1}=f_{t,1}^{0}$ & 4.9 & 4.9 & 4.9 & 5.1 & 5.0 & 5.1 \\
\ \ $H_{0}:f_{t,2}=f_{t,2}^{0}$ & 4.9 & 4.9 & 5.0 & 4.9 & 5.2 & 5.0 \\
\multicolumn{7}{l}{} \\[-7pt]
\ \ $H_{0}:b_{i,1}=b_{i,1}^{0}$ & 4.9 & 5.1 & 5.2 & 5.1 & 5.2 & 5.4 \\
\ \ $H_{0}:b_{i,2}=b_{i,2}^{0}$ & 5.1 & 5.0 & 5.3 & 5.2 & 5.8 & 6.5 \\
\multicolumn{7}{l}{} \\[-7pt]
\ \ $H_{0}:c_{t,i}=c_{t,i}^{\ast}$ & 4.9 & 5.0 & 5.1 & 5.1 & 5.3 & 5.1 \\
\hline
\multicolumn{7}{c}{Factor Augmented Regressions} \\ \hline
\multicolumn{7}{l}{Estimation Accuracy (x1000)} \\
\ \ $||\boldsymbol{\hat{\gamma}}^{0}-\boldsymbol{\gamma}^{0}||_{\text{F}}$ & 67.0 & 67.1 & 67.4 & 67.2 & 67.4 & 67.3 \\
\ \ $||\boldsymbol{\hat{\gamma}}-\boldsymbol{\gamma}^{0}||_{\text{F}}$ & 67.0 & 67.1 & 67.4 & 67.5 & 68.3 & 70.0 \\
\ \ $||\boldsymbol{\hat{\gamma}}-\mathbf{\hat{H}}^{-1}\boldsymbol{\gamma}^{\ast}||_{\text{F}}$ & 67.2 & 67.6 & 68.9 & 71.5 & 79.7 & 98.8 \\
\multicolumn{7}{l}{Normal Approximation (rejection frequencies of 5\% $\chi_{r}^{2}$ tests)} \\
\ \ $\boldsymbol{\hat{\gamma}}^{0}-\boldsymbol{\gamma}^{0}$ & 5.0 & 5.1 & 5.2 & 5.2 & 5.1 & 5.0 \\
\ \ $\boldsymbol{\hat{\gamma}}-\boldsymbol{\gamma}^{0}$ & 5.1 & 5.1 & 5.2 & 5.3 & 5.4 & 6.2 \\
\ \ $\boldsymbol{\hat{\gamma}}-\mathbf{\hat{H}}^{-1}\boldsymbol{\gamma}^{\ast}$ & 5.1 & 5.4 & 5.8 & 7.3 & 11.6 & 24.0 \\
\multicolumn{7}{l}{Test for linear restrictions (rejection frequencies of 5\% $N(0,1)$ tests)} \\
\ \ $\hat{\gamma}_{1}^{0}-\gamma_{1}^{0}$ & 5.1 & 4.9 & 5.1 & 5.0 & 5.0 & 4.9 \\
\ \ $\hat{\gamma}_{1}-\gamma_{1}^{0}$ & 5.1 & 5.0 & 5.1 & 5.2 & 5.1 & 5.4 \\
\multicolumn{7}{l}{} \\[-7pt]
\ \ $\hat{\gamma}_{2}^{0}-\gamma_{2}^{0}$ & 5.0 & 5.1 & 5.1 & 5.0 & 5.1 & 5.2 \\
\ \ $\hat{\gamma}_{2}-\gamma_{2}^{0}$ & 5.0 & 5.2 & 5.1 & 5.1 & 5.7 & 6.7 \\
\multicolumn{7}{l}{Empirical coverage} \\
\ \ \ 95\% $CI(y_{T+h|T})$ & 95.0 & 95.0 & 94.8 & 94.5 & 94.4 & 93.4 \\
\ \ \ 95\% $CI(y_{T+h})$ & 94.9 & 94.9 & 95.2 & 94.9 & 95.2 & 95.0 \\
\multicolumn{7}{l}{} \\[-7pt]
\ \ \ 90\% $CI(y_{T+h|T})$ & 90.1 & 90.1 & 89.9 & 89.3 & 89.3 & 87.8 \\
\ \ \ 90\% $CI(y_{T+h})$ & 90.0 & 90.0 & 90.2 & 89.8 & 90.2 & 90.0 \\
\hline
\end{tabular}
\label{table:sparseNT500}
\end{table}

\section{Full Results under the Rotation of \citet{FanEtAl2024}}
\label{sec:D_fan}
\setcounter{figure}{0}
\setcounter{table}{0}
\renewcommand{\thefigure}{D.\arabic{figure}}
\renewcommand{\thetable}{D.\arabic{table}}

This section reports the complete set of Monte Carlo results for the
Procrustes rotation of \citet{FanEtAl2024}, defined in Section
\ref{mc: fac des} of the main text, complementing the summary in Section
\ref{sec: MC}. The data generating process and the number of replications
(50{,}000) are identical to those in the main text; we report the
non-sparse loading design, the results for the sparse design being very
similar.

\subsection{Size of the $z$-tests under all centerings}

\begin{figure}[htbp!]
\centering
\includegraphics[width=\textwidth]{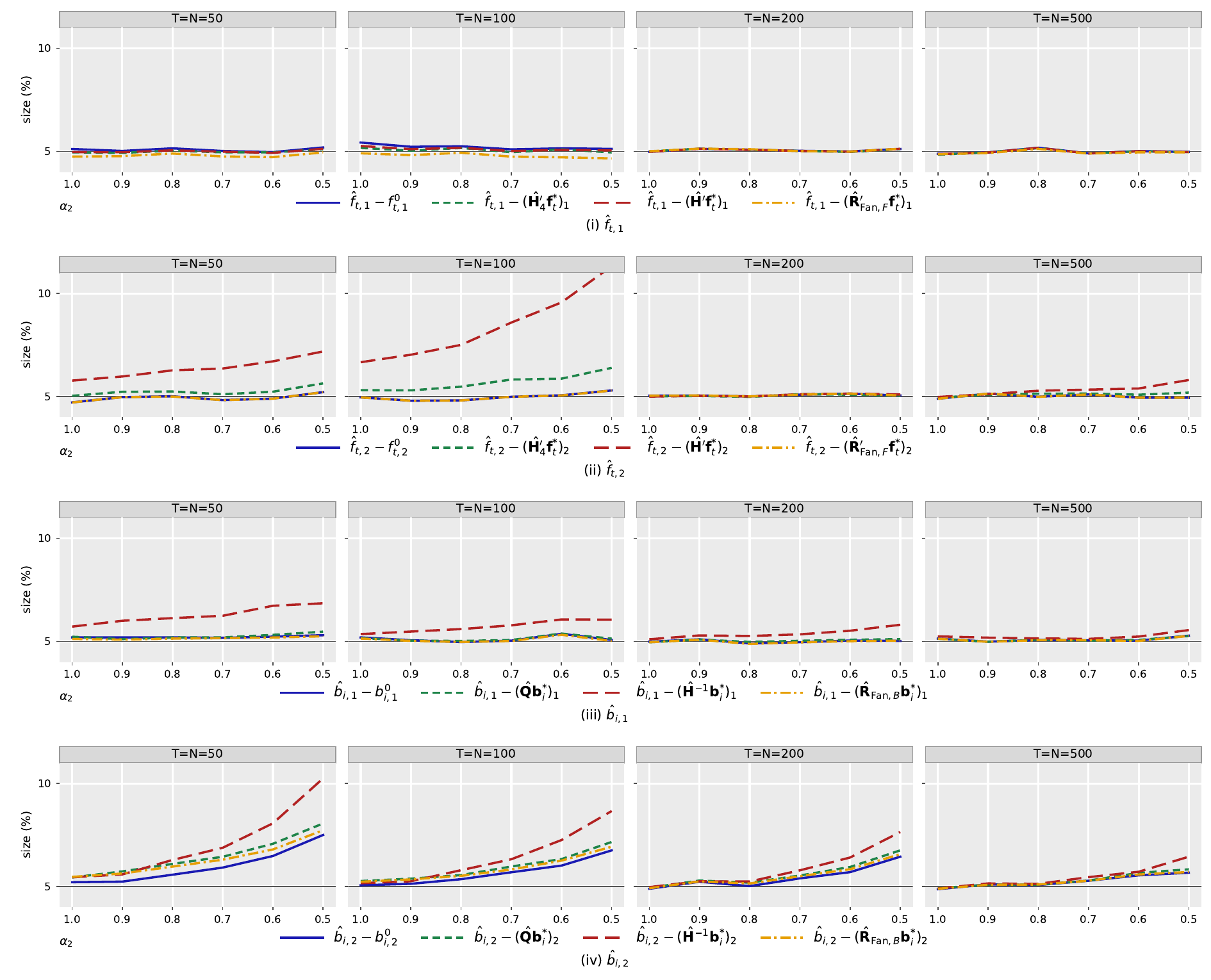}
\caption{Size of the 5\% level $z$-tests under all centerings:
(i) $\hat{f}_{t,1}$; (ii) $\hat{f}_{t,2}$; (iii) $\hat{b}_{i,1}$;
(iv) $\hat{b}_{i,2}$. The legend indicates the centering of the test
statistic.}
\label{fig:D_ztests}
\end{figure}

Figure \ref{fig:D_ztests} reports the empirical size of the 5\% level
$z$-tests for each element of the factors and the loadings under four
centerings: the pseudo-true parameters $(f_{t,k}^{0}, b_{i,k}^{0})$, and
the data-dependent targets implied by $(\hat{\bH}_4,\hat{\bQ})$,
$\hat{\bH}$, and
$(\hat{\bR}_{\mathrm{Fan},F},\hat{\bR}_{\mathrm{Fan},B})$. Three features
stand out. First, the \citet{FanEtAl2024} centering behaves essentially
identically to the pseudo-true centering for both factors throughout, and
very similarly for the loadings. Second, centering at
$\hat{\bH}^{\prime}\bff_t^{*}$ or $\hat{\bH}^{-1}\bb_i^{*}$ leads to
substantial over-rejection for the second (weaker) factor and its
loadings, which grows as the model weakens; the
$(\hat{\bH}_4,\hat{\bQ})$ centering lies in between. Third, all
centerings improve as the sample size grows. We stress that only the
pseudo-true centering corresponds to a fixed null hypothesis; the
data-dependent centerings change with the sample and are reported purely
as approximation diagnostics.

\subsection{Average rotation matrices}

Table~\ref{tab:D_rotations} reports the Monte Carlo averages of the
data-dependent rotation matrices $\hat{\bH}$, $\hat{\bH}_4$,
$\hat{\bH}_3=\hat{\bQ}^{-1}$, and $\hat{\bR}_{\mathrm{Fan},F}$, all of which
estimate the fixed rotation $\bH$; Table~\ref{tab:D_rotations_sd} reports the
corresponding Monte Carlo standard deviations. Three features stand out.

First, the Monte Carlo mean of $\hat{\bR}_{\mathrm{Fan},F}$ is numerically
very close to $\bH$ in every model and at every sample size, including
$T=N=50$: the largest entrywise deviation of the mean across all designs is
below $5\times10^{-4}$, while the across-replication standard deviations of
its entries range from about $5\times10^{-3}$ at $T=N=500$ to about
$4\times10^{-2}$ at $T=N=50$ in the weakest design. 
That is, the Procrustes alignment is, in
effect, recovering the fixed rotation that maps the original representation
$(\bF^{*},\bB^{*})$ to the PC-normalized population representation
$(\bF^{0},\bB^{0})$. This is why the targets of \citet{FanEtAl2024} track the
pseudo-true targets so closely in Figures~\ref{fig:1}, \ref{fig:2}, and
\ref{fig:D_ztests}. The reason is that, in our design, the Fan et al.\ alignment
and the fixed rotation $\bH$ align the original representation with the same
PC-normalized population representation. Writing the scaled common component as
$\bC^{*}/\sqrt{T}=\bU^0\bLambda^{1/2}\bV^{0\prime}$ with
$\bU^{0\prime}\bU^0=\bV^{0\prime}\bV^0=\bI_r$, the singular value decomposition
recovers, up to column sign changes, the PC-normalized representation
$\bF^0=\sqrt{T}\,\bU^0$ and $\bB^0=\bV^0\bLambda^{1/2}$, which satisfies
$T^{-1}\bF^{0\prime}\bF^0=\bI_r$ and $\bB^{0\prime}\bB^0=\bLambda$. Since
$\bF^0=\bF^*\bH$, the population Procrustes rotation on the factor side is $\bH$,
so $\hat{\bR}_{\mathrm{Fan},F}$ is close to $\bH$ whenever the sample singular
vectors are accurately estimated.

Second, the other three rotations are biased in the entries associated with the
weak factor, and the biases are systematic in sign. Both $\hat{\bH}$ and
$\hat{\bH}_4$ are biased \emph{downward}, while $\hat{\bH}_3$ is biased
\emph{upward}; in all three cases the distortion is concentrated in the $(2,2)$
element and, to a lesser extent, the $(1,2)$ element, and it grows as the model
weakens. The magnitudes are ordered $\hat{\bH}>\hat{\bH}_4>\hat{\bH}_3$
uniformly: in the weakest design at $T=N=50$, the $(2,2)$ entry deviates from its
target of $2$ by $-0.46$ for $\hat{\bH}$, $-0.29$ for $\hat{\bH}_4$, and $+0.12$
for $\hat{\bH}_3$. These biases shrink only slowly with the sample size---roughly
halving as $N$ rises tenfold---so they remain visible even at $T=N=500$. This is
the source of the poorer performance of the $\hat{\bH}$- and $\hat{\bH}_4$-based
targets in the weak models, and explains why the Procrustes losses are almost
identical to the pseudo-true losses in Figure~\ref{fig:1} while the $\hat{\bH}$-
and $\hat{\bH}_4$-based losses can behave differently.

Third, $\hat{\bR}_{\mathrm{Fan},F}$ is also the most stable of the four across
replications. Its overall standard deviation is the smallest in every design, and
the dispersion ordering $\hat{\bR}_{\mathrm{Fan},F}<\hat{\bH}_3<\hat{\bH}_4<\hat{\bH}$
mirrors the bias ordering above: the rotation closest to $\bH$ on average is also
the one that varies least around it. The only entry in which
$\hat{\bR}_{\mathrm{Fan},F}$ is not the tightest is the $(2,1)$ element, where the
more concentrated $\hat{\bH}$ and $\hat{\bH}_4$ have smaller standard deviations;
in summed and worst-entry terms $\hat{\bR}_{\mathrm{Fan},F}$ remains the most
concentrated throughout.

Finally, since the columns of $\hat{\bF}$ are sign- and order-adjusted to
$\bF^{0}$ as described in Section~\ref{mc: fac des}, the rotations are reported
under the same column-sign convention as $\bH$; without this adjustment, the
near-equality holds only up to the usual column sign indeterminacy.

\begin{table}[htbp!]
\centering\scriptsize
\setlength{\tabcolsep}{3pt}
\caption{Monte Carlo averages of the data-dependent rotation matrices (non-sparse loadings, 50{,}000 replications). All of $\hat{\bH}$, $\hat{\bH}_4$, $\hat{\bH}_3=\hat{\bQ}^{-1}$, and $\hat{\bR}_{\mathrm{Fan},F}$ estimate the fixed rotation $\bH=\bigl(\begin{smallmatrix}1&0.5\\0.5&2\end{smallmatrix}\bigr)$.}
\label{tab:D_rotations}
\resizebox{\textwidth}{!}{%
\begin{tabular}{clcccccc}
\hline\hline
$T=N$ & $(\alpha_1,\alpha_2)$: & (1.0,1.0) & (1.0,0.9) & (1.0,0.8) & (0.9,0.7) & (0.8,0.6) & (0.7,0.5) \\
\hline
50 & $\overline{\hat{\bH}}$ & $\bigl(\begin{smallmatrix}0.990&0.481\\0.495&1.926\end{smallmatrix}\bigr)$ & $\bigl(\begin{smallmatrix}0.990&0.473\\0.495&1.893\end{smallmatrix}\bigr)$ & $\bigl(\begin{smallmatrix}0.990&0.460\\0.495&1.844\end{smallmatrix}\bigr)$ & $\bigl(\begin{smallmatrix}0.986&0.444\\0.493&1.775\end{smallmatrix}\bigr)$ & $\bigl(\begin{smallmatrix}0.979&0.419\\0.490&1.677\end{smallmatrix}\bigr)$ & $\bigl(\begin{smallmatrix}0.970&0.387\\0.485&1.545\end{smallmatrix}\bigr)$ \\[3pt]
 & $\overline{\hat{\bH}}_{4}$ & $\bigl(\begin{smallmatrix}0.994&0.489\\0.497&1.955\end{smallmatrix}\bigr)$ & $\bigl(\begin{smallmatrix}0.994&0.483\\0.497&1.934\end{smallmatrix}\bigr)$ & $\bigl(\begin{smallmatrix}0.994&0.476\\0.497&1.903\end{smallmatrix}\bigr)$ & $\bigl(\begin{smallmatrix}0.992&0.465\\0.496&1.859\end{smallmatrix}\bigr)$ & $\bigl(\begin{smallmatrix}0.987&0.449\\0.494&1.795\end{smallmatrix}\bigr)$ & $\bigl(\begin{smallmatrix}0.981&0.427\\0.491&1.706\end{smallmatrix}\bigr)$ \\[3pt]
 & $\overline{\hat{\bH}}_{3}\,(=\hat{\bQ}^{-1})$ & $\bigl(\begin{smallmatrix}1.002&0.504\\0.501&2.015\end{smallmatrix}\bigr)$ & $\bigl(\begin{smallmatrix}1.002&0.506\\0.501&2.023\end{smallmatrix}\bigr)$ & $\bigl(\begin{smallmatrix}1.002&0.509\\0.501&2.034\end{smallmatrix}\bigr)$ & $\bigl(\begin{smallmatrix}1.003&0.513\\0.501&2.052\end{smallmatrix}\bigr)$ & $\bigl(\begin{smallmatrix}1.004&0.519\\0.502&2.078\end{smallmatrix}\bigr)$ & $\bigl(\begin{smallmatrix}1.006&0.530\\0.503&2.119\end{smallmatrix}\bigr)$ \\[3pt]
 & $\overline{\hat{\bR}}_{\mathrm{Fan},F}$ & $\bigl(\begin{smallmatrix}1.000&0.500\\0.500&2.000\end{smallmatrix}\bigr)$ & $\bigl(\begin{smallmatrix}1.000&0.500\\0.500&2.000\end{smallmatrix}\bigr)$ & $\bigl(\begin{smallmatrix}1.000&0.500\\0.500&2.000\end{smallmatrix}\bigr)$ & $\bigl(\begin{smallmatrix}1.000&0.500\\0.500&2.000\end{smallmatrix}\bigr)$ & $\bigl(\begin{smallmatrix}1.000&0.500\\0.500&2.000\end{smallmatrix}\bigr)$ & $\bigl(\begin{smallmatrix}1.000&0.500\\0.500&2.000\end{smallmatrix}\bigr)$ \\[3pt]
\hline
100 & $\overline{\hat{\bH}}$ & $\bigl(\begin{smallmatrix}0.995&0.490\\0.498&1.962\end{smallmatrix}\bigr)$ & $\bigl(\begin{smallmatrix}0.995&0.485\\0.498&1.940\end{smallmatrix}\bigr)$ & $\bigl(\begin{smallmatrix}0.995&0.476\\0.498&1.906\end{smallmatrix}\bigr)$ & $\bigl(\begin{smallmatrix}0.992&0.463\\0.496&1.853\end{smallmatrix}\bigr)$ & $\bigl(\begin{smallmatrix}0.988&0.444\\0.494&1.773\end{smallmatrix}\bigr)$ & $\bigl(\begin{smallmatrix}0.981&0.414\\0.490&1.654\end{smallmatrix}\bigr)$ \\[3pt]
 & $\overline{\hat{\bH}}_{4}$ & $\bigl(\begin{smallmatrix}0.997&0.494\\0.499&1.977\end{smallmatrix}\bigr)$ & $\bigl(\begin{smallmatrix}0.997&0.491\\0.499&1.963\end{smallmatrix}\bigr)$ & $\bigl(\begin{smallmatrix}0.997&0.486\\0.499&1.942\end{smallmatrix}\bigr)$ & $\bigl(\begin{smallmatrix}0.995&0.477\\0.498&1.910\end{smallmatrix}\bigr)$ & $\bigl(\begin{smallmatrix}0.993&0.465\\0.496&1.859\end{smallmatrix}\bigr)$ & $\bigl(\begin{smallmatrix}0.988&0.445\\0.494&1.781\end{smallmatrix}\bigr)$ \\[3pt]
 & $\overline{\hat{\bH}}_{3}\,(=\hat{\bQ}^{-1})$ & $\bigl(\begin{smallmatrix}1.001&0.502\\0.501&2.008\end{smallmatrix}\bigr)$ & $\bigl(\begin{smallmatrix}1.001&0.503\\0.500&2.013\end{smallmatrix}\bigr)$ & $\bigl(\begin{smallmatrix}1.001&0.505\\0.500&2.020\end{smallmatrix}\bigr)$ & $\bigl(\begin{smallmatrix}1.002&0.508\\0.501&2.032\end{smallmatrix}\bigr)$ & $\bigl(\begin{smallmatrix}1.002&0.513\\0.501&2.051\end{smallmatrix}\bigr)$ & $\bigl(\begin{smallmatrix}1.004&0.521\\0.502&2.083\end{smallmatrix}\bigr)$ \\[3pt]
 & $\overline{\hat{\bR}}_{\mathrm{Fan},F}$ & $\bigl(\begin{smallmatrix}1.000&0.500\\0.500&2.000\end{smallmatrix}\bigr)$ & $\bigl(\begin{smallmatrix}1.000&0.500\\0.500&2.000\end{smallmatrix}\bigr)$ & $\bigl(\begin{smallmatrix}1.000&0.500\\0.500&2.000\end{smallmatrix}\bigr)$ & $\bigl(\begin{smallmatrix}1.000&0.500\\0.500&2.000\end{smallmatrix}\bigr)$ & $\bigl(\begin{smallmatrix}1.000&0.500\\0.500&2.000\end{smallmatrix}\bigr)$ & $\bigl(\begin{smallmatrix}1.000&0.500\\0.500&2.000\end{smallmatrix}\bigr)$ \\[3pt]
\hline
200 & $\overline{\hat{\bH}}$ & $\bigl(\begin{smallmatrix}0.998&0.495\\0.499&1.980\end{smallmatrix}\bigr)$ & $\bigl(\begin{smallmatrix}0.998&0.492\\0.499&1.967\end{smallmatrix}\bigr)$ & $\bigl(\begin{smallmatrix}0.998&0.486\\0.499&1.944\end{smallmatrix}\bigr)$ & $\bigl(\begin{smallmatrix}0.996&0.476\\0.498&1.906\end{smallmatrix}\bigr)$ & $\bigl(\begin{smallmatrix}0.993&0.461\\0.496&1.844\end{smallmatrix}\bigr)$ & $\bigl(\begin{smallmatrix}0.988&0.436\\0.494&1.744\end{smallmatrix}\bigr)$ \\[3pt]
 & $\overline{\hat{\bH}}_{4}$ & $\bigl(\begin{smallmatrix}0.999&0.497\\0.499&1.988\end{smallmatrix}\bigr)$ & $\bigl(\begin{smallmatrix}0.999&0.495\\0.499&1.980\end{smallmatrix}\bigr)$ & $\bigl(\begin{smallmatrix}0.999&0.492\\0.499&1.966\end{smallmatrix}\bigr)$ & $\bigl(\begin{smallmatrix}0.997&0.486\\0.499&1.943\end{smallmatrix}\bigr)$ & $\bigl(\begin{smallmatrix}0.996&0.476\\0.498&1.904\end{smallmatrix}\bigr)$ & $\bigl(\begin{smallmatrix}0.993&0.460\\0.496&1.840\end{smallmatrix}\bigr)$ \\[3pt]
 & $\overline{\hat{\bH}}_{3}\,(=\hat{\bQ}^{-1})$ & $\bigl(\begin{smallmatrix}1.000&0.501\\0.500&2.004\end{smallmatrix}\bigr)$ & $\bigl(\begin{smallmatrix}1.000&0.502\\0.500&2.007\end{smallmatrix}\bigr)$ & $\bigl(\begin{smallmatrix}1.000&0.503\\0.500&2.012\end{smallmatrix}\bigr)$ & $\bigl(\begin{smallmatrix}1.001&0.505\\0.500&2.020\end{smallmatrix}\bigr)$ & $\bigl(\begin{smallmatrix}1.001&0.508\\0.501&2.034\end{smallmatrix}\bigr)$ & $\bigl(\begin{smallmatrix}1.002&0.515\\0.501&2.058\end{smallmatrix}\bigr)$ \\[3pt]
 & $\overline{\hat{\bR}}_{\mathrm{Fan},F}$ & $\bigl(\begin{smallmatrix}1.000&0.500\\0.500&2.000\end{smallmatrix}\bigr)$ & $\bigl(\begin{smallmatrix}1.000&0.500\\0.500&2.000\end{smallmatrix}\bigr)$ & $\bigl(\begin{smallmatrix}1.000&0.500\\0.500&2.000\end{smallmatrix}\bigr)$ & $\bigl(\begin{smallmatrix}1.000&0.500\\0.500&2.000\end{smallmatrix}\bigr)$ & $\bigl(\begin{smallmatrix}1.000&0.500\\0.500&2.000\end{smallmatrix}\bigr)$ & $\bigl(\begin{smallmatrix}1.000&0.500\\0.500&2.000\end{smallmatrix}\bigr)$ \\[3pt]
\hline
500 & $\overline{\hat{\bH}}$ & $\bigl(\begin{smallmatrix}0.999&0.498\\0.500&1.992\end{smallmatrix}\bigr)$ & $\bigl(\begin{smallmatrix}0.999&0.496\\0.499&1.985\end{smallmatrix}\bigr)$ & $\bigl(\begin{smallmatrix}0.999&0.493\\0.500&1.973\end{smallmatrix}\bigr)$ & $\bigl(\begin{smallmatrix}0.998&0.487\\0.499&1.950\end{smallmatrix}\bigr)$ & $\bigl(\begin{smallmatrix}0.997&0.477\\0.498&1.907\end{smallmatrix}\bigr)$ & $\bigl(\begin{smallmatrix}0.994&0.458\\0.497&1.831\end{smallmatrix}\bigr)$ \\[3pt]
 & $\overline{\hat{\bH}}_{4}$ & $\bigl(\begin{smallmatrix}0.999&0.499\\0.500&1.995\end{smallmatrix}\bigr)$ & $\bigl(\begin{smallmatrix}0.999&0.498\\0.500&1.991\end{smallmatrix}\bigr)$ & $\bigl(\begin{smallmatrix}0.999&0.496\\0.500&1.984\end{smallmatrix}\bigr)$ & $\bigl(\begin{smallmatrix}0.999&0.492\\0.499&1.970\end{smallmatrix}\bigr)$ & $\bigl(\begin{smallmatrix}0.998&0.486\\0.499&1.944\end{smallmatrix}\bigr)$ & $\bigl(\begin{smallmatrix}0.996&0.474\\0.498&1.896\end{smallmatrix}\bigr)$ \\[3pt]
 & $\overline{\hat{\bH}}_{3}\,(=\hat{\bQ}^{-1})$ & $\bigl(\begin{smallmatrix}1.000&0.500\\0.500&2.002\end{smallmatrix}\bigr)$ & $\bigl(\begin{smallmatrix}1.000&0.501\\0.500&2.003\end{smallmatrix}\bigr)$ & $\bigl(\begin{smallmatrix}1.000&0.501\\0.500&2.006\end{smallmatrix}\bigr)$ & $\bigl(\begin{smallmatrix}1.000&0.503\\0.500&2.010\end{smallmatrix}\bigr)$ & $\bigl(\begin{smallmatrix}1.001&0.505\\0.500&2.019\end{smallmatrix}\bigr)$ & $\bigl(\begin{smallmatrix}1.001&0.509\\0.501&2.037\end{smallmatrix}\bigr)$ \\[3pt]
 & $\overline{\hat{\bR}}_{\mathrm{Fan},F}$ & $\bigl(\begin{smallmatrix}1.000&0.500\\0.500&2.000\end{smallmatrix}\bigr)$ & $\bigl(\begin{smallmatrix}1.000&0.500\\0.500&2.000\end{smallmatrix}\bigr)$ & $\bigl(\begin{smallmatrix}1.000&0.500\\0.500&2.000\end{smallmatrix}\bigr)$ & $\bigl(\begin{smallmatrix}1.000&0.500\\0.500&2.000\end{smallmatrix}\bigr)$ & $\bigl(\begin{smallmatrix}1.000&0.500\\0.500&2.000\end{smallmatrix}\bigr)$ & $\bigl(\begin{smallmatrix}1.000&0.500\\0.500&2.000\end{smallmatrix}\bigr)$ \\[3pt]
\hline\hline
\end{tabular}}
\end{table}

\begin{table}[htbp!]
\centering\scriptsize
\setlength{\tabcolsep}{3pt}
\caption{Monte Carlo standard deviations of the entries of the data-dependent rotation matrices (non-sparse loadings, 50{,}000 replications), in units of $10^{-2}$, matching the layout of Table~\ref{tab:D_rotations}. Each cell is the across-replication standard deviation of the corresponding entry of $\hat{\bH}$, $\hat{\bH}_4$, $\hat{\bH}_3=\hat{\bQ}^{-1}$, or $\hat{\bR}_{\mathrm{Fan},F}$. All standard deviations are below $0.15$ and shrink with the sample size. The overall dispersion ordering is $\hat{\bR}_{\mathrm{Fan},F}<\hat{\bH}_3<\hat{\bH}_4<\hat{\bH}$ in every design: $\hat{\bR}_{\mathrm{Fan},F}$ has both the smallest summed and the smallest largest-entry standard deviation throughout (the lone exception is the $(2,1)$ entry, where the more concentrated $\hat{\bH}$ and $\hat{\bH}_4$ are tighter).}
\label{tab:D_rotations_sd}
\resizebox{\textwidth}{!}{%
\begin{tabular}{clcccccc}
\hline\hline
$T=N$ & $(\alpha_1,\alpha_2)$: & (1.0,1.0) & (1.0,0.9) & (1.0,0.8) & (0.9,0.7) & (0.8,0.6) & (0.7,0.5) \\
\hline
50 & $\overline{\mathrm{sd}}(\hat{\bH})$ & $\bigl(\begin{smallmatrix}1.8&5.4\\1.1&7.1\end{smallmatrix}\bigr)$ & $\bigl(\begin{smallmatrix}1.8&6.9\\1.0&8.6\end{smallmatrix}\bigr)$ & $\bigl(\begin{smallmatrix}1.8&9.1\\0.9&10.3\end{smallmatrix}\bigr)$ & $\bigl(\begin{smallmatrix}2.1&10.6\\1.1&11.8\end{smallmatrix}\bigr)$ & $\bigl(\begin{smallmatrix}2.6&12.3\\1.3&13.0\end{smallmatrix}\bigr)$ & $\bigl(\begin{smallmatrix}3.1&14.1\\1.6&14.1\end{smallmatrix}\bigr)$ \\[3pt]
 & $\overline{\mathrm{sd}}(\hat{\bH}_{4})$ & $\bigl(\begin{smallmatrix}0.9&2.7\\1.4&3.6\end{smallmatrix}\bigr)$ & $\bigl(\begin{smallmatrix}0.9&3.0\\1.1&4.3\end{smallmatrix}\bigr)$ & $\bigl(\begin{smallmatrix}0.9&3.3\\0.8&5.0\end{smallmatrix}\bigr)$ & $\bigl(\begin{smallmatrix}1.1&3.9\\1.0&5.8\end{smallmatrix}\bigr)$ & $\bigl(\begin{smallmatrix}1.3&4.6\\1.2&6.6\end{smallmatrix}\bigr)$ & $\bigl(\begin{smallmatrix}1.6&5.4\\1.5&7.3\end{smallmatrix}\bigr)$ \\[3pt]
 & $\overline{\mathrm{sd}}(\hat{\bH}_{3})\,(=\hat{\bQ}^{-1})$ & $\bigl(\begin{smallmatrix}0.7&1.3\\2.7&0.7\end{smallmatrix}\bigr)$ & $\bigl(\begin{smallmatrix}0.6&1.2\\2.4&0.8\end{smallmatrix}\bigr)$ & $\bigl(\begin{smallmatrix}0.5&1.1\\2.2&0.9\end{smallmatrix}\bigr)$ & $\bigl(\begin{smallmatrix}0.7&1.4\\2.6&1.3\end{smallmatrix}\bigr)$ & $\bigl(\begin{smallmatrix}0.8&1.7\\3.3&1.9\end{smallmatrix}\bigr)$ & $\bigl(\begin{smallmatrix}1.0&2.2\\4.1&2.9\end{smallmatrix}\bigr)$ \\[3pt]
 & $\overline{\mathrm{sd}}(\hat{\bR}_{\mathrm{Fan},F})$ & $\bigl(\begin{smallmatrix}0.7&1.3\\2.7&0.7\end{smallmatrix}\bigr)$ & $\bigl(\begin{smallmatrix}0.6&1.2\\2.3&0.6\end{smallmatrix}\bigr)$ & $\bigl(\begin{smallmatrix}0.5&1.1\\2.1&0.5\end{smallmatrix}\bigr)$ & $\bigl(\begin{smallmatrix}0.6&1.3\\2.6&0.7\end{smallmatrix}\bigr)$ & $\bigl(\begin{smallmatrix}0.8&1.6\\3.2&0.8\end{smallmatrix}\bigr)$ & $\bigl(\begin{smallmatrix}1.0&2.0\\3.9&1.0\end{smallmatrix}\bigr)$ \\[3pt]
\hline
100 & $\overline{\mathrm{sd}}(\hat{\bH})$ & $\bigl(\begin{smallmatrix}0.9&2.7\\0.6&3.7\end{smallmatrix}\bigr)$ & $\bigl(\begin{smallmatrix}0.9&3.7\\0.5&4.6\end{smallmatrix}\bigr)$ & $\bigl(\begin{smallmatrix}0.9&5.2\\0.5&5.8\end{smallmatrix}\bigr)$ & $\bigl(\begin{smallmatrix}1.1&6.4\\0.6&7.0\end{smallmatrix}\bigr)$ & $\bigl(\begin{smallmatrix}1.4&7.8\\0.7&8.3\end{smallmatrix}\bigr)$ & $\bigl(\begin{smallmatrix}1.7&9.2\\0.9&9.4\end{smallmatrix}\bigr)$ \\[3pt]
 & $\overline{\mathrm{sd}}(\hat{\bH}_{4})$ & $\bigl(\begin{smallmatrix}0.5&1.4\\0.7&1.9\end{smallmatrix}\bigr)$ & $\bigl(\begin{smallmatrix}0.5&1.5\\0.5&2.3\end{smallmatrix}\bigr)$ & $\bigl(\begin{smallmatrix}0.5&1.8\\0.4&2.8\end{smallmatrix}\bigr)$ & $\bigl(\begin{smallmatrix}0.6&2.2\\0.5&3.4\end{smallmatrix}\bigr)$ & $\bigl(\begin{smallmatrix}0.7&2.7\\0.6&4.0\end{smallmatrix}\bigr)$ & $\bigl(\begin{smallmatrix}0.9&3.2\\0.8&4.7\end{smallmatrix}\bigr)$ \\[3pt]
 & $\overline{\mathrm{sd}}(\hat{\bH}_{3})\,(=\hat{\bQ}^{-1})$ & $\bigl(\begin{smallmatrix}0.3&0.7\\1.3&0.4\end{smallmatrix}\bigr)$ & $\bigl(\begin{smallmatrix}0.3&0.6\\1.2&0.3\end{smallmatrix}\bigr)$ & $\bigl(\begin{smallmatrix}0.3&0.5\\1.1&0.4\end{smallmatrix}\bigr)$ & $\bigl(\begin{smallmatrix}0.3&0.7\\1.3&0.6\end{smallmatrix}\bigr)$ & $\bigl(\begin{smallmatrix}0.4&0.9\\1.7&0.9\end{smallmatrix}\bigr)$ & $\bigl(\begin{smallmatrix}0.5&1.1\\2.2&1.4\end{smallmatrix}\bigr)$ \\[3pt]
 & $\overline{\mathrm{sd}}(\hat{\bR}_{\mathrm{Fan},F})$ & $\bigl(\begin{smallmatrix}0.3&0.7\\1.3&0.3\end{smallmatrix}\bigr)$ & $\bigl(\begin{smallmatrix}0.3&0.6\\1.1&0.3\end{smallmatrix}\bigr)$ & $\bigl(\begin{smallmatrix}0.3&0.5\\1.0&0.3\end{smallmatrix}\bigr)$ & $\bigl(\begin{smallmatrix}0.3&0.7\\1.3&0.3\end{smallmatrix}\bigr)$ & $\bigl(\begin{smallmatrix}0.4&0.8\\1.7&0.4\end{smallmatrix}\bigr)$ & $\bigl(\begin{smallmatrix}0.5&1.1\\2.1&0.5\end{smallmatrix}\bigr)$ \\[3pt]
\hline
200 & $\overline{\mathrm{sd}}(\hat{\bH})$ & $\bigl(\begin{smallmatrix}0.4&1.4\\0.3&1.9\end{smallmatrix}\bigr)$ & $\bigl(\begin{smallmatrix}0.4&2.0\\0.2&2.4\end{smallmatrix}\bigr)$ & $\bigl(\begin{smallmatrix}0.4&3.0\\0.2&3.2\end{smallmatrix}\bigr)$ & $\bigl(\begin{smallmatrix}0.6&3.8\\0.3&4.1\end{smallmatrix}\bigr)$ & $\bigl(\begin{smallmatrix}0.8&4.8\\0.4&5.1\end{smallmatrix}\bigr)$ & $\bigl(\begin{smallmatrix}1.0&6.0\\0.5&6.1\end{smallmatrix}\bigr)$ \\[3pt]
 & $\overline{\mathrm{sd}}(\hat{\bH}_{4})$ & $\bigl(\begin{smallmatrix}0.2&0.7\\0.4&0.9\end{smallmatrix}\bigr)$ & $\bigl(\begin{smallmatrix}0.2&0.8\\0.2&1.2\end{smallmatrix}\bigr)$ & $\bigl(\begin{smallmatrix}0.2&0.9\\0.2&1.5\end{smallmatrix}\bigr)$ & $\bigl(\begin{smallmatrix}0.3&1.2\\0.2&1.9\end{smallmatrix}\bigr)$ & $\bigl(\begin{smallmatrix}0.4&1.5\\0.3&2.4\end{smallmatrix}\bigr)$ & $\bigl(\begin{smallmatrix}0.5&1.9\\0.4&2.9\end{smallmatrix}\bigr)$ \\[3pt]
 & $\overline{\mathrm{sd}}(\hat{\bH}_{3})\,(=\hat{\bQ}^{-1})$ & $\bigl(\begin{smallmatrix}0.2&0.3\\0.7&0.2\end{smallmatrix}\bigr)$ & $\bigl(\begin{smallmatrix}0.1&0.3\\0.6&0.2\end{smallmatrix}\bigr)$ & $\bigl(\begin{smallmatrix}0.1&0.3\\0.5&0.2\end{smallmatrix}\bigr)$ & $\bigl(\begin{smallmatrix}0.2&0.3\\0.7&0.3\end{smallmatrix}\bigr)$ & $\bigl(\begin{smallmatrix}0.2&0.5\\0.9&0.4\end{smallmatrix}\bigr)$ & $\bigl(\begin{smallmatrix}0.3&0.6\\1.2&0.7\end{smallmatrix}\bigr)$ \\[3pt]
 & $\overline{\mathrm{sd}}(\hat{\bR}_{\mathrm{Fan},F})$ & $\bigl(\begin{smallmatrix}0.2&0.3\\0.7&0.2\end{smallmatrix}\bigr)$ & $\bigl(\begin{smallmatrix}0.1&0.3\\0.6&0.1\end{smallmatrix}\bigr)$ & $\bigl(\begin{smallmatrix}0.1&0.3\\0.5&0.1\end{smallmatrix}\bigr)$ & $\bigl(\begin{smallmatrix}0.2&0.3\\0.7&0.2\end{smallmatrix}\bigr)$ & $\bigl(\begin{smallmatrix}0.2&0.4\\0.9&0.2\end{smallmatrix}\bigr)$ & $\bigl(\begin{smallmatrix}0.3&0.6\\1.1&0.3\end{smallmatrix}\bigr)$ \\[3pt]
\hline
500 & $\overline{\mathrm{sd}}(\hat{\bH})$ & $\bigl(\begin{smallmatrix}0.2&0.6\\0.1&0.8\end{smallmatrix}\bigr)$ & $\bigl(\begin{smallmatrix}0.2&0.9\\0.1&1.0\end{smallmatrix}\bigr)$ & $\bigl(\begin{smallmatrix}0.2&1.4\\0.1&1.5\end{smallmatrix}\bigr)$ & $\bigl(\begin{smallmatrix}0.2&1.9\\0.1&2.0\end{smallmatrix}\bigr)$ & $\bigl(\begin{smallmatrix}0.3&2.6\\0.2&2.6\end{smallmatrix}\bigr)$ & $\bigl(\begin{smallmatrix}0.4&3.4\\0.2&3.3\end{smallmatrix}\bigr)$ \\[3pt]
 & $\overline{\mathrm{sd}}(\hat{\bH}_{4})$ & $\bigl(\begin{smallmatrix}0.1&0.3\\0.1&0.4\end{smallmatrix}\bigr)$ & $\bigl(\begin{smallmatrix}0.1&0.3\\0.1&0.5\end{smallmatrix}\bigr)$ & $\bigl(\begin{smallmatrix}0.1&0.4\\0.1&0.7\end{smallmatrix}\bigr)$ & $\bigl(\begin{smallmatrix}0.1&0.6\\0.1&0.9\end{smallmatrix}\bigr)$ & $\bigl(\begin{smallmatrix}0.2&0.7\\0.1&1.2\end{smallmatrix}\bigr)$ & $\bigl(\begin{smallmatrix}0.2&1.0\\0.2&1.6\end{smallmatrix}\bigr)$ \\[3pt]
 & $\overline{\mathrm{sd}}(\hat{\bH}_{3})\,(=\hat{\bQ}^{-1})$ & $\bigl(\begin{smallmatrix}0.1&0.1\\0.3&0.1\end{smallmatrix}\bigr)$ & $\bigl(\begin{smallmatrix}0.1&0.1\\0.2&0.1\end{smallmatrix}\bigr)$ & $\bigl(\begin{smallmatrix}0.0&0.1\\0.2&0.1\end{smallmatrix}\bigr)$ & $\bigl(\begin{smallmatrix}0.1&0.1\\0.3&0.1\end{smallmatrix}\bigr)$ & $\bigl(\begin{smallmatrix}0.1&0.2\\0.4&0.2\end{smallmatrix}\bigr)$ & $\bigl(\begin{smallmatrix}0.1&0.3\\0.5&0.3\end{smallmatrix}\bigr)$ \\[3pt]
 & $\overline{\mathrm{sd}}(\hat{\bR}_{\mathrm{Fan},F})$ & $\bigl(\begin{smallmatrix}0.1&0.1\\0.3&0.1\end{smallmatrix}\bigr)$ & $\bigl(\begin{smallmatrix}0.1&0.1\\0.2&0.1\end{smallmatrix}\bigr)$ & $\bigl(\begin{smallmatrix}0.0&0.1\\0.2&0.0\end{smallmatrix}\bigr)$ & $\bigl(\begin{smallmatrix}0.1&0.1\\0.3&0.1\end{smallmatrix}\bigr)$ & $\bigl(\begin{smallmatrix}0.1&0.2\\0.4&0.1\end{smallmatrix}\bigr)$ & $\bigl(\begin{smallmatrix}0.1&0.3\\0.5&0.1\end{smallmatrix}\bigr)$ \\[3pt]
\hline\hline
\end{tabular}}
\end{table}

\section{Additional Empirical Application}
\label{subsec:addapplication}
\setcounter{figure}{0}
\setcounter{table}{0}
\renewcommand{\thefigure}{E.\arabic{figure}}
\renewcommand{\thetable}{E.\arabic{table}}
In this section, we consider an additional empirical example that is justified only by our results. The number of factors $r$ should be estimated while allowing $\lambda_{k}\asymp N^{\alpha_{k}}$, $\alpha_{k}\in(0,1]$ for $k=1,\dots,r$. The methods of \cite{BaiNg2002} and \cite{AhnHorenstein2013} are not asymptotically valid in general, but \cite{Onatski2010} edge distribution (ED) method based on a calibration is, and is used below.

\subsection{Uncovering the characteristics of factors by multiple testing}
\label{subsubsec:application1}
The true factor, $\mathbf{F}^{\ast}$, is not identifiable in general,
however, $\mathbf{F}^{0}=\mathbf{F}^{\ast}\mathbf{H}$ is. It is therefore of interest to researchers to study the characteristics of each column of
$\mathbf{F}^{0}=(\mathbf{f}_{1}^{0},\dots ,\mathbf{f}_{r}^{0})$. To find out the property of $k^{th}$ factor $\mathbf{f}_{k}^{0}$, following \cite{UY2019inference}, we
perform multiple testing of $H_{0}^{i,k}:b_{i,k}^{0}=0$ versus $H_{1}%
^{i,k}:b_{i,k}^{0}\neq0$, $i=1,\dots ,N$, $k=1,\dots ,r$, to find out the
characteristics of the cross-sectional units $i$ that are significantly
affected by the $k^{th}$ common factor $\mathbf{f}_{k}^{0}$. Given that the
number of hypotheses tends to infinity as $N\rightarrow\infty$, \cite{UY2019inference}
propose a \cite{JJ2019} type algorithm to control the false discovery rate (FDR) below a predetermined level $q\in(0,1)$. In short, construct the t-ratio with appropriate standard error\ such that T$_{ik}=\hat{b}_{i,k}/se_{i,k}%
\overset{d}{\rightarrow}N(0,1)$ under the null, and obtain the threshold
$t_{0}$ such that $\FDR\leq q+o(1)$ as $N,T\rightarrow\infty$, where $\FDR=\E[| \hat{\mathcal{S}}\cap \mathcal{S}^{c}| / |\hat{\mathcal{S}}|]$ with $\mathcal{S}=\{(i,k):b_{ik}^{0}\neq0\}$ and
$\hat{\mathcal{S}}=\{(i,k):\left\vert \text{T}_{ik}\right\vert >t_{0}\}$; see \cite{UY2019inference} for more details.

We consider the monthly excess returns of the firms that constitute the
Standard \& Poor's 500 stock index (S\&P500) in April 2018 over the period
from May 1998 to April 2018, leaving 376 series. The monthly excess returns of
the $i^{th}$ security for month $t$ are computed as $100(P_{ti}-P_{t-1,i}%
)/P_{t-1,i}+DY_{ti}/12-r_{ft}$, where $P_{ti}$ is the month-end price,
$DY_{ti}$ is the annual dividend yield in percent, and $r_{ft}$ is the
one-month US treasury bill rate. $P_{ti}$ and $DY_{ti}$ are obtained from
Datastream, and $r_{ft}$ is obtained from Ken French's data library website.
We standardize the excess returns before the analysis. The ED estimation with
$r_{\max}=9$ gave $\hat{r}=3$.

\begin{figure}[h!]
	\centering
	\begin{minipage}{1\hsize}
		\centering
		\includegraphics[width=0.8\linewidth]{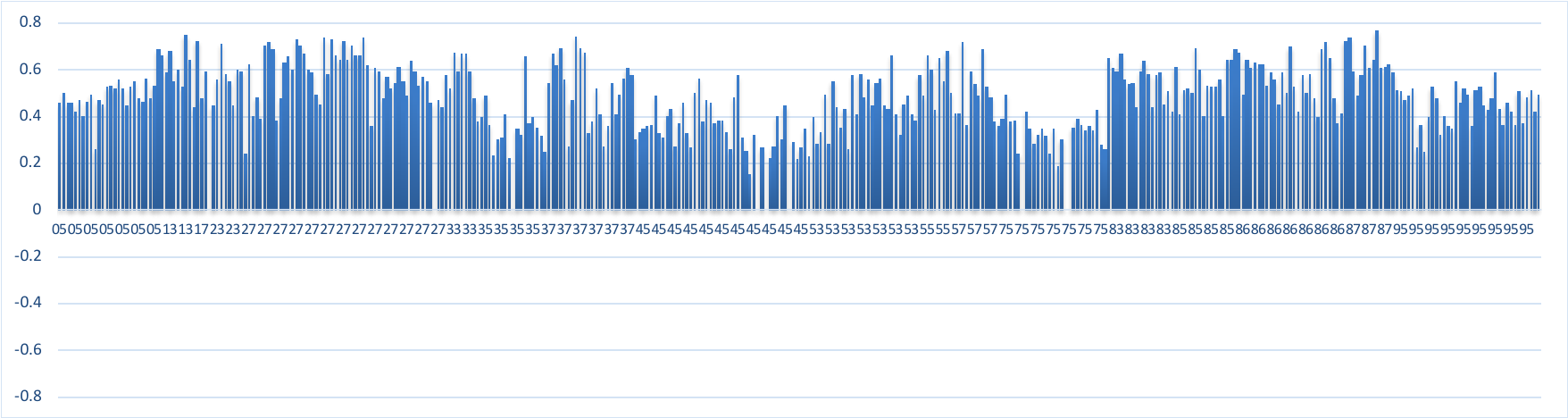}
		\subcaption{\footnotesize{$\hat{\bb}_1$ replacing insignificant values with zeros}}
	\end{minipage}
	\begin{minipage}{1\hsize}
		\centering

  \includegraphics[width=0.8\linewidth]{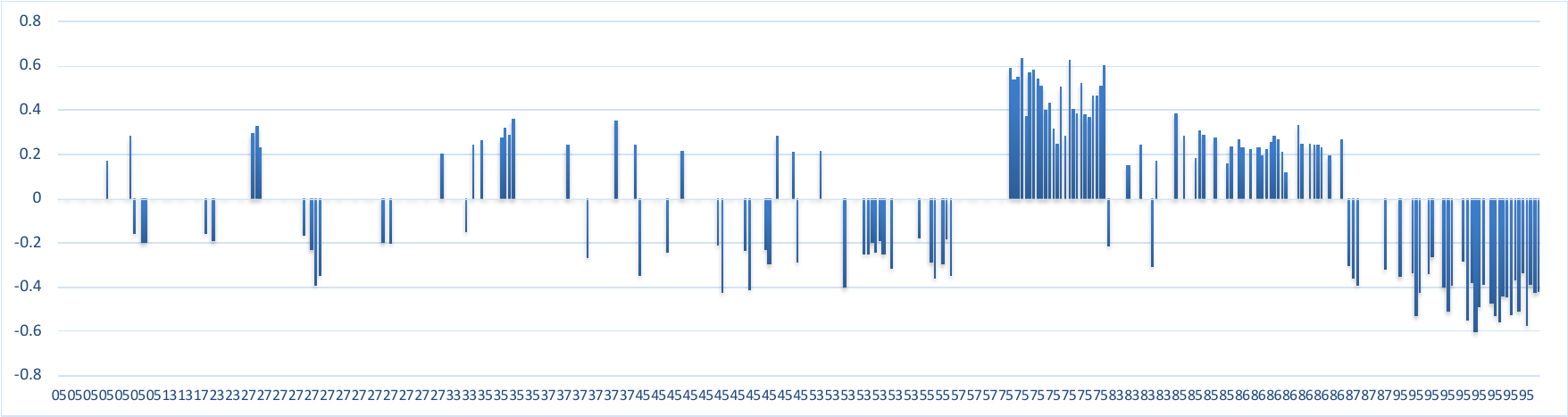}
\subcaption{\footnotesize{$\hat{\bb}_2$ replacing insignificant values with zeros}}
	\end{minipage}
	\begin{minipage}{1\hsize}
		\centering
		\includegraphics[width=0.8\linewidth]{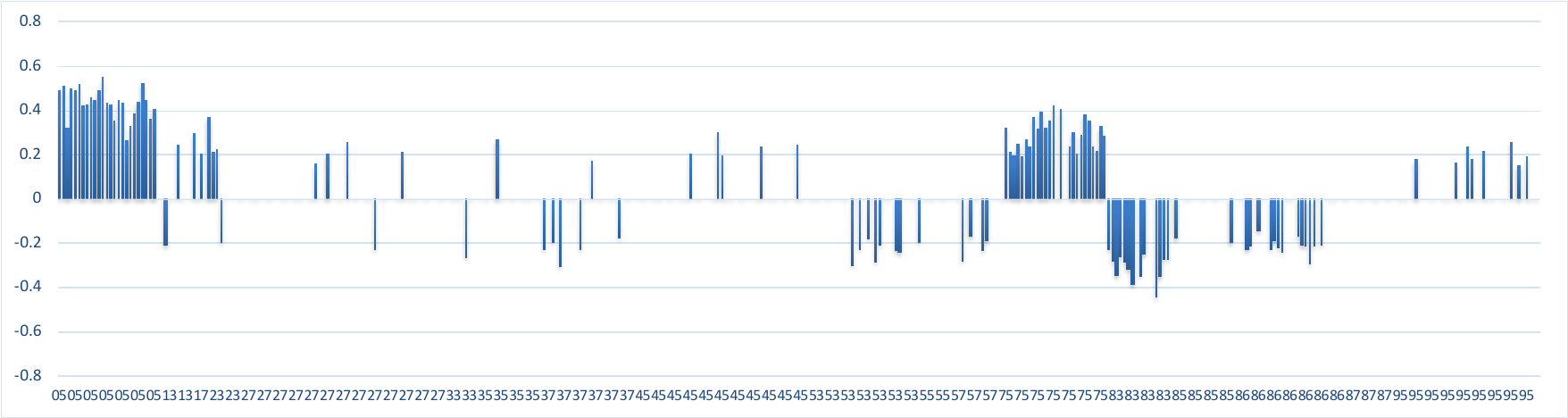}
\subcaption{\footnotesize{$\hat{\bb}_3$ replacing insignificant values with zeros}}
	\end{minipage} 
	\caption{Bar charts of S\&P500 $\hat{\bB}$ with $\hat{r}=3$, replacing insignificant values with zeros}
	\label{fig:sp500}
\end{figure}

We obtained the PC estimator, then computed the t-ratio for each $\hat{b}%
_{ik}$ using the heteroskedasticity and autocorrelation consistent (HAC)
standard error of Newey-West with the threshold value, integer part of
$T^{1/4}$. We emphasize that hypothesis testing such as $H_{0}:b_{ik}^{0}=0$
is only justified with our approach. With these t-ratios, we conducted the FDR controlled multiple test with level $q=0.01$. The computed threshold was $t_{0}=2.77$. Figure \ref{fig:sp500} shows the values of the loadings of each of the three PC factors, where the values of insignificant loadings are replaced by zeros. The firm securities are ordered and categorized into 17 industry sectors according to the Industry Classification Benchmark (ICB). Define an estimate of the divergence rate of the $k^{th}$ factor under the choice of $q$ as $\hat{\alpha}_{k}=\ln(\sum_{i=1}^{N}I[\left\vert \text{T}_{ik}\right\vert >t_{0}])/\ln(N)$. We found that $(\hat{\alpha}_{1},\hat{\alpha}_{2},\hat{\alpha}_{3})=(1.00,0.84,0.81)$, which can be taken as evidence for weak factor models, hence, supports our approach. The first factor can be considered as the market factor, as virtually all the securities are significantly affected by it. The characteristics of the second and the third factors are visualized in the figure. The significant loadings of the second factor are concentrated on the ICB two digits-industry sectors 75, 85, 86, and 95, hence we may call it the Utility, Insurance, Real Estate and Technology factor. The significant loadings of the third factor are concentrated on 05, 75, 83 and 86, so we may call it the Oil\&Gas, Utility, Banks and Real Estate factor. The overlap of 75 and 86, Utility and Real Estate, for these two factors indicates the significant presence of these industries in the stock market during the data period.

\section{Coverage of confidence and prediction intervals}
\label{app:forecast-interval-mc}

\setcounter{table}{0}
\setcounter{figure}{0}
\renewcommand{\thefigure}{F.\arabic{figure}}
\renewcommand{\thetable}{F.\arabic{table}}
We examine the empirical coverage of the intervals for the conditional mean
$y_{T+h|T}=\mathbf f_T^{0\prime}\bgamma^0+\mathbf w_T'\boldsymbol\beta$
and for the actual value $y_{T+h}$. The feasible forecast is
$\hat y_{T+h|T}=\hat{\mathbf f}_T'\hat{\bgamma}
+\mathbf w_T'\hat{\boldsymbol\beta}$. Denoting the corresponding quantities
in the $s$-th Monte Carlo replication ($s\in[R]$) by $y_{T+h|T}^{[s]}$,
$y_{T+h}^{[s]}$, and $\hat y_{T+h|T}^{[s]}$, respectively, we compute
empirical coverage as $
CI_{\tau}(y)
=
R^{-1}\sum_{s=1}^{R}
1\{y^{[s]}\in[LB_{\tau,y}^{[s]},UB_{\tau,y}^{[s]}]\}$, 
for $y=y_{T+h|T},\,y_{T+h}$. 
For the conditional mean, 
$[LB_{\tau,y_{T+h|T}}^{[s]},UB_{\tau,y_{T+h|T}}^{[s]}]
=
\hat y_{T+h|T}^{[s]}
\pm c_{\tau}\sigma_{T+h|T}$, 
whereas for the actual value, 
$[LB_{\tau,y_{T+h}}^{[s]},UB_{\tau,y_{T+h}}^{[s]}]
=
\hat y_{T+h|T}^{[s]}
\pm c_{\tau}
(\sigma_{T+h|T}^{2}+\sigma_{\epsilon}^{2})^{1/2}$. 
Here $c_{\tau}$ is defined by
$\Pr(|Z|\leq c_{\tau})=\tau$ for $Z\sim N(0,1)$, and 
$\sigma_{T+h|T}^{2}
=
\sigma_{\epsilon}^{2}
\mathbf z_T^{0\prime}
(\mathbf Z^{0\prime}\mathbf Z^{0})^{-1}
\mathbf z_T^{0}
+
\sigma_e^{2}
\bgamma^{0\prime}
(\mathbf B^{0\prime}\mathbf B^{0})^{-1}
\bgamma^{0}$, 
where
$\mathbf Z^{0}=(\mathbf F^{0},\mathbf W)
=(\mathbf z_1^{0},\ldots,\mathbf z_T^{0})'$. We consider
$\tau=0.90$ and $0.95$.
\begin{figure}[h!]	
	\centering
	\includegraphics[clip,width=14.3cm]{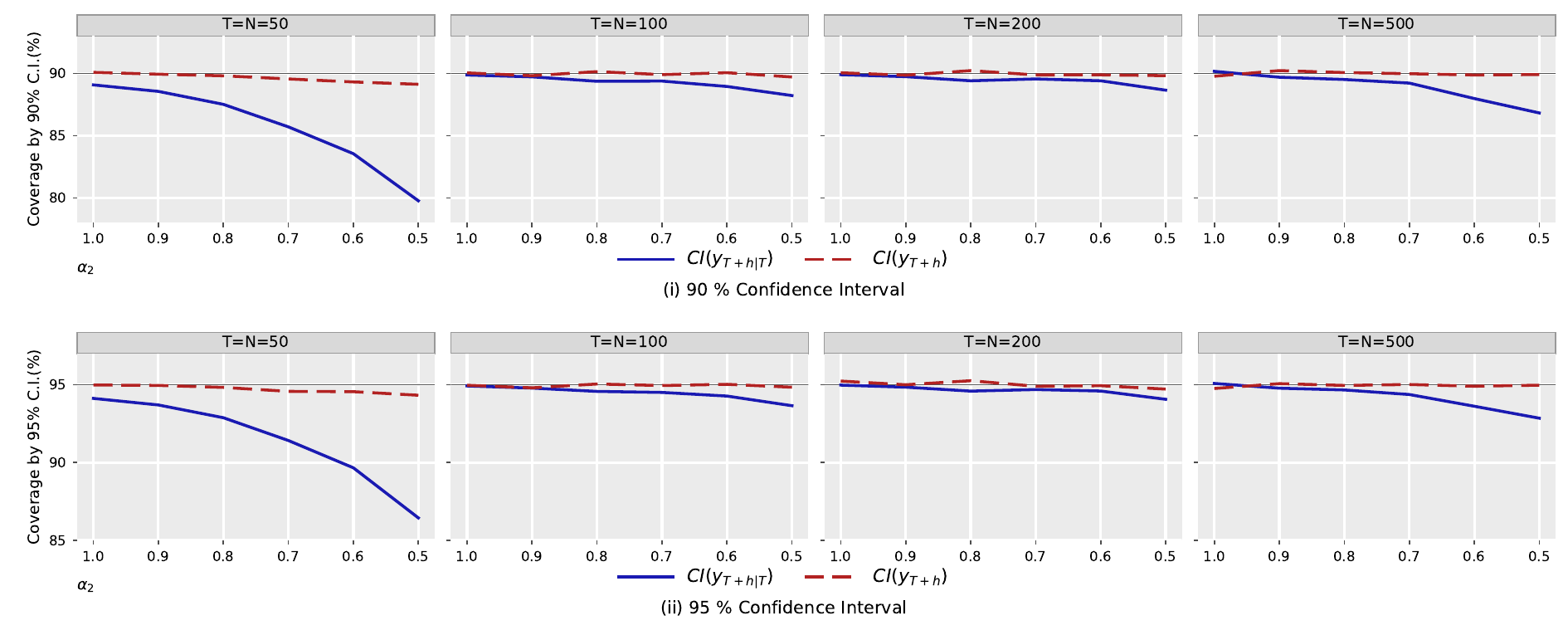}
	\vspace{-5mm}
\caption{Empirical coverages of the 90\% and 95\% intervals for the conditional mean $y_{T+h|T}$ and the actual value $y_{T+h}$.}
 \label{fig:7}
\end{figure}
Figure~\ref{fig:7} summarizes the empirical coverages. The coverage of the actual value, $CI_{\tau}(y_{T+h})$, is close to the nominal level in every model and at every sample size: 89.1\%--90.2\% for the 90\% intervals and 94.3\%--95.3\% for the 95\% intervals. By contrast, the coverage of the conditional mean, $CI_{\tau}(y_{T+h|T})$, is accurate in the stronger designs but deteriorates in the weakest designs when the sample is small. At $T=N=50$, the 90\% and 95\% intervals cover $y_{T+h|T}$ only 79.8\% and 86.4\% of the time, respectively, in the weakest design. This is consistent with Theorem~\ref{thm:forecast2}: when $N=T$, its rate conditions reduce to $\alpha_r>1/2$, so the weakest design, with $\alpha_2=0.5$, lies on the boundary of the admissible region, while the second weakest design, with $\alpha_2=0.6$, satisfies the condition only marginally. Accordingly, the coverage of $y_{T+h|T}$ is essentially nominal for all designs with $\alpha_2\geq0.7$ once $T=N\geq100$, whereas convergence is visibly slower in the two weakest designs, with mild undercoverage remaining even at $T=N=500$. The same pattern is obtained under the sparse loadings reported in Section~\ref{subsec:addMC} of the Supplementary Materials.


The contrast between the two intervals is straightforward. The width of
$CI_{\tau}(y_{T+h})$ is dominated by the innovation variance
$\sigma_{\epsilon}^{2}$, which is accurately estimated in all designs. Hence
the slower convergence of $\hat{\sigma}_{T+h|T}^{2}$ in the weakest designs has
little effect on coverage of the actual value, but matters for coverage of the
conditional mean. Overall, the results support the usefulness of the forecast
intervals in Theorem~\ref{thm:forecast2}, while showing that inference on the
conditional mean requires the factor-strength conditions to hold with some
margin.

\section{Joint normal approximation under alternative centerings}
\label{app:joint-normal-centerings}

\setcounter{table}{0}
\setcounter{figure}{0}
\renewcommand{\thetable}{G.\arabic{table}}
\renewcommand{\thefigure}{G.\arabic{figure}}

We first check the quality of the joint normal approximations of the PC
estimators, by comparing the following statistics with the $\chi_{r}^{2}$
distribution:
$
Q_f^{2}(\Delta\hat{\bff}_{t}) := 
\Delta\hat{\bff}_{t}^{\prime}(\bB^{0\prime}\bB^{0})^{\frac{1}{2}}
\mathbf{\Gamma}_{t}^{-1}
(\bB^{0\prime}\bB^{0})^{\frac{1}{2}}\Delta\hat{\bff}_{t},
~
Q_b^{2}(\Delta\hat{\bb}_{i}) := T(\Delta\hat{\bb}_{i})^{\prime}
\mathbf{\Phi}_{i}^{-1}(\Delta\hat{\bb}_{i})
$
with
$\Delta\hat{\bff}_{t}=\hat{\bff}_{t}-\bff_{t}^{0}$,
$\hat{\bff}_{t}-\hat{\bH}_{4}^{\prime}\bff_{t}^{*}$,
$\hat{\bff}_{t}-\hat{\bH}^{\prime}\bff_{t}^{*}$,
$\hat{\bff}_{t}-\hat{\bR}_{\mathrm{Fan},F}^{\prime}\bff_{t}^{*}$,
$\mathbf{\Gamma}_{t}=\sigma_{e}^{2}\bI_{r}$,
$\Delta\hat{\bb}_{i}=\hat{\bb}_{i}-\bb_{i}^{0}$,
$\hat{\bb}_{i}-\hat{\bQ}\bb_{i}^{*}$,
$\hat{\bb}_{i}-\hat{\bH}^{-1}\bb_{i}^{*}$,
$\hat{\bb}_{i}-\hat{\bR}_{\mathrm{Fan},B}'\bb_{i}^{*}$, and
$\mathbf{\Phi}_{i}=\sigma_{e}^{2}\bI_{r}$. We assess the accuracy of the
tail behavior by computing the frequencies of
$Q_\cdot^{2}>\chi_{r,0.95}^{2}$ over the replications, where
$\chi_{r,0.95}^{2}$ is the $95$-percentile of a $\chi_{r}^{2}$
distribution. The closer the size (frequency) is to 5\%, the more
accurate the approximation by the statistics.

\begin{figure}[h!]
	\centering
	\includegraphics[clip,width=14.3cm]{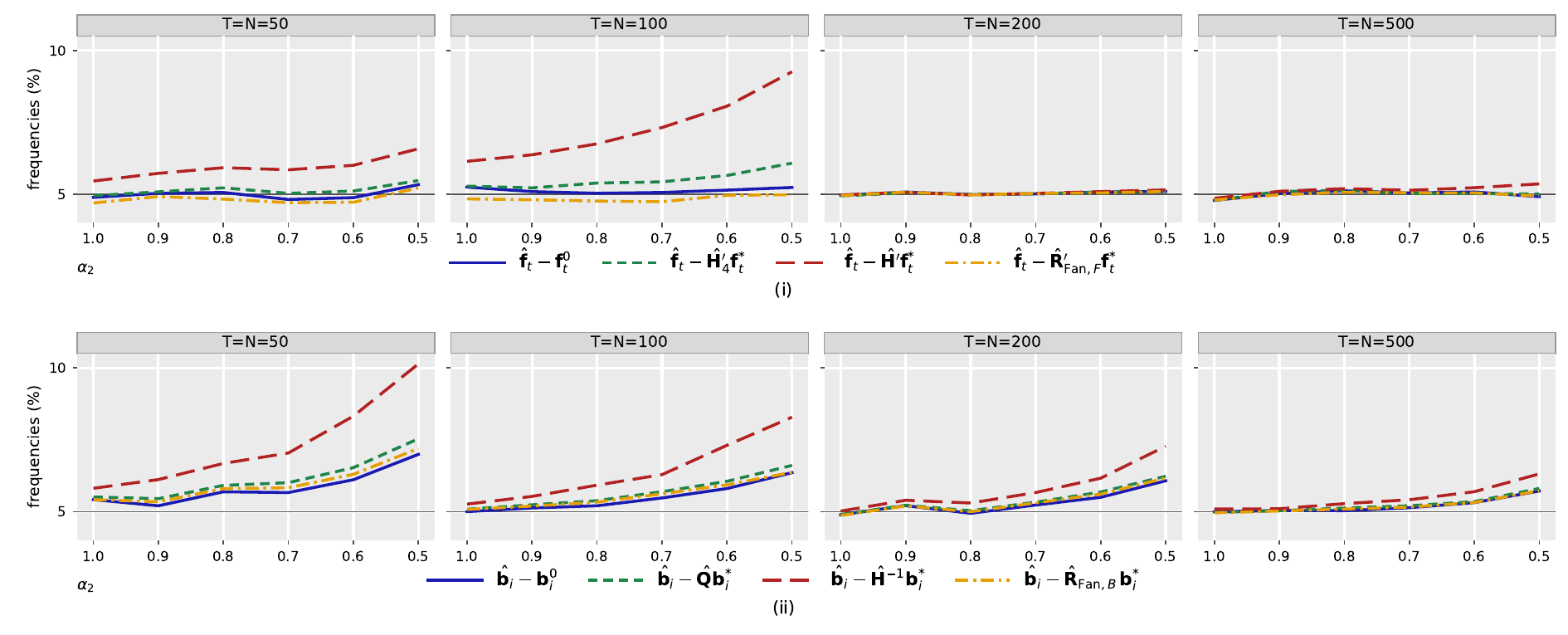}
	\vspace{-5mm}
\caption{Rejection frequencies based on the chi-square statistics
(i)$Q_f^{2}(\Delta\hat{\bff}_{t})$ and
(ii)$Q_b^{2}(\Delta\hat{\bb}_{i})$, using centerings at the pseudo-true
parameters and at the data-dependent rotations, including the Procrustes
rotation of \citet{FanEtAl2024}.}
 \label{fig:2}
\end{figure}

Figure \ref{fig:2} shows the frequencies. The PC estimator performs the
best in the cases when centered at $(\bff_t^0,\bb_i^0)$ and at the targets
$(\hat{\bR}_{\mathrm{Fan},F}^{\prime}\bff_t^{*},
\hat{\bR}_{\mathrm{Fan},B}'\bb_i^{*})$ of \citet{FanEtAl2024}, which behave almost identically, compared to those when centered at
$(\hat{\bH}_4'\bff_t^*,\hat{\bQ}\bb_i^*)$ and, especially,
$(\hat{\bH}'\bff_t^*,\hat{\bH}^{-1}\bb_i^*)$. The over-rejection of the statistics with $\hat{\bH}$ and $\hat{\bH}_4$ grows as the model
weakens, reflecting that these center at a sample-dependent rotation rather
than at a fixed parameter, whereas the pseudo-true and \citet{FanEtAl2024}
centerings remain close to nominal throughout.
The near-coincidence of the
pseudo-true and \citet{FanEtAl2024}'s centerings again reflects the fact that
$\hat{\bR}_{\mathrm{Fan},F}$ essentially equals $\bH$ in this setting. We stress, however, that only the pseudo-true centering corresponds to a fixed null hypothesis; the data-dependent centerings change with the sample and are reported as approximation diagnostics.

\section{Additional diagnostics for factor-augmented regressions}
\label{app:far-joint-normal}

\setcounter{table}{0}
\setcounter{figure}{0}
\renewcommand{\thetable}{H.\arabic{table}}
\renewcommand{\thefigure}{H.\arabic{figure}}

We assess the chi-square approximations for 
$Q_\gamma^{2}(\Delta\hat{\bgamma})
=
T(\Delta\hat{\bgamma})'
\boldsymbol{\Sigma}_{\gamma^{0}}^{-1}
(\Delta\hat{\bgamma})$,
where
$\Delta\hat{\bgamma}=\hat{\bgamma}^{0}-\bgamma^{0}$,
$\hat{\bgamma}-\bgamma^{0}$, or
$\hat{\bgamma}-\hat{\bH}^{-1}\bgamma^{*}$, and 
$\boldsymbol{\Sigma}_{\gamma^{0}}
=
\sigma_{\epsilon}^{2}
\left(T^{-1}\bF^{0\prime}\bM_w\bF^{0}\right)^{-1}$. 
We evaluate tail accuracy by computing the frequencies with which
$Q_\gamma^{2}(\Delta\hat{\bgamma})>\chi_{r,0.95}^{2}$ over the Monte Carlo
replications. Frequencies closer to 5\% indicate a more accurate chi-square
approximation.
\begin{figure}[htp!]	
	\centering
	\includegraphics[clip,width=14.3cm]{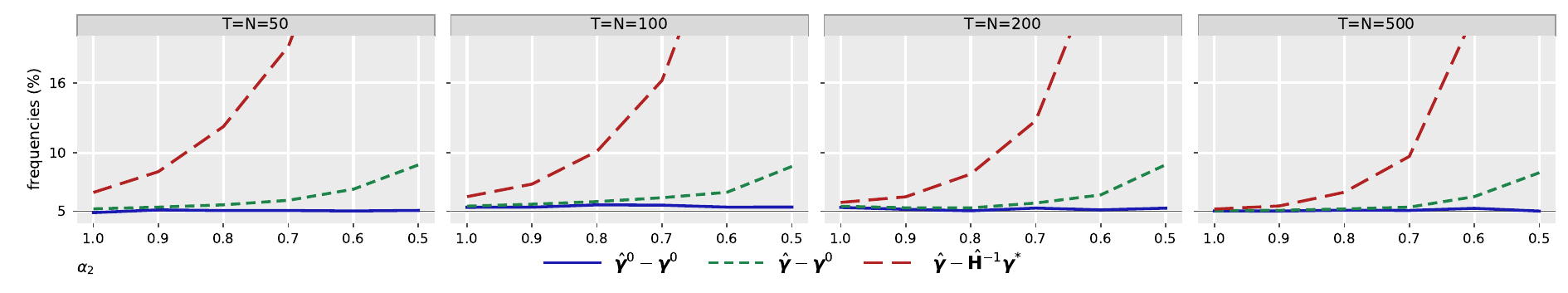}
	\vspace{-4mm}
\caption{Rejection frequencies based on
$Q_\gamma^{2}(\Delta\hat{\bgamma})>\chi_{r,0.95}^{2}$ under the three
coefficient centerings.}
 \label{fig:5}
\end{figure}
Figure \ref{fig:5} summarizes the results. The centering based on
$\hat{\bH}^{-1}\bgamma^{*}$ leads to severe over-rejection in the weakest
model, with rejection frequencies around 60\%. In contrast, the approximation
relative to the pseudo-true parameter is accurate except in the weakest
model.

\end{document}